\documentclass{birkmult}
\usepackage{amsmath,amssymb,amsfonts}
\usepackage{verbatim}

\newtheorem{theorem*}{Theorem}[section]{\bf}{\it}
\newtheorem{lemma*}[theorem*]{Lemma}{\bf}{\it}
\newtheorem{proposition*}[theorem*]{Proposition}{\bf}{\it}
\newtheorem{corollary*}[theorem*]{Corollary}{\bf}{\it}
\newtheorem{definition*}[theorem*]{Definition}{\bf}{\it}

\def\a{\mathfrak{t}}
\def\A{{\rm A}}
\def\ab{\bar{a}}
\def\al{\alpha}
\def\At{\widetilde A}

\def\b{\mathfrak{b}}
\def\B{{\rm B}}
\def\Bb{\ts\overline{\ns\rm B\!}\,}
\def\be{\beta}
\def\bet{\widetilde{\be}}
\def\bg{\tilde g}
\def\bi{\tilde\imath}
\def\bj{\tilde\jmath}
\def\bh{\tilde h}
\def\bk{\tilde k}
\def\bl{\tilde l}
\def\bp{\tilde p}
\def\bt{\ts\,\raise-0.5pt\hbox{\small$\boxtimes$}\,\,}
\def\Bt{\widetilde B}
\def\C{{\rm C}}
\def\cb{\bar{c}}
\def\CC{\mathbb{C}}
\def\com{\ts,\hskip-.5pt}
\def\Ct{\widetilde C}

\def\d{\partial}
\def\db{\bar{d}}
\def\D{\mathfrak{A}}
\def\de{\delta}
\def\De{\Delta}
\def\Dt{\widetilde D}

\def\e{\mathfrak{e}}
\def\E{\mathcal{E}}
\def\Eb{\widetilde{\mathcal{E}}}
\def\End{\operatorname{End}\ts}
\def\ep{\varepsilon}
\def\Ep{E^{\ts\prime}}
\def\Et{\widetilde E}

\def\f{\mathfrak{f}}
\def\F{\mathcal{F}}
\def\Fp{F^{\,\prime}}
\def\Ft{\widetilde{F}}
\def\Ftp{\widetilde{F}^{\,\prime}}

\def\g{\mathfrak{g}}
\def\ga{\gamma}
\def\ge{\geqslant}
\def\gl{\mathfrak{gl}}
\def\Gt{\widetilde G}

\def\h{\mathfrak{h}}
\def\H{\mathfrak{H}}
\def\Hh{\mathfrak{B}}
\def\Hom{\operatorname{Hom}}
\def\Hp{\H^{\ts\prime}}
\def\Ht{\widetilde H}

\def\I{{\rm I}}
\def\Ib{\tts\overline{\nns\rm I\nns}\tts}
\def\io{\iota}
\def\It{\widetilde I}

\def\J{{\rm J}}
\def\Jb{\,\overline{\!\rm J\ns}\ts}
\def\Jp{{\rm J}^{\ts\prime}}
\def\Jpb{\,\overline{\!\rm J\ns}\ts^{\ts\prime}}
\def\Jt{\widetilde J}

\def\ka{\kappa}

\def\la{\lambda}
\def\las{\la^{\ts\prime}}
\def\lat{\bar\la}
\def\lcd{\ts,\ldots,}
\def\le{\leqslant}

\def\mus{\mu^{\ts\prime}}
\def\mut{\bar\mu}

\def\n{\mathfrak{n}}
\def\nns{\hskip-.5pt}
\def\np{\n^{\ts\prime}}
\def\ns{\hskip-1pt}
\def\nus{\nu^{\ts\prime}}
\def\nut{\bar{\nu}}

\def\om{\omega}
\def\op{\oplus}
\def\ot{\otimes}

\def\p{\mathfrak{p}}
\def\P{\mathcal{P}}
\def\Pb{\bar P}
\def\PD{\mathcal{PD}}
\def\ph{\chi}
\def\Pp{P^{\ts\prime}}
\def\Pt{\widetilde P}

\def\q{\mathfrak{q}}
\def\qp{\q^{\ts\prime}}
\def\Qb{\bar Q}
\def\Qt{\widetilde Q}

\def\r{\mathfrak{r}}
\def\rp{\r^{\ts\prime}}
\def\Rp{R^{\ts\prime}}

\def\s{\mathfrak{s}}
\def\S{\operatorname{S}}
\def\si{\sigma}
\def\sib{\ts\overline{\ns\si\ns}\ts}
\def\sih{\widetilde{\si}}
\def\sip{\si^{\ts\prime}}
\def\so{\mathfrak{so}}
\def\sp{\mathfrak{sp}}
\def\Sp{S^{\ts\prime}}
\def\St{\widetilde{S}}
\def\Sym{\mathfrak{S}}

\def\th{\theta}
\def\Tp{T^{\ts\prime}}
\def\ts{\hskip1pt}
\def\tts{\hskip.5pt}

\def\uo{\ts\overline{\ns u\nns}\,}
\def\U{\operatorname{U}}
\def\Uhb{\,\overline{\!\U(\h)\!\!\!}\,\,\,}

\def\vuo{\ts\overline{\ns v\ot u\ns}\,}

\def\Wb{\ts\overline{\ns W^{\phantom{!}}\!\!}}
\def\Wt{\widetilde{W}}

\def\X{\operatorname{X}}
\def\xic{\check\xi}

\def\Y{\operatorname{Y}}

\def\Z{\operatorname{Z}}
\def\ZZ{\mathbb{Z}}

\begin{document} 
 

\title{Twisted Yangians and Mickelsson Algebras I}

\author{Sergey Khoroshkin}
\address{%
Institute for Theoretical and Experimental Physics\\
Moscow 117259, Russia}
\email{\texttt{khor@itep.ru}}

\author{Maxim Nazarov}
\address{%
Department of Mathematics, University of York\\
York YO10 5DD, England}
\email{\texttt{mln1@york.ac.uk}}


\begin{abstract}
We introduce an analogue of the composition
of the Cherednik and Drinfeld functors
for twisted Yangians. Our definition is based on 
the Howe duality, and originates from the centralizer
construction of twisted Yangians due to Olshanski.
Using our functor, 
we establish a correspondence
between intertwining operators on the tensor products of 
certain modules over twisted Yangians,
and the extremal cocycle on the hyperoctahedral group.
\end{abstract}


\maketitle


\renewcommand{\theequation}{\thesection.\arabic{equation}} 
\makeatletter                                          
\@addtoreset{equation}{section}                        
\makeatother    


\noindent{\Large\bf 0.\ Introduction}
\setcounter{equation}{0}

\vspace{6pt}
\noindent
This article is a continuation of our work \cite{KN1} which
concerned two well known functors. The definition of one of these two
functors belongs to V.\,Drinfeld \cite{D2}. Let $\D_N$ be the
\textit{degenerate affine Hecke algebra\/} corresponding to the
general linear group $GL_N$ over a non-Archimedean local field. This
is an associative algebra over the 
field $\CC$
which contains the symmetric group ring $\CC\,\Sym_N$ as a subalgebra.
Let $\Y(\gl_n)$ be the 
\textit{Yangian\/} of the general linear Lie algebra $\gl_n\ts$.
This is a deformation of the universal enveloping algebra
of the polynomial current Lie algebra 
$\gl_n[u]$ in the class of Hopf algebras \cite{D1}.
It contains the universal enveloping algebra $\U(\gl_n)$ as a subalgebra.
There is also a homomorphism of associative algebras 
$\Y(\gl_n)\to\U(\gl_n)$ identical on the subalgebra
$\U(\gl_n)\subset\Y(\gl_n)\ts$.
In \cite{D2} for any $\D_N\ts$-module $M\ts$, an action of the algebra
$\Y(\gl_n)$ was defined on the vector space
$(\ts M\ot(\CC^{\ts n})^{\ot N}){}^{\ts\Sym_N}$
of diagonal $\Sym_N$-invariants in the tensor 
product of the vector spaces $M$ and $(\CC^{\ts n})^{\ot N}\ts$.
So one gets a functor from the category of all 
$\D_N\ts$-modules to the category of $\Y(\gl_n)\ts$-modules, 
called the \textit{Drinfeld functor}.

In \cite{KN1} we studied the composition of the Drinfeld functor
with another functor, introduced by I.\,Cherednik \cite{C}.
This second functor was also studied by
T.\,Arakawa, T.\,Suzuki and A.\,Tsuchiya \cite{A,AS,AST}. 
For any module $U$ over the Lie algebra $\gl_{\ts l}\ts$,
an action of the algebra $\D_N$ can be defined on the
tensor product $U\ot(\CC^{\ts l})^{\ot N}$
of $\gl_{\ts l}\ts$-modules. This action of $\D_N$
commutes with the diagonal action of $\gl_{\ts l}$ on the tensor product.
So one gets a functor from the category of all 
$\gl_{\ts l}\ts$-modules to the category of bimodules over $\gl_{\ts l}$ 
and $\D_N\ts$, called the \textit{Cherednik functor}.
By applying the Drinfeld functor to the $\D_N\ts$-module 
$M=U\ot(\CC^{\ts l})^{\ot N}\ts$, one turns to an $\Y(\gl_n)\ts$-module
the vector space
$$
(\ts U\ot(\CC^{\ts l})^{\ot N}\ot(\CC^{\ts n})^{\ot N}){}^{\ts\Sym_N}
=
U\ns\ot\ts\operatorname{S\ts}^N(\CC^{\ts l}\ot\CC^{\ts n})\ts.
$$
The action of the associative algebra 
$\Y(\gl_n)$ on this vector space
commutes with the action of 
$\gl_{\ts l}\ts$. 
By taking the direct sum of these $\Y(\gl_n)\ts$-modules
over $N=0,1,2,\ldots$
we turn into an $\Y(\gl_n)\ts$-module
the space $U\ns\ot\ts\S\ts(\CC^{\ts l}\ot\CC^{\ts n})\ts$.
It is also a $\gl_{\ts l}\ts$-module; denote this bimodule by
$\E_{\ts l}\ts(U)\ts$. We identify the symmetric algebra
$\S\ts(\CC^{\ts l}\ot\CC^{\ts n})$ with the ring
$\P\ts(\CC^{\ts l}\ot\CC^{\ts n})$ of polynomial functions on
$\CC^{\ts l}\ot\CC^{\ts n}$, and denote by 
$\PD\ts(\CC^{\ts l}\ot\CC^{\ts n})$ the algebra
of the differential operators on the vector space
$\CC^{\ts l}\ot\CC^{\ts n}$ with polynomial coefficients.
The action of the Yangian
$\Y(\gl_n)$ on its module 
$\E_{\ts l}\ts(U)$ is then determined by a homomorphism 
$\al_{\ts l}:\Y(\gl_n)\to
\U(\gl_{\ts l})\ot\PD\ts(\CC^{\ts l}\ot\CC^{\ts n})\ts$,
see Proposition~\ref{dast} below.

Now let $\f_m$ be either the symplectic Lie algebra
$\sp_{2m}$ or the orthogonal Lie algebra $\so_{2m}\ts$.
The first objective of the present article
is to define analogues of the functor $\E_{\ts l}$
and of the homomorphism $\al_{\ts l}$
for the Lie algebra $\f_m$ instead of $\gl_{\ts l}\ts$.
The role of the Yangian $\Y(\gl_n)$ is played here by the
\textit{twisted Yangian\/} $\Y(\g_n)\ts$, which is
a right coideal subalgebra of the Hopf algebra $\Y(\gl_n)\ts$.
Here $\g_n$ is a Lie subalgebra of $\gl_n\ts$, orthogonal
in the case $\f_m=\sp_{2m}$ and symplectic in the case $\f_m=\so_{2m}\ts$;
in the latter case $n$ has to be even.
Let the superscript ${}^{\ts\prime}$ indicate the transposition
in $\gl_n$ relative to the bilinear form on $\CC^n$
preserved by the subalgebra $\g_n\subset\gl_n\ts$, so that
$
\g_n=\{\ts A\in\gl_n\,\ts|\,A^{\ts\prime}=-A\ts\}\ts.
$
As an associative algebra, 
$\Y(\g_n)$ is a deformation of the
universal enveloping algebra of the 
\textit{twisted polynomial current Lie algebra\/}
$$
\{\ts A(u)\in\gl_n[u]\,\ts|\,A^{\ts\prime}(u)=-A(-u)\ts\}\ts.
$$

Twisted Yangians were introduced
by G.\,Olshanski \cite{O2}.
In Section 2 of the present article we introduce a homomorphism
$\Y(\g_n)\to\U(\ts\f_m)\ot\PD\ts(\CC^{\ts m}\ot\CC^{\ts n})\ts$,
see our Propositions~\ref{xb} and \ref{tb}.
The image of $\Y(\g_n)$ under this homomorphism commutes with the image
of the algebra $\U(\ts\f_m)$ under its diagonal embedding 
\eqref{xact}
to the tensor product $\U(\ts\f_m)\ot\PD\ts(\CC^{\ts m}\ot\CC^{\ts n})\,$;
here we use the homomorphism
$\zeta_{\ts n}:\U(\ts\f_m)\to\PD\ts(\CC^{\ts m}\ot\CC^{\ts n})$
defined by \eqref{gan}.
The twisted Yangian $\Y(\g_n)$
contains the universal enveloping algebra $\U(\g_n)$ as a subalgebra.
Moreover, there is a homomorphism
$\pi_n:\Y(\g_n)\to\U(\g_n)$ identical on the subalgebra
$\U(\g_n)\subset\Y(\g_n)\ts$.
Our results extend the classical theorem \cite{H}
stating that the image of $\U(\ts\f_m)$ in
$\PD\ts(\CC^{\ts m}\ot\CC^{\ts n})$
under the homomorphism
$\zeta_n$ consists of all $G_n$-invariant elements.
Here $G_n$ is either the orthogonal or the
symplectic group, so that $\g_n$ is its Lie algebra;
the group $G_n$ acts on 
$\PD\ts(\CC^{\ts m}\ot\CC^{\ts n})$
via its natural action on $\CC^{\ts n}\ts$.

A certain analogue of the algebra $\D_N$ 
was defined in \cite{N0}, it may be called \cite{GH} the 
\textit{degenerate affine Birman-Wenzl-Murakami algebra}. 
In \cite{N0} the role of the subalgebra $\CC\,\Sym_N\subset\D_N$ was played 
by the Brauer algebra \cite{B}. 
By using that analogue, R.\,Orellana and A.\,Ram \cite{OR}
have found counterparts of the Cherednik functor 
for the categories of modules over
the Lie algebras $\sp_{2m}\com\so_{2m}$ and $\so_{2m+1}\ts$.
Still there is no analogue of the Drinfeld functor
available for the category of modules over the 
degenerate affine Birman-Wenzl-Murakami algebra.

In this article, we do not use the functors from \cite{OR}.
Instead, we introduce an analogue of the functor $\E_{\ts l}$
for the Lie algebra $\f_m\ts$. The functor $\E_{\ts l}$ was defined 
as a composition of two functors, due to Cherednik and Drinfeld. But here
we do not work with the two functors separately from each other.
Unlike the natural action of
$\gl_{\ts l}$ on $\P\ts(\CC^{\ts l}\ot\CC^{\ts n})\ts$,
the action of $\f_m$ on $\P\ts(\CC^{\ts m}\ot\CC^{\ts n})$
defined by the homomorphism $\zeta_{\ts n}$ does not preserve 
the subspace of polynomials of any degree $N$,
so the results of \cite{N0,OR} do not apply.
Our approach was inspired by a joint work of
I.\,Grojnowski and M.\,Kleber
which did not involve twisted Yangians however\tts;~see~\cite{G,K}.
\\
\indent
We prefer to work with a certain central extension
$\X(\g_n)$ of the algebra 
$\Y(\g_n)\ts$, called the 
\textit{extended twisted Yangian}. Central elements
$\,O^{(1)}\ns\com O^{(2)}\ns\com\,\ldots\,$ of the algebra $\X(\g_n)$
generating the kernel of the canonical homomorphism $\X(\g_n)\to\Y(\g_n)$
are given in Section 1, together with the definitions of 
$\X(\g_n)$ and $\Y(\g_n)\ts$. There is also a homomorphism
$\X(\g_n)\to\X(\g_n)\ot\Y(\gl_n)\ts$.
Using it, the tensor product
of any modules over the algebras $\X(\g_n)$ and $\Y(\gl_n)$
becomes another module over $\X(\g_n)\ts$.
Moreover, this homomorphism is a coaction
of the Hopf algebra $\Y(\gl_n)$ on the algebra $\X(\g_n)\ts$.
We define a homomorphism
$\be_m:\X(\g_n)\to\U(\ts\f_m)\ot\PD\ts(\CC^{\ts m}\ot\CC^{\ts n})$
which is our analogue of the homomorphism $\al_{\ts l}\ts$, see
Proposition~\ref{xb}. 
The image of $\X(\g_n)$ under $\be_m$ commutes with the image
of the algebra $\U(\ts\f_m)$ under its diagonal embedding to
$\U(\ts\f_m)\ot\PD\ts(\CC^{\ts m}\ot\CC^{\ts n})\ts$.
The reason why we work with $\X(\g_n)$ rather than with 
$\Y(\g_n)$ is explained in Section~2.
\\
\indent
The generators of the algebra $\X(\g_n)$ appear as coefficients
of certain series $S_{ij}(u)$ in the variable $u$
where $i\com j=1\lcd n\ts$.
We define the homomorphism $\be_m$ by applying it to these
coefficients, and by giving the resulting series
$\be_m(S_{ij}(u))$ with the coefficients in the algebra 
$\U(\ts\f_m)\ot\PD\ts(\CC^{\ts m}\ot\CC^{\ts n})\ts$. 
Then we define another homomorphism
$\bet_m:\X(\g_n)\to\U(\ts\f_m)\ot\PD\ts(\CC^{\ts m}\ot\CC^{\ts n})$
which factors through canonical homomorphism
$\X(\g_n)\to\Y(\g_n)\ts$. 
Hence we obtain the homomorphism
$\Y(\g_n)\to\U(\ts\f_m)\ot\PD\ts(\CC^{\ts m}\ot\CC^{\ts n})$ mentioned above.
Every series $\bet_m(S_{ij}(u))$ is the product
of $\be_m(S_{ij}(u))$ with a certain series
with coefficients from 
$\Z(\ts\f_m)\ot1\ts$,
where $\Z(\ts\f_m)$ is the centre of the algebra $\U(\ts\f_m)\ts$.
Now let $V$ be any $\f_m\ts$-module. Using the homomorphism $\be_m\ts$,
we turn the vector space $V\ns\ot\ts\P\ts(\CC^{\ts m}\ot\CC^{\ts n})$
into a bimodule over $\f_m$ and $\X(\g_n)\ts$.
We denote this bimodule by $\F_m(V)\ts$.
The functor $\F_m$ is our analogue of the functor $\E_{\ts l}$ for
$\f_m$ instead of $\gl_{\ts l}\ts$. 
When $m=0\ts$, we set $\F_{\ts0}(V)=\CC$ so that
$\be_{\ts0}$ is the composition of the canonical
homomorphism $\X(\g_n)\to\Y(\g_n)$ with the restriction of 
the counit homomorphism $\Y(\gl_n)\to\CC$ to $\Y(\g_n)\ts$.
\\
\indent
Here we show that the functor $\F_m$
shares the three fundamental properties of the functor $\E_{\ts l}$
considered in \cite{KN1}. The first of these properties of
$\E_{\ts l}$ concerns \textit{parabolic induction\/}
from the direct sum of Lie algebras
$\gl_m\op\gl_{\ts l}$ to $\gl_{m+l}\ts$.
Let $\p$ be the maximal parabolic subalgebra of $\gl_{m+l}$
containing the direct sum $\gl_m\op\gl_{\ts l}\ts$.
Let $\q\subset\gl_{m+l}$ be the Abelian subalgebra with
$
\gl_{m+l}=\q\op\p\ts.
$
For any $\gl_m\ts$-module $W$ let $W\bt U$ be the
$\gl_{m+l}\ts$-module 
parabolically induced from the 
$\gl_m\op\gl_{\ts l}\ts$-module $W\ot U\ts$.
This is a module induced from the subalgebra $\p\ts$.
Consider the space $\E_{\ts m+l}\ts(\ts W\ns\bt U\ts)_{\ts\q}$
of $\q\ts$-coinvariants of the $\gl_{m+l}\ts$-module 
$\E_{\ts m+l}\ts(\ts W\ns\bt U\ts)\ts$.
This space is an $\Y(\gl_n)\ts$-module, which also
inherits the action of the Lie algebra $\gl_m\op\gl_{\ts l}\ts$.
The additive group $\CC$ acts on 
the Hopf algebra $\Y(\gl_n)$ by automorphisms.
Denote by $\E_{\ts l}^{\,z}\ts(U)$ the $\Y(\gl_n)\ts$-module obtained
from $\E_{\ts l}\ts(U)$ by pulling it back through the automorphism
of $\Y(\gl_n)$ corresponding to the element $z\in\CC\ts$.
The automorphism itself is denoted by $\tau_z\ts$, see \eqref{tauz}. 
As a $\gl_{\ts l}\ts$-module $\E_{\ts l}^{\,z}\ts(U)$
coincides with $\E_{\ts l}\ts(U)\ts$.
In \cite{KN1} we proved that
the bimodule $\E_{\ts m+l}\ts(\ts W\ns\bt U\ts)_{\ts\q}$
of $\Y(\gl_n)$ and $\gl_m\op\gl_{\ts l}$
is equivalent to 
$\E_{\ts m}(W)\ot\E_{\ts l}^{\ts m}\ts(U)\ts$.
There we used the comultiplication on $\Y(\gl_n)\ts$.

Our Theorem~\ref{parind} is an analogue of this 
comultiplicative property of 
$\E_{\ts l}\ts$. Take the maximal parabolic subalgebra
of the Lie algebra $\f_{m+l}$ containing the direct sum
$\f_m\op\gl_{\ts l}\,$; we do not exclude the case $m=0$ here.
Using that subalgebra, determine
the $\f_{m+l}\ts$-module $V\bt U$ parabolically induced from
the $\f_m\op\gl_{\ts l}\ts$-module $V\ot U\ts$.
Consider the space of coinvariants of the
$\f_{m+l}\ts$-module $\F_{m+l}\ts(\ts V\ns\bt U\ts)$
relative to the nilpotent subalgebra of $\f_{m+l}$
complementary to our parabolic subalgebra. 
This space is a bimodule over $\f_m\op\gl_{\ts l}$
and $\X(\g_n)\ts$. We prove that this bimodule
is essentially equivalent to the tensor product
$\F_m(V)\ot\E_{\ts l}^{\ts z}\ts(U)$ with
$z=m+\frac12$ for $\f_m=\sp_{2m}\ts$, and
$z=m-\frac12$ for $\f_m=\so_{2m}\ts$.
Namely, the action of $\X(\g_n)$ on this module
$\F_m(V)\ot\E_{\ts l}^{\ts z}\ts(U)$ is
defined by a composition of two homomorphisms,
$$
\X(\g_n)\to
\X(\g_n)\ot\Y(\gl_n)\to
\U(\ts\f_m)\ot\PD\ts(\CC^{\ts m}\ot\CC^{\ts n})\ot
\U(\ts\gl_{\ts l})\ot\PD\ts(\CC^{\ts l}\ot\CC^{\ts n})\,.
$$
By replacing the image of $S_{ij}(u)$
under this composition by a product of this image
with a certain series with the coefficients from the subalgebra
$$
1\ot1\ot\Z(\gl_{\ts l})\ot1
\ts\subset\ts
\U(\ts\f_m)\ot\PD\ts(\CC^{\ts m}\ot\CC^{\ts n})\ot
\U(\ts\gl_{\ts l})\ot\PD\ts(\CC^{\ts l}\ot\CC^{\ts n})\,,
$$
we define another action of $\X(\g_n)$ on the 
vector space of $\F_m(V)\ot\E_{\ts l}^{\ts z}\ts(U)\ts$. 
This is equivalent to the action of $\X(\g_n)$
on the space of coinvariants of $\F_{m+l}\ts(\ts V\ns\bt U\ts)\ts$.
The actions of the direct summand $\f_m$ of $\ts\f_m\op\gl_{\ts l}$
on $\F_m(V)\ot\E_{\ts l}^{\ts z}\ts(U)$
and on the space of coinvariants of 
$\F_{m+l}\ts(\ts V\ns\bt U\ts)$ are the same, while
the actions of the direct summand $\gl_{\ts l}$
differ only by the automorphism \eqref{autol} of 
the Lie algebra $\gl_{\ts l}\ts$.
Our proof of Theorem~\ref{parind} is rather technical,
and is organized as a sequence of three lemmas and five propositions.

Using the functor $\E_{\ts l}\ts,$
in \cite{KN1} we gave a representation theoretic 
explanation of the correspondence between intertwining operators 
of the tensor products of certain $\Y(\gl_n)\ts$-modules, and
the \textit{extremal cocycle\/} 
on the Weyl group $\Sym_{\ts l}$ of the reductive Lie algebra $\gl_{\ts l}$ 
defined by D.\,Zhelobenko \cite{Z1,Z2}.
This correspondence, discovered
by V.\,Tarasov and A.\,Varchenko \cite{TV}, was one of the motivations
of our work. 
The arguments of \cite{TV}, inspired by
the results of V.\,Toledano-Laredo \cite{T}, were based on
the classical duality theorem \cite{H}
stating that the images of 
$\U(\gl_m)$ and $\U(\gl_n)$ in the ring
$\PD\ts(\CC^{\ts m}\otimes\CC^{\ts n})$
of differential operators on $\CC^{\ts m}\otimes\CC^{\ts n}$
with polynomial coefficients
are the commutants of each other.
Our explanation of this correspondence was based on the theory of
Mickelsson algebras \cite{M} as developed in \cite{KO}.
From a different perspective, connections between the 
representation theory of the Yangian $\Y(\gl_n)$ and 
the theory of Mickelsson algebras were studied by A.\,Molev~\cite{M1}.
We considered the diagonal embedding \eqref{eabact}
of the algebra $\U(\gl_{\ts l})$
to the tensor product 
$\U(\gl_{\ts l})\ot\PD\ts(\CC^{\ts l}\ot\CC^{\ts n})\,$;
there we used the homomorphism 
$\U(\gl_{\ts l})\to\PD\ts(\CC^{\ts l}\ot\CC^{\ts n})$
corresponding to the natural action of 
$\gl_{\ts l}$ on $\P\ts(\CC^{\ts l}\ot\CC^{\ts n})\ts$.
The Mickelsson algebra which we employed was determined by the pair
formed by this tensor product algebra, and its subalgebra 
$\U(\gl_{\ts l})$ relative to this embedding.
We regard the role played by $\E_{\ts l}$
as another fundamental property of that functor. 

Here we show that the functor $\F_m$ plays a similar role
for the Lie algebra $\f_m$ instead of $\gl_{\ts l}\ts$.
We establish a correspondence between intertwining operators
of certain $\X(\g_n)\ts$-modules, and the 
extremal cocycle
on the hyperoctahedral group $\H_m$ corresponding
to the reductive Lie algebra $\f_m\ts$.
Here $\H_m$ is regarded as the Weyl group of $\f_m=\sp_{2m}\ts$,
and as an extension of the Weyl group of $\f_m=\so_{2m}$
by a Dynkin diagram automorphism. In both cases, the definition
of the extremal cocycle
is essentially due to D.\,Zhelobenko \cite{Z1,Z2}.
However, the original extremal cocycle
has been defined on the Weyl group of $\f_m$ itself,
which in the case $\f_m=\so_{2m}$ is a subgroup of $\H_m$ of index~2.  
An extension of the original definition to the
whole group $\H_m$ is a new feature in this case. 
All necessary details on the extremal cocycle corresponding to $\f_m$
are given in Section~4.

For any positive integer $N$ and for any $z\in\CC\ts$ let $P_z^{\ts N}$ be
the $\Y(\gl_n)\ts$-module obtained by pulling back the standard action
of $\U(\gl_n)$ on the space of polynomials on $\CC^{\ts n}$
of degree $N$ through the homomorphism $\Y(\gl_n)\to\U(\gl_n)$ 
and then through the automorphism $\tau_z$ of $\Y(\gl_n)\ts$.
In \cite{KN1}, the value of the extremal cocycle at a permutation
from $\Sym_{\ts l}$ corresponds to an intertwining operator between
two $l\ts$-fold tensor products of $\Y(\gl_n)\ts$-modules of
the form $P_z^{\ts N}\ts$. We assume that
the differences between the parameters $z$ for the
$l$ tensor factors of the source product are not in $\ZZ\ts$.
The $l$ tensor factors of the target product are obtained 
from those of the source product
by that permutation from $\Sym_{\ts l}\ts$,
where the extremal cocycle has been evaluated. 
It is well known that under our assumption on the parameters $z$
both source and target tensor products are
irreducible $\Y(\gl_n)\ts$-modules, equivalent to each other\ts;
see \cite{NT} for a more general result. Hence an
intertwining operator between these two tensor products
is unique up to a multiplier from $\CC\ts$.

The twisted Yangian $\Y(\g_n)$ is determined \cite{N2} by
a distinguished involutive automorphism
of the algebra $\Y(\gl_n)\ts$. It corresponds to the 
automorphism $A(u)\mapsto-\ts A^{\ts\prime}(-u)$
of the Lie algebra $\gl_n[u]\ts$,
when $\Y(\gl_n)$ is regarded as a deformation of 
the universal enveloping algebra of $\gl_n[u]\ts$.
By pulling the $\Y(\gl_n)\ts$-module $P_z^{\ts N}$ back through that
automorphism of $\Y(\gl_n)$ we get another
$\Y(\gl_n)\ts$-module, which we denote by $P_z^{\ts-N}\ts$.
Now take any positive integers $\nu_1\lcd\nu_m$ 
and any $z_1\lcd z_m\in\CC$ such that $\ts z_a-z_b\notin\ZZ\ts$ 
and $\ts z_a+z_b\notin\ZZ\ts$ whenever $a\neq b\ts$. In the case
$\f_m=\sp_{2m}$ we also assume that $2\ts z_a\notin\ZZ$
for any index~$a\ts$.

The hyperoctahedral group $\H_m$ can be realized
as the group of permutations $\si$ of
$-\ts m\lcd-1\com1\lcd m\ts$
such that $\si\ts(-c)=-\ts\si\ts(c)$ for any $c\ts$. 
In Section~5 we show how the value of the extremal cocycle 
related to the Lie algebra $\f_m$ at any element $\si\in\H_m$ 
determines an intertwining operator of $\X(\g_n)\ts$-modules 
\begin{equation}
\label{nuz}
P_{z_m}^{\,\nu_m}
\ot\ldots\ot
P_{z_1}^{\,\nu_1}
\,\to\,
P_{\widetilde{z}_m}^{\,\de_m\ts\widetilde{\nu}_m}
\ot\ldots\ot
P_{\widetilde{z}_1}^{\,\de_1\ts\widetilde{\nu}_1}
\end{equation}
where
$$
\widetilde\nu_a=\nu_{\ts|\si^{-1}(a)|}\ts,
\quad 
\widetilde{z}_a=z_{\ts|\si^{-1}(a)|}
\quad\textrm{and}\quad
\de_a=\operatorname{sign}\si^{-1}(a)
$$
for each $a=1\lcd m\ts$.
The tensor products in \eqref{nuz} are those of
$\Y(\gl_n)\ts$-modules. By restricting both
tensor products to the subalgebra 
$\Y(\g_n)\subset\Y(\gl_n)$ and by pulling the 
restrictions back through the canonical homomorphism 
$\X(\g_n)\to\Y(\g_n)$,
both tensor products in \eqref{nuz} become 
$\X(\g_n)\ts$-modules. The same $\ts\X(\g_n)\ts$-modules
can also be obtained by using the coaction of the Hopf algebra
$\Y(\gl_n)$ on $\X(\g_n)\ts$.
It has been proved in \cite{MN} that
under our assumptions on $z_1\lcd z_m$
both the source and the target tensor products in \eqref{nuz}
are irreducible $\X(\g_n)\ts$-modules, equivalent to each other.
Hence an intertwining operator between them 
is unique up to a multiplier from $\CC\ts$.
For our intertwining operator,
this multiplier is determined by Proposition~\ref{isis}.
Our proof of that proposition employs a well known
formula \eqref{gaga} for the value of \textit{hypergeometric function}
$\mathrm{F}\tts(\tts u\com v\com w\ts;z\tts)$ at $z=1\ts$.

To obtain our intertwining 
operator \eqref{nuz}
we use the theory of Mickelsson algebras, like we did in \cite{KN1}.
Our particular Mickelsson algebra is determined by the pair
formed by the tensor product
$\U(\ts\f_m)\ot\PD\ts(\CC^{\ts m}\ot\CC^{\ts n})$
and by its subalgebra 
$\U(\ts\f_m)$ relative to the 
embedding \eqref{xact}.
The extended twisted Yangian $\X(\g_n)$ appears here naturally, 
because its image relative to the homomorphism 
$\be_m$ commutes with the image of $\U(\ts\f_m)$ in the tensor product.
Another expression for an intertwining operator
of the $\X(\g_n)\ts$-modules 
\eqref{nuz} was given in \cite{N2}.

In Section 2 we choose a triangular decomposition \eqref{tridec}
of the Lie algebra $\f_m$ into a direct sum
of a Cartan subalgebra $\h$ and of two maximal nilpotent subalgebras
$\n\,,\np\ts$. 
For any formal power series $f(u)$ in $u^{-1}$ with coefficients 
from $\CC$ and leading term $1$, the assignments \eqref{fus}
define an automorphism of the algebra $\X(\g_n)\ts$. 
Up to pulling it back through such an automorphism,
the source $\X(\g_n)\ts$-module in \eqref{nuz} 
arises as the space of $\n\ts$-coinvariants of weight $\la$
for the $\f_m\ts$-module $\F_m(M_\mu)\ts$,
where $M_\mu$ is the Verma module over $\f_m$ with the
highest vector of weight $\mu$ annihilated by the
action of the subalgebra $\np\subset\f_m\ts$.
The weights $\la$ and $\mu$
relative to the Cartan subalgebra $\h$ 
can be determined by the parameters $\nu_1\lcd\nu_m$ 
and $z_1\lcd z_m$ taken here\ts; see Corollary \ref{verma} for details.

To get the target $\X(\g_n)\ts$-module in \eqref{nuz},
we generalize our definition of the functor $\F_m\ts$.
In the beginning of Section 5, 
for any sequence $\de=(\ts\de_1\lcd\de_m)$ of $m$
elements of the set $\{1\com-1\}$ we define a functor
$\F_\de\,$, with the same source and target categories as
the functor $\F_m$ has. Moreover, for any $\f_m\ts$-module
$V$ the underlying vector spaces 
of the bimodules $\F_\de(V)$ and $\F_m(V)$ are the same,
that is $V\ot\P\ts(\CC^{\ts m}\ot\CC^{\ts n})\ts$.
The actions of $\f_m$ and $\X(\g_n)$ on $\F_\de(V)$ 
are obtained by pushing forward the defining homomorphisms 
$$
\zeta_n:\U(\ts\f_m)\to\PD\ts(\CC^{\ts m}\ot\CC^{\ts n})
\quad\text{and}\quad
\be_m:\X(\g_n)\to\U(\ts\f_m)\ot\PD\ts(\CC^{\ts m}\ot\CC^{\ts n})
$$
through a certain authomorphism of the ring
$\PD\ts(\CC^{\ts m}\ot\CC^{\ts n})$ depending on the sequence $\de\ts$.
In particular, we have $\F_\de\ts(V)=\F_m(V)$ for 
the sequence $\de=(1\lcd1)\ts$.
Up to pulling it back through an automorphism of the form \eqref{fus},
the target $\X(\g_n)\ts$-module in \eqref{nuz} 
arises as the space of $\n\ts$-coinvariants of weight $\si\ts\circ\la$
for the $\f_m\ts$-module $\F_\de\ts(M_{\ts\si\ts\circ\ts\mu})\ts$.
The symbol $\circ$ indicates the \textit{shifted action\/}
of the group $\H_m$ on the weights of $\h\ts$.
Moreover, the authomorphism \eqref{fus} here does not depend
on $\si\in\H_m\ts$, and hence is the same for the source and the target 
modules in \eqref{nuz}\ts.

The third fundamental property of the functor $\E_{\ts l}$ considered
in \cite{KN1} is its connection with the 
\textit{centralizer construction\/} 
of the Yangian $\Y(\gl_n)$ proposed by G.\,Olshanski \cite{O1}.
For any two irreducible polynomial modules
$U$ and $U^{\ts\prime}$ over the Lie algebra $\gl_{\ts l}\ts$,
the results of \cite{O1} provide an action
of $\Y(\gl_n)$ on the vector space
$$
\Hom_{\,\gl_{\ts l}}(\ts U^{\ts\prime}\ts,
U\ot\S\ts(\ts\CC^{\ts l}\ot\CC^{\ts n}\ts))\ts.
$$
Moreover, this action is irreducible. In \cite{KN1}
we proved that the same action is inherited from the bimodule 
$\E_{\ts l}\ts(U)=U\ns\ot\ts\S\ts(\CC^{\ts l}\ot\CC^{\ts n})$
over $\Y(\gl_n)$ and $\gl_{\ts l}\ts$.

There is a centralizer construction of 
$\Y(\g_n)$ again due
to G.\,Olshanski \cite{O2}, see also \cite{MO} and Section 6 here.
That construction
served as a motivation for introducing the twisted Yangians.
Let $V$ and $V^{\prime}$ be irreducible highest weight modules
of the Lie algebra $\f_m\ts$,
corresponding to unitary
representations of the real metaplectic group
$Mp_{\ts 2m}$ if $\f_m=\sp_{2m}\ts$, 
or the real group $SO_{2m}^\ast$ if $\f_m=\so_{2m}\ts$.
The results of \cite{O2} provide an action
of the algebra $\X(\g_n)$ on the vector 
space
\begin{equation}
\label{vvp}
\Hom_{\,\f_m}(\ts V^{\prime}\ts,
V\ot\P\ts(\ts\CC^{\ts m}\ot\CC^{\ts n}\ts))\,.
\end{equation}
Here we employ the theory of \textit{reductive dual pairs}
due to R.\,Howe, see \cite{KV} and \cite{EHW,EP}.
The group $G_n$ also acts on this vector space,
via its natural action on $\CC^{\ts n}$. When $\g_n$ is an
orthogonal Lie algebra, the space \eqref{vvp} is irreducible
under the joint action of $\X(\g_n)$ and $G_n\ts$. 
When $\g_n$ is symplectic, 
\eqref{vvp} is irreducible under the action of the algebra $\X(\g_n)$ alone. 
Our Theorem~\ref{5.1} asserts that the action
of $\X(\g_n)$ on \eqref{vvp} is essentially the same as
the action inherited from the bimodule
$\F_m(V)=V\ns\ot\ts\P\ts(\CC^{\ts m}\ot\CC^{\ts n})$
over $\X(\g_n)$ and $\f_m\ts$. More precisely, the action
of $\X(\g_n)$ on \eqref{vvp}
provided by \cite{O2} can also be obtained by
pulling the action of $\X(\g_n)$ on $\F_m(V)$
back through another automorphism of $\X(\g_n)$ of the form \eqref{fus}.
This third property of $\F_m$ was the origin
of our definition of this functor. 
Thus here we have two different descriptions of
the same action of $\X(\g_n)$ on \eqref{vvp}.
Another two, still different descriptions of the same action of $\X(\g_n)$
on the vector space \eqref{vvp} were provided in
\cite{M2} and \cite{N2}~respectively.

Finally, we note that there is an \lq\lq\ts antisymmetric\ts\rq\rq\
version of the functor $\E_{\ts l}$ for the Lie algebra
$\gl_{\ts l}\ts$, this version has been used by T.\,Arakawa \cite{A}.
There the symmetric algebra 
$\S\ts(\CC^{\ts l}\ot\CC^{\ts n})$ is replaced by the
exterior algebra $\mathrm{\Lambda}\ts(\CC^{\ts l}\ot\CC^{\ts n})\ts$.
The corresponding analogues of the three fundamental properties of
$\E_{\ts l}$ were considered in \cite{KN2}. 
An \hbox{\lq\lq\ts antisymmetric\ts\rq\rq}
version of the functor $\F_m$ for the Lie algebra 
$\f_m$ will be studied in a sequel to this article.


\enlargethispage{20pt}

\section{Twisted Yangians}
\setcounter{section}{1}
\setcounter{equation}{0}
\setcounter{theorem*}{0}

Let $G_n$ be one of the complex Lie groups $O_n$ and $Sp_n\ts$. 
We regard $G_n$ as the subgroup  of the general linear Lie group 
$GL_n\ts$, preserving a non-degenerate bilinear form $\langle\ ,\,\rangle$ 
on the vector space $\CC^{\ts n}\ts$. This form is symmetric in the case
$G_n=O_n\ts$, and alternating in the case $G_n=Sp_n\ts$. In the latter case
$n$ has to be even. We always assume that the integer $n$ is positive.
Throughout this article, we will use the following convention. Whenever
the double sign $\ts\pm\ts$ or $\ts\mp\ts$ appears, the upper sign corresponds
to the case of a symmetric form on $\CC^{\ts n}$ while the lower sign 
corresponds to the case of an alternating form on $\CC^{\ts n}$.

Let $i$ be any of the indices $1\lcd n\ts$. If $i$ is even, put
$\bi=i-1\ts$. If $i$ is odd and $i<n$, put $\bi=i+1\ts$. Finally,
if $i=n$ and $n$ is odd, put $\bi=i\ts$.
Let $e_1\lcd e_n$ be the vectors of the standard basis in $\CC^{\ts n}$.
Choose the bilinear form on $\CC^{\ts n}$ so that for any two basis vectors 
$e_i$ and $e_j$ we have
$\langle\ts e_i\com e_j\ts\rangle=\th_i\,\de_{\ts\bi j}$ where 
$\th_i=1$ or $\th_i=(-1)^{\ts i-1}$
in the case of the symmetric or alternating form.

Let $E_{ij}\in\End(\CC^{\ts n})$ be the standard matrix units.
We will also regard these matrix units as basis elements of
the general linear Lie algebra $\gl_n\ts$.
Let $\g_n$ be the Lie algebra of the group $G_n\ts$,
so that $\g_n=\so_n$ or $\g_n=\sp_n$
in the case of the symmetric or alternating form on $\CC^{\ts n}$.
The Lie subalgebra $\g_n\subset\gl_n$ is spanned by the elements
$
E_{ij}-\th_i\ts\th_j\ts E_{\ts\bj\ts\bi}\,.
$

Take the \textit{Yangian\/} $\Y(\gl_n)$ of the Lie algebra $\gl_n\ts$.
The unital associative algebra $\Y(\gl_n)$ over $\CC$
has a family of generators 
$
T_{ij}^{\ts(1)},T_{ij}^{\ts(2)},\ts\ldots
$
where $i\com j=1\lcd n\ts$.
Defining relations for these generators
can be written using the series
\begin{equation*}
T_{ij}(u)=
\de_{ij}+T_{ij}^{\ts(1)}u^{-\ns1}+T_{ij}^{\ts(2)}u^{-\ns2}+\,\ldots
\end{equation*}
where $u$ is a formal parameter. Let $v$ be another formal parameter.  
Then the defining relations in the associative algebra $\Y(\gl_n)$
can be written as
\begin{equation}
\label{yrel}
(u-v)\,[\ts T_{ij}(u)\ts,T_{kl}(v)\ts]\ts=\;
T_{kj}(u)\ts T_{il}(v)-T_{kj}(v)\ts T_{il}(u)\,.
\end{equation}

The algebra $\Y(\gl_n)$ is commutative if $n=1\ts$.
By \eqref{yrel}, for any $z\in\CC$ 
\begin{equation}
\label{tauz}
\tau_z\ts:\,T_{ij}(u)\ts\mapsto\,T_{ij}(u-z)
\end{equation}
defines an automorphism $\tau_z$ of the algebra $\Y(\gl_n)\ts$. 
Here each of the formal
power series $T_{ij}(u-z)$ in $(u-z)^{-1}$ should be re-expanded 
in $u^{-1}$, and every assignment \eqref{tauz} is a correspondence
between the respective coefficients of series in $u^{-1}$.
Relations \eqref{yrel} also show that for any 
formal power series $g(u)$ in $u^{-1}$ with coefficients from 
$\CC$ and leading term $1$, the assignments
\begin{equation}
\label{fut}
T_{ij}(u)\ts\mapsto\,g(u)\,T_{ij}(u)
\end{equation}
define an automorphism of the algebra $\Y(\gl_n)\ts$. 
Using \eqref{yrel}, one can 
verify that 
\begin{equation}
\label{eval}
T_{ij}(u)\ts\mapsto\,\de_{ij}+E_{ij}\ts u^{-1}
\end{equation}
defines a homomorphism of unital associative algebras 
$\Y(\gl_n)\to\U(\gl_n)\ts$. 

There is an embedding $\U(\gl_n)\to\Y(\gl_n)\ts$,
defined by mapping $E_{ij}\mapsto T_{ij}^{\ts(1)}$. Hence
$\Y(\gl_n)$ contains the universal enveloping
algebra $\U(\gl_n)$ as a subalgebra. The homomorphism \eqref{eval} is
identical on the subalgebra $\U(\gl_n)\subset\Y(\gl_n)\ts$.

Let $T(u)$ be the $n\times n$ matrix 
whose $i\com j$ entry is the series $T_{ij}(u)\ts$. 
The relations \eqref{yrel} can be rewritten
by using the \textit{Yang R-matrix\/}.
This is the $n^2\times n^2$ matrix
\begin{equation}
\label{ru}
R(u)\,=\,u-\sum_{i,j=1}^n\,E_{ij}\ot E_{ji}
\end{equation}
where the tensor factors $E_{ij}$ and $E_{ji}$
are regarded as $n\times n$ matrices. 
Note that
\begin{equation}
\label{rur}
R(u)\,R(-u)=1-u^2\ts.
\end{equation}
Take $n^2\times n^2$ matrices whose entries are
series with coefficients from $\Y(\gl_n)\ts$, 
$$
T_1(u)=T(u)\ot1
\ \quad\text{and}\ \quad
T_2(v)=1\ot T(v)\,.
$$
The collection of relations \eqref{yrel} for all possible indices
$i\com j\com k\com l$ can be written as 
\begin{equation}
\label{rtt}
R(u-v)\,T_1(u)\,T_2(v)\,=\,T_2(v)\,T_1(u)\,R(u-v)\,.
\end{equation}

Using this form of the defining relations together with \eqref{rur}, 
one shows that 
\begin{equation}
\label{tin}
T(u)\mapsto T(-u)^{-1}
\end{equation}
defines an involutive automorphism of the algebra $\Y(\gl_n)\ts$.
Here each entry of the inverse matrix $T(-u)^{-1}$
is a formal power series in $u^{-1}$ with coefficients
from the algebra $\Y(\gl_n)\ts$, and the assignment \eqref{tin} 
is as a correspondence between the respective matrix entries.

The Yangian $\Y(\gl_n)$ is a Hopf algebra over the field $\CC\ts$.
The comultiplication $\De:\Y(\gl_n)\to\Y(\gl_n)\ot\Y(\gl_n)$ is defined by
the assignment
\begin{equation}\label{1.33}
\De:\,T_{ij}(u)\ts\mapsto\ts\sum_{k=1}^n\ T_{ik}(u)\ot T_{kj}(u)\,.
\end{equation}
When taking tensor products of $\Y(\gl_n)\ts$-modules, 
we use the comultiplication \eqref{1.33}.
The counit homomorphism 
$\Y(\gl_n)\to\CC$ is defined by the assignment
$
T_{ij}(u)\ts\mapsto\ts\de_{ij}\ts.
$
The antipodal map $\Y(\gl_n)\to\Y(\gl_n)$ is defined by the assignment
$
T(u)\mapsto T(u)^{-1}.
$
This map is an anti-automorphism of the associative algebra $\Y(\gl_n)\ts$. 
For further details on the Hopf algebra
structure on $\Y(\gl_n)$ see \cite[Chapter 1]{MNO}.

Let $\Tp(u)$ be the transpose to the matrix $T(u)$
relative to the form $\langle\ ,\,\rangle$ on $\CC^{\ts n}\ts$. 
The $i\com j$ entry of the matrix $\Tp(u)$ is
$\ts\th_i\ts\th_j\ts T_{\ts\bj\ts\bi\ts}(u)\ts$. 
Define the $n^2\times n^2$ matrices 
$$
\Tp_1(u)=\Tp(u)\ot1
\ \quad\text{and}\ \quad
\Tp_2(v)=1\ot\Tp(v)\,.
$$
Note that the Yang $R$-matrix \eqref{ru} is invariant under applying 
the transposition relative to $\langle\ ,\,\rangle$ to both tensor factors.
Hence the relation \eqref{rtt} implies that
$$
\Tp_1(u)\,\Tp_2(v)\,R(u-v)
\,=\,
R(u-v)\,\Tp_2(v)\,\Tp_1(u)\,,
$$
\begin{equation}
\label{rttp}
R(u-v)\,\Tp_1(-u)\,\Tp_2(-v)
\,=\,
\Tp_2(-v)\,\Tp_1(-u)\,R(u-v)\,.
\vspace{4pt}
\end{equation}
To obtain the latter relation, we used \eqref{rur}.
By comparing the relations \eqref{rtt} and \eqref{rttp}, 
an involutive automorphism of the algebra $\Y(\gl_n)$
can be defined by the assignment $\ts T(u)\mapsto\Tp(-u)\ts$.
This assignments is
understood as a correspondence between the respective matrix entries.  

Now take the product $\Tp(-u)\,T(u)\ts$.
The $i\com j$ entry of this 
matrix is the series
\begin{equation}
\label{yser}
\sum_{k=1}^n\,\th_i\ts\th_k\,T_{\,\bk\ts\bi\ts}(-u)\,T_{kj}(u)\,.
\end{equation}
The \textit{twisted Yangian\/} 
corresponding to the form $\langle\ ,\,\rangle$
is the subalgebra of $\Y(\gl_n)$ generated by
coefficients of all series \eqref{yser}.
We denote this subalgebra by $\Y(\g_n)\ts$.

To give defining relations for these generators of
$\Y(\g_n)\ts$, let us introduce the \textit{extended twisted Yangian}
$\X(\ts\g_n)\ts$. The unital associative algebra $\X(\ts\g_n)$ has a 
family of generators 
$
S_{ij}^{\ts(1)},S_{ij}^{\ts(2)},\ts\ldots
$
where $i\com j=1\lcd n\ts$. 
Put
\begin{equation*}
S_{ij}(u)=
\de_{ij}+S_{ij}^{\ts(1)}u^{-\ns1}+S_{ij}^{\ts(2)}u^{-\ns2}+\,\ldots
\end{equation*}
and let $S(u)$ be the $n\times n$ matrix 
whose $i\com j$ entry is $S_{ij}(u)\ts$.
Take the $n^2\times n^2$ matrix
\begin{equation}
\label{rup}
\Rp(u)\,=\,u-\sum_{i,j=1}^n\,\th_i\,\th_j\,E_{ij}\ot E_{\ts\bi\ts\bj}
\end{equation}
which is obtained from the Yang $R$-matrix \eqref{ru} by applying 
to any of the two tensor factors
the transposition relative to the form $\langle\ ,\,\rangle$ on $\CC^n\ts$. 
Note the relation
\begin{equation}
\label{rurp}
\Rp(u)\,\Rp(\ts n-u)=u\ts(\ts n-u)\ts.
\end{equation}
Take $n^2\times n^2$ matrices whose entries are
series with coefficients from 
$\X(\g_n)\ts$, 
$$
S_1(u)=S(u)\ot1
\ \quad\text{and}\ \quad
S_2(v)=1\ot S(v)\,.
$$
Defining relations in the algebra $\X(\g_n)$ can be written
as a single matrix relation
\begin{equation}
\label{rsrs}
R(u-v)\,S_1(u)\,\Rp(-u-v)\,S_2(v)\,=\,S_2(v)\,\Rp(-u-v)\,S_1(u)\,R(u-v)\,.
\end{equation}
It is equivalent to the collection of relations
$$
(u^2-v^2)\,[\ts S_{ij}(u)\ts,S_{kl}(v)\ts]\ts=\ts
(u+v)(\ts S_{kj}(u)\,S_{il}(v)-S_{kj}(v)\,S_{il}(u))
$$
$$
\mp\,(u-v)\,(\,
\th_k\ts\th_j\,S_{i\ts\bk}(u)\,S_{\ts\bj\ts l}(v)-
\th_i\ts\th_l\,S_{k\ts\bi}(v)\,S_{\ts\bl\ts j}(u))
$$
\vspace{-8pt}
\begin{equation}
\label{xrel}
\pm\,\ts\th_i\ts\th_j\,
(\ts S_{k\ts\bi\ts}(u)\,S_{\ts\bj\ts l}(v)-
S_{k\ts\bi\ts}(v)\,S_{\ts\bj\ts l}(u))\,.
\vspace{4pt}
\end{equation}
Similarly to \eqref{fut}, this collection of
relations shows that for any 
formal power series $f(u)$ in $u^{-1}$ with the coefficients from 
$\CC$ and leading term $1$, the assignments
\begin{equation}
\label{fus}
S_{ij}(u)\ts\mapsto\,f(u)\,S_{ij}(u)
\end{equation}
define an automorphism of the algebra $\X(\g_n)\ts$. 

\begin{proposition*}
\label{xyp}
One can define a homomorphism $\X(\g_n)\to\Y(\g_n)$ by assigning
\begin{equation}
\label{xy}
S(u)\,\mapsto\,\Tp(-u)\,T(u)\ts.
\end{equation}
\end{proposition*}

\begin{proof}
We have to establish the equality of $n^2\times n^2$ matrices whose 
entries are series with the coefficients from the algebra $\Y(\gl_n)\ts$, 
$$
R(u-v)\,\Tp_1(-u)\,T_1(u)\,\Rp(-u-v)\,\Tp_2(-v)\,T_2(v)
$$
\begin{equation}
\label{1.0.1}
\,=\,\Tp_2(-v)\,T_2(v)\,\Rp(-u-v)\,\Tp_1(-u)\,T_1(u)\,R(u-v)\,.
\end{equation}
But the relations \eqref{rtt},\eqref{rttp} respectively imply
\begin{equation}
\label{1.0.2}
\Tp_1(-u)\,\Rp(-u-v)\,T_2(v)\,=\,T_2(v)\,\Rp(-u-v)\,\Tp_1(-u)\,,
\end{equation}
\begin{equation}
\label{1.0.3}
T_1(u)\,\Rp(-u-v)\,\Tp_2(-v)\,=\,\Tp_2(-v)\,\Rp(-u-v)\,T_1(u)\,.
\vspace{4pt}
\end{equation}
By using 
\eqref{1.0.3},\eqref{rttp},\eqref{rtt}
and \eqref{1.0.2}
we get the relation \eqref{1.0.1}.
Here we also used the commutativity of the
$n^2\times n^2$ matrices \eqref{ru} and \eqref{rup} with each other.
\end{proof}

By definition, the homomorphism \eqref{xy} is surjective. 
Further, the algebra $\X(\g_n)$ has a distinguished family 
of central elements. Indeed,
by dividing each side of the equality \eqref{rsrs}
by $S_2(v)$ on the left and right and 
setting $v=-u\ts$, we get
$$
\Rp(0)\,S_1(u)\,R(2u)\,S_2(-u)^{-1}=\ts
S_2(-u)^{-1}\ts R(2u)\,S_1(u)\,\Rp(0)\,.
$$
The rank of the matrix $\Rp(0)$ equals $1$. So
the last displayed equality implies existence of a formal power series
$O(u)$ in $u^{-1}$ with the coefficients in $\X(\g_n)$ and leading term $1$,
such that
\begin{equation}
\label{5.2}
\Rp(0)\,S_1(u)\,R(2u)\,S_2(-u)^{-1}=\ts
(2u\mp1)\,O(u)\,\Rp(0)\,.
\end{equation}
By \cite[Theorem 6.3]{MNO} all coefficients of the series $O(u)$
belong to the centre of $\X(\g_n)\ts$.
Let us write
$$
O(u)\,=\,1+O^{(1)}u^{\ns-1}+O^{(2)}u^{-\ns2}+\,\ldots\,.
$$
By \cite[Theorem 6.4]{MNO}
the kernel of the homomorphism \eqref{xy}
coincides with the (two-sided) ideal generated by the central elements
$O^{\ts(1)}\ns\com O^{\ts(2)}\ns\com\,\ldots$
defined as coefficients of the series $O(u)\ts$. Using \eqref{rur}, 
one derives from \eqref{5.2} the relation $O(u)\ts O(-u)=1\ts$.

Thus the twisted Yangian $\Y(\g_n)$ can be defined as the 
associative algebra with the generators
$
S_{ij}^{\ts(1)},S_{ij}^{\ts(2)},\ts\ldots
$
which satisfy the relation $O(u)=1$ and the 
\textit{reflection equation\/} \eqref{rsrs}.
This terminology has been used by physicists;
see \cite[Chapter~3]{MNO} for the references,
and for more details on the definition of the algebra $\Y(\g_n)$.
In the present article we need the algebra $\X(\g_n)$ which is
determined by \eqref{rsrs} alone, because this algebra admits
an analogue of the automorphism \eqref{tin} of the Yangian $\Y(\gl_n)\ts$.
Indeed, using \eqref{rsrs} together with the relations 
\eqref{rur} and \eqref{rurp}, one shows that the assignment
\begin{equation}
\label{sin}
\om_n:\ts S(u)\ts\mapsto{S\ts\Bigl(-\ts u-\frac{n}2\,\Bigr)}^{-1}
\end{equation}
defines an involutive automorphism $\om_n$ of $\X(\g_n)\ts$.
However, $\om_n$ does not determine an automorphism of the algebra 
$\Y(\g_n)\ts$, because the map $\om_n$ does
not preserve the ideal of $\X(\g_n)$ generated by the elements
$O^{(1)}\ns\com O^{(2)}\ns\com\,\ldots\ts\,$; 
see \cite[Section 6.6]{MNO}. 
Note that by multiplying
\eqref{5.2} on the right by $S_2(-u)\ts$, the relation $O(u)=1$
can be rewritten as
\begin{equation}
\label{spu}
\Sp(u)\ts=\ts S(-u)\ts\pm\ts\frac{S(u)-S(-u)}{2u}
\end{equation}
where $\Sp(u)$ is the transpose to the matrix $S(u)$
relative to the form $\langle\ ,\,\rangle$ on $\CC^{\ts n}\ts$. 

The definition \eqref{5.2} of the series 
$O(u)$ implies that 
\eqref{fus} determines an automorphism
of the quotient algebra $\Y(\g_n)$ of $\X(\g_n)\ts$, 
if and only if $f(u)=f(-u)$.
If $z\neq0$, the automorphism $\tau_z$ of
$\Y(\gl_n)$ does not preserve the subalgebra 
$\Y(\g_n)$ of $\Y(\gl_n)\ts$. 
There is no analogue of the automorphism $\tau_z$
for the algebra $\X(\g_n)$.

However, there is an analogue of the 
homomorphism $\Y(\gl_n)\to\U(\gl_n)$ defined by \eqref{eval}.
Namely, one can define a homomorphism $\pi_n:\X(\g_n)\to\U(\g_n)$
by the assignments
\begin{equation}
\label{pin}
\pi_n:\,S_{ij}(u)\,\mapsto\,\de_{ij}+
\frac{E_{ij}-\th_i\ts\th_j\ts E_{\ts\bj\ts\bi}}
{\textstyle u\pm\frac12}
\end{equation}
This can be proved using the defining relations \eqref{xrel},
see \cite[Proposition~3.11]{MNO}. 
Furthermore, the central elements
$O^{\ts(1)}\ns\com O^{\ts(2)}\ns\com\,\ldots$ of $\X(\g_n)$
belong to the kernel of $\pi_n\ts$.
Thus $\pi_n$ factors through the homomorphism $\X(\g_n)\to\Y(\g_n)$
defined by \eqref{xy}. 

Further, there is an embedding $\U(\g_n)\to\Y(\g_n)$ defined by mapping
each element $E_{ij}-\,\th_i\ts\th_j\ts E_{\ts\bj\ts\bi}$ of 
the Lie algebra $\g_n$
to the coefficient at $u^{-1}$ of the series \eqref{yser}. Hence
the twisted Yangian $\Y(\g_n)$ contains the universal enveloping
algebra $\U(\g_n)$ as a subalgebra.
The homomorphism $\Y(\g_n)\to\U(\g_n)$ corresponding to $\pi_n$
is evidently identical on the subalgebra $\U(\g_n)\subset\Y(\g_n)\ts$.

Note that the homomorphism $\Y(\g_n)\to\U(\g_n)$ corresponding to $\pi_n$
cannot be obtained by restriction of 
\eqref{eval} to the subalgebra $\Y(\g_n)\subset\Y(\gl_n)\ts$,
because the image of $\Y(\g_n)$ relative to 
the homomorphism 
\eqref{eval} is not contained in the subalgebra $\U(\g_n)\subset\U(\gl_n)\ts$. 
A link between the homomorphisms
\eqref{eval} and \eqref{pin} is given by 
\cite[Lemma 5.4]{N2}.

For any positive integer $l\ts$, consider 
the vector space $\CC^{\ts l}$ and the corresponding Lie algebra
$\gl_{\ts l}\ts$. Let $E_{ab}\in\End(\CC^{\ts l})$ with
$a\com b=1\lcd l$ be the standard matrix units.
When regarding these matrix units as generators of
the universal enveloping algebra $\U(\gl_{\ts l})\ts$, introduce the
$l\times l$ matrix $E$ whose $a\com b$ entry is the
generator $E_{ab}\ts$. Denote by $\Ep$ the 
$l\times l$ matrix whose $a\com b$ entry is the
generator $E_{\ts ba}\ts$. Then consider the matrix inverse
$(\ts u+\Ep\ts)^{-1}\ts$. The $a\com b$ entry 
$(\ts u+\Ep\ts)^{-1}_{\ts ab}$ of the inverse matrix
is a formal power series in $u^{-1}$
with the leading term $\de_{ab}\,u^{-1}$ and the
coefficients from the algebra $\U(\gl_{\ts l})\ts$.

Now take the tensor product of vector spaces
$\CC^{\ts l}$ and $\CC^{\ts n}$. Consider 
the ring $\P\ts(\CC^{\ts l}\ot\CC^{\ts n})$
of polynomial functions on $\CC^{\ts l}\ot\CC^{\ts n}\ts$.
Here have the coordinate functions
$x_{ai}$ where $a=1\lcd l$ and $i=1\lcd n\ts$.
The ring of differential operators on $\P\ts(\CC^{\ts l}\ot\CC^{\ts n})$
with polynomial coefficients will be denoted by 
$\PD\ts(\CC^{\ts l}\ot\CC^{\ts n})\ts$.
Let $\d_{ai}$ be the partial derivation on $\P\ts(\CC^{\ts l}\ot\CC^{\ts n})$
corresponding to $x_{ai}\ts$.

The Lie algebra $\gl_{\ts l}$ acts on the vector space
$\P\ts(\CC^{\ts l}\ot\CC^{\ts n})$ so that the generator
$E_{ab}$ acts as the differential operator
\begin{equation}
\label{glmact}
\sum_{k=1}^n\,
x_{ak}\,\d_{\ts bk}\ts.
\end{equation}
Denote by $\A_{\ts l}$ the tensor product of associative algebras 
$\U(\gl_{\ts l})\ot\PD\ts(\CC^{\ts l}\ot\CC^{\ts n})\ts$. We have an 
embedding $\U(\gl_{\ts l})\to\A_{\ts l}$ defined for $a\com b=1\lcd l$
by the mappings
\begin{equation}
\label{eabact}
E_{ab}\mapsto 
E_{ab}\ot1\,+\,
\sum_{k=1}^n\,
1\ot x_{ak}\,\d_{\ts bk}\ts.
\end{equation}
The following proposition was proved in \cite[Section 1]{KN1},
see also \cite[Section 3]{A}.

\begin{proposition*}
\label{dast}
{\rm\,\,(i)} 
One can define a homomorphism $\al_{\ts l}:\Y(\gl_n)\to\A_{\ts l}$ 
by mapping
\begin{equation}
\label{ehom}
\al_{\ts l}:\ 
T_{ij}(u)\ts\mapsto\,\de_{ij}+
\sum_{a,b=1}^l(\ts u+\Ep\ts)^{-1}_{\ts ab}\ot x_{ai}\,\d_{\ts bj}
\end{equation}
{\rm(ii)} 
The image of\/ $\Y(\gl_n)$ in\/ $\A_{\ts l}$ relative to this 
homomorphism commutes with the image of\/ $\U(\gl_{\ts l})$ in\/ $\A_{\ts l}$
relative to the embedding \eqref{eabact}.
\end{proposition*}

Note that
$$
\al_{\ts l}:\ 
T_{ij}^{\ts(1)}\,\mapsto\,\sum_{c=1}^l\,1\ot x_{ci}\,\d_{cj}\ts.
$$
Hence the restriction of the homomorphism $\al_{\ts l}$
to the subalgebra $\U(\gl_n)\subset\Y(\gl_n)$ corresponds
to the natural action of the Lie algebra $\gl_n$ on
$\P\ts(\CC^{\ts l}\ot\CC^{\ts n})$.

Denote by $Z(u)$ the trace of the inverse matrix
$(u-E)^{-1}\ts$, so that
\begin{equation}
\label{zu}
Z(u)\,=\,\sum_{c=1}^l\,(\ts u-E\ts)^{-1}_{\ts cc}\,.
\end{equation}
Then $Z(u)$ is a formal power series in $u^{-1}$
with the coefficients from the algebra  $\U(\gl_{\ts l})\ts$.
It is well known that these coefficients actually belong to the
centre $\Z\ts(\gl_{\ts l})$ of $\U(\gl_{\ts l})\ts$.
Note that the leading term of this series is $l\ts u^{-1}$.

Let us choose the Borel subalgebra $\b$ of the Lie algebra $\gl_{\ts l}$
spanned by the elements $E_{ab}$ where $a\le b\ts$.
Let $\a\subset\b$ be the Cartan subalgebra of $\gl_{\ts l}$
with the basis $(E_{11}\lcd E_{ll})\ts$. Consider
the corresponding \textit{Harish-Chandra homomorphism\/}
$\varphi_{\ts l}:\ts\U(\gl_{\ts l})^\a\to\U(\a)\ts$. By definition,
for any $\a\ts$-invariant element $X\in\U(\gl_{\ts l})$ the difference
$X-\ts\varphi_{\ts l}\ts(X)$ belongs to the left ideal of
$\U(\gl_{\ts l})$ generated by the elements $E_{ab}$ where $a<b\ts$.
Restriction of the homomorphism $\varphi_{\ts l}$ to
$\Z\ts(\gl_{\ts l})\subset\U(\gl_{\ts l})^\a$ is injective. 
It is well known that 
\begin{equation}
\label{hczu}
1-\varphi_{\ts l}\ts(Z(u))\,=\,\prod_{a=1}^l\,\ts
\Bigl(1-\frac1{u-l+a-E_{aa}}\ts\Bigr)\,,
\end{equation}
see \cite[Theorem 3]{PP}.
The following lemma also involves the series $1-Z(u)\ts$.

\begin{lemma*}
\label{eep}
For any\/ $a\com d=1\lcd l$ we have the equality
$$
(u-E\ts)^{-1}_{\ts da}=(1-Z(u))\,(\ts u-l-\Ep\ts)^{-1}_{\ts ad}\,.
$$
\end{lemma*}

\begin{proof}
First, let us prove for any $a\com b\com c\com d=1\lcd l$
the commutation relation
\begin{equation}
\label{eeu}
[\ts (u-E)_{ab}\,,(u-E\ts)^{-1}_{\ts cd}\ts]=
\de_{ad}\,(u-E\ts)^{-1}_{\ts cb}-
\de_{cb}\,(u-E\ts)^{-1}_{\ts ad}\,.
\end{equation}
Indeed, for any $e\com f=1\lcd l$
by multiplying \eqref{eeu} on the left by $(u-E)_{ec}\ts$,
on the right by $(u-E)_{df}\ts$, and taking the sums over the indices
$c\com d$ we get the relation
$$
[\ts (u-E\ts)_{ef}\,,(u-E)_{ab}\ts]=
\de_{eb}\,(u-E\ts)_{af}-
\de_{af}\,(u-E\ts)_{eb}\,,
$$
which is evidently valid. By setting $b=c$ in \eqref{eeu} and
taking the sum over the index $c$ of the resulting relations, we obtain
the equality
$$
\de_{ad}\,-\,\sum_{c=1}^l\,
(u-E\ts)^{-1}_{\ts cd}\,(u-E)_{ac}\,=\,
\de_{ad}\,Z(u)-l\,(u-E\ts)^{-1}_{\ts ad}\,\,,
\vspace{-8pt}
$$
or
$$
\de_{ad}\,(1-Z(u))\,=\,
\sum_{c=1}^l\,(u-E\ts)^{-1}_{\ts cd}\,(\ts u-l-E)_{\ts ac}\,.
$$
At the right hand side of the latter equality
we have the factor
$$
(\ts u-l-E)_{\ts ac}=(\ts u-l-\Ep\ts)_{\ts ca}\ts.
$$
Thus Lemma \ref{eep} follows from that equality.
\end{proof}

Now let $U$ be a module of the Lie algebra $\gl_{\ts l}\ts$.
Using the homomorphism \eqref{ehom} we can turn the tensor
product of $\gl_{\ts l}$-modules $U\ot\P\ts(\CC^{\ts l}\ot\CC^{\ts n})$
to a bimodule over $\gl_{\ts l}$ and $\Y(\gl_n)\ts$.
This bimodule is denoted by $\E_{\ts l}\ts(U)\ts$.
More generally, for $z\in\CC$ denote by
$\E_{\ts l}^{\ts z}\ts(U)$ the $\Y(\gl_n)\ts$-module obtained
from  $\E_{\ts l}\ts(U)$ via pull-back through the automorphism
$\tau_z$ of $\Y(\gl_n)\ts$, see \eqref{tauz}. It is determined by 
the homomorphism $\Y(\gl_n)\to\A_{\ts l}$ such that
$$
T_{ij}(u)\ts\mapsto\,\de_{ij}\,+
\sum_{a,b=1}^l(\ts u-z+\Ep\ts)^{-1}_{\ts ab}\ot x_{ai}\,\d_{\ts bj}
$$
for any $i\com j=1\lcd n\ts$. 
As a $\gl_{\ts l}\ts$-module $\E_{\ts l}^{\ts z}\ts(U)$
coincides with $\E_{\ts l}\ts(U)$ by definition.
In the next section we will introduce analogues
of the homomorphism \eqref{eabact} 
and of the correspondence $U\mapsto\E_{\ts l}\ts(U)$ for the twisted Yangian
$\Y(\g_n)$ instead of $\Y(\gl_n)\ts$.


\section{Howe duality}
\setcounter{section}{2}
\setcounter{equation}{0}
\setcounter{theorem*}{0}

Take the even-dimensional vector space $\CC^{\ts 2m}$.
Equip $\CC^{\ts 2m}$ with a non-degenerate bilinear form,
alternating in the case $G_n=O_n$ and symmetric in the case $G_n=Sp_n\ts$. 
Thus the tensor product of vector spaces
$\CC^{\ts 2m}$ and $\CC^{\ts n}$
will always come with an alternating bilinear form.
Let $\f_m$ be the subalgebra of the general Lie algebra
$\gl_{\ts 2m}$ preserving our bilinear form on $\CC^{\ts 2m}$. 
We have $\f_m=\sp_{2m}$ or $\f_m=\so_{2m}$ respectively
in the case of an alternating or a symmetric form on $\CC^{\ts 2m}$.
Let us label the standard basis vectors of $\CC^{\ts 2m}$ by the
numbers $-\ts m\lcd-1\com1\lcd m\ts$.
Let $E_{ab}\in\End(\CC^{\ts2m})$ be the standard matrix units,
where the indices $a$ and $b$ run through these numbers.
We will also regard these matrix units as basis elements of 
$\gl_{\ts2m}\ts$. Put 
\begin{equation}
\label{eab}
\ep_{ab}=\operatorname{sign}a\cdot\operatorname{sign}b
\ \quad\text{or}\ \quad
\ep_{ab}=1
\end{equation}
respectively in the case of an alternating or a symmetric form on 
$\CC^{\ts 2m}$.
Then choose the form on $\CC^{\ts2m}$ so that
the Lie subalgebra $\f_m\subset\gl_{\ts2m}$ is spanned
by the elements
\begin{equation}
\label{fab}
F_{ab}=E_{ab}-\ep_{ab}\,E_{-b,-a}\,.
\end{equation}
In the universal enveloping algebra $\U(\ts\f_m)$
we have the commutation relations
\begin{equation}
\label{ufmrel}
[\ts F_{ab}\com F_{cd}\ts]=
\de_{cb}\,F_{ad}-
\de_{ad}\,F_{cb}-
\ep_{ab}\,\de_{c,-a}\,F_{-b,d}+
\ep_{ab}\,\de_{-b,d}\,F_{c,-a}\,.
\end{equation}

Let $F$ be the $2\ts m\times2\ts m$ matrix 
whose $a\com b$ entry is the element $F_{ab}\ts$. 
Denote by $F(u)$ the inverse to the matrix $u-F\ts$. 
Let $F_{ab}(u)$ be the $a\com b$ entry of the inverse matrix.
Any of these entries may be regarded as a formal power series in $u^{-1}$
with the coefficients from the algebra $\U(\ts\f_m)\ts$. Then 
\begin{equation}
\label{fabu}
F_{ab}(u)=\ts\de_{ab}\,u^{-1}\,+\,
\sum_{s=0}^\infty\,
\sum_{\ |c_1|,\ldots,|c_s|=1}^m
F_{ac_1}\ts F_{c_1c_2}\ldots\ts F_{c_{s-1}c_s}\ts F_{c_sb}\,\,
u^{-s-2}\,.
\end{equation}
When $s=0\ts$, the sum over 
$c_1\lcd c_s$ in \eqref{fabu} is understood as 
$F_{ab}\,u^{-2}\ts$.

\begin{lemma*}
\label{1.1}
For\/ $a,b,c,d=-\ts m\lcd-1\com1\lcd m$
we have the relation
$$
[\ts F_{ab}\com F_{cd\ts}(u)\ts]=
\de_{cb}\,F_{ad}(u)-
\de_{ad}\,F_{cb}(u)-
\ep_{ab}\,\de_{c,-a}\,F_{-b,d}(u)+
\ep_{ab}\,\de_{-b,d}\,F_{c,-a}(u)\,.
$$
\end{lemma*}

\begin{proof}
For any $e\com f=-\ts m\lcd-1\com1\lcd m$
by multiplying this relation on the left by $(u-F)_{ec}\ts$,
on the right by $(u-F)_{df}\ts$, and taking the sums over the indices
$c\com d$ we get
$$
[\ts F_{ab}\com F_{ef}\ts]=
\de_{eb}\,F_{af}-
\de_{af}\,F_{eb}-
\ep_{ab}\,\de_{a,-e}\,F_{-b,f}+
\ep_{ab}\,\de_{-b,f}\,F_{e,-a}\,.
$$
This relation can also be obtained by
changing in \eqref{ufmrel} the indices $c\com d\ts$ to $e\com f\ts$.
\end{proof}

This lemma is well known.
Denote by $W(u)$ the trace of the matrix $F(u)\ts$,
\begin{equation}
\label{wu}
W(u)=\sum_{|c|=1}^m F_{cc}(u)\,.
\end{equation}
Lemma \ref{1.1} implies that the coefficients
of the series $W(u)$ belong to the centre $\Z\ts(\f_m)$
of the algebra $\U(\ts\f_m)\ts$. 
Now consider the \textit{triangular decomposition}
of the Lie algebra $\f_m\ts$,
\begin{equation}
\label{tridec}
\f_m=\n\op\h\op\np
\end{equation}
where $\h$ is the Cartan subalgebra of $\f_m$ with the basis
$(\ts F_{-m,-m}\lcd F_{-1,-1})\ts$. Further,
$\n$ and $\np$ are the nilpotent 
subalgebras of $\f_m$ spanned by elements $F_{ab}$ where
$a>b$ and $a<b$ respectively; the indices $a\com b$ here
can be positive or negative. Consider the Harish-Chandra homomorphism
$\psi_m:\ts\U(\ts\f_m)^\h\to\U(\h)$
corresponding to the decomposition \eqref{tridec}.
By definition, for any $\h\ts$-invariant element $Y\in\U(\ts\f_m)$
the difference $Y-\ts\psi_m\ts(Y)$ belongs to the left ideal in 
$\U(\ts\f_m)$ generated by all the elements of $\np\ts$.
Restriction of the homomorphism $\psi_m$ to
$\Z\ts(\f_m)\subset\U(\ts\f_m)^\h$ is injective. 
By \cite[Theorem 3]{PP} we have
$$
1-\Bigl(1\ts\mp\ts\frac1{2u-2m\mp1}\ts\Bigr)^{\ns-1}\psi_m\ts(\ts W(u))\,=
$$
\begin{equation}
\label{hcwu}
\prod_{a=1}^m\,\ts
\Bigl(1-\frac1{u-m+a-F_{aa}}\ts\Bigr)
\Bigl(1-\frac1{m\pm1-u+a-F_{aa}}\ts\Bigr)^{\ns-1}\,.
\end{equation}
Here and in what follows, the upper signs in $\ts\pm\ts$ and $\ts\mp\ts$
correspond to the case of an alternating form on $\CC^{\ts 2m}$  
while the lower signs correspond to the case of a symmetric form
on $\CC^{\ts 2m}$. In these cases 
we have respectively a symmetric or an alternating form 
on $\CC^n\ts$. Thus the choice of signs in
$\ts\pm\ts$ and $\ts\mp\ts$ here agrees with our general convention
on double signs.
Let $\Fp(u)$ be the transpose to the matrix $F(u)$
relative to our 
bilinear form on $\CC^{\ts2m}$ so that
the $a\com b$ entry $\Fp_{ab}(u)$ of the matrix $\Fp(u)$ equals
$\ts\ep_{ab}\,F_{-b,-a}(u)\ts$.

\begin{proposition*}
\label{1.2}
We have the equality of\/ $2m\times2m$ matrices
$$
\Fp(u)\,=\,\bigl(\ts W(u)\ts\pm\ts\frac1{2u-2m\mp1}-1\tts\bigr)
\,F(\tts2m\pm1-u)\,\mp\,\frac{F(u)}{2u-2m\mp1}\,\,.
\vspace{4pt}
$$
\end{proposition*}

\begin{proof}
Let us multiply both sides of the commutation relation in Lemma \ref{1.1}
by $\ep_{a,-d}$ and set $c=-a$ there. 
By taking the sum over the index $a\ts$ of resulting equalities, we get
$$
-\!\sum_{|a|=1}^m
\ep_{\ts b,-d}\,F_{-b,-a}\,F_{-ad}(u)
\,-\!\sum_{|a|=1}^m
\ep_{a,-d}\,F_{-ad}(u)\,F_{ab}\,=
$$
\begin{equation}
\label{1.2.1}
\mp\,\,\ep_{\ts b,-d}\,F_{-b,d}(u)
\ts\pm\ts F_{-db}(u)-
2m\,\ep_{\ts b,-d}\,F_{-b,d}(u)\ts+\ts
\de_{-b,d}\,W(u)\,.
\vspace{10pt}
\end{equation}
Here we used the relation $F_{ab}=-\,\ep_{ab}\,F_{-b,-a}\ts$.
The collection of equalities \eqref{1.2.1} for all indices $b\com d$
can be written as the matrix equality
\begin{equation}
\label{1.2.2}
(\ts-\ts F\,F(u)\tts)^{\ts\prime}-\Fp(u)\,F=
\ts\mp\ts\,\Fp(u)\pm F(u)-2m\ts\Fp(u)+W(u)\,.
\end{equation}
Since $-\ts F\,F(u)=1-u\ts F(u)\ts$,
the equality \eqref{1.2.2} can be rewritten as  
\begin{equation}
\label{1.2.3}
\Fp(u)\cdot(\tts2m\pm1-u-F\tts)=W(u)-1\pm F(u)\,.
\end{equation}
Proposition \ref{1.2} follows from the latter equality
by using the relation
$$
F(u)\,F(\tts2m\pm1-u\tts)
\,=\,
\frac{F(\tts2m\pm1-u\tts)-F(u)}{2u-2m\mp1}\,\,.
\eqno\qedhere
$$
\end{proof}

The next result can be derived either from
Proposition \ref{1.2} or from \eqref{hcwu}.
\begin{corollary*}
\label{1.3}
We have the equality
$$
\bigl(\ts W(u)\ts\pm\ts\frac1{2u-2m\mp1}-1\tts\bigr)
\bigl(\ts W(\tts 2m\pm1-u\tts)\ts\mp\ts\frac1{2u-2m\mp1}-1\tts\bigr)
$$
$$
=1-\frac1{(2u-2m\mp1)^2}\,\,.
$$
\end{corollary*}

\medskip
Now consider the
ring $\P\ts(\CC^{\ts m}\ot\CC^{\ts n})$ of polynomial functions
on $\CC^{\ts m}\ot\CC^{\ts n}\ts$.
Here we have the standard coordinate functions
$x_{ai}$ where $a=1\lcd m$ and $i=1\lcd n\ts$.
Let $\d_{ai}$ be the partial derivation on 
$\P\ts(\CC^{\ts m}\ot\CC^{\ts n})$ corresponding to $x_{ai}\ts$.
There is an action of the Lie algebra $\f_m$
on $\P\ts(\CC^{\ts m}\ot\CC^{\ts n})\ts$
by differential operators,
commuting with the natural action of the group 
$G_n\ts$. The corresponding homomorphism
$\zeta_{\ts n}:\U(\ts\f_m)\to\PD\ts(\CC^{\ts m}\ot\CC^{\ts n})$ 
is defined by the following mappings for $a\com b=1\lcd m\ts$:
$$
\zeta_{\ts n}:\ F_{ab}\,\mapsto\,
\de_{ab}\,\frac{n}2\,\,+\,
\sum_{k=1}^n\,x_{ak}\,\d_{\ts bk}\,,
$$
\begin{equation}
\label{gan}
F_{a,-b}\,\mapsto\,-\,\sum_{k=1}^n\,
\th_k\,x_{a\bk}\,x_{\ts bk}\,,
\ \quad
F_{-a,b}\,\mapsto\,\sum_{k=1}^n\,
\th_k\,\d_{ak}\,\d_{\ts b\bk}\,.
\end{equation}
The homomorphism property here
can be verified by using the relations \eqref{ufmrel}.
Moreover, the image of the homomorphism $\zeta_{\ts n}$ coincides with the
subring of all $G_n$-invariants in $\PD\ts(\CC^{\ts m}\ot\CC^{\ts n})\ts$; 
see \cite[Sections 3.5 and 3.8]{H}.
Let $\B_m$ be the tensor product of associative algebras 
$\U(\ts\f_m)\ot\PD\ts(\CC^{\ts m}\ot\CC^{\ts n})\ts$. Take the 
embedding $\U(\ts\f_m)\to\B_m\ts$ defined by mapping
\begin{equation}
\label{xact}
X\,\mapsto\, 
X\ot1\,+\,
1\ot\zeta_{\ts n}\ts(X)
\ \quad\text{for each}\ \quad
X\in\f_m\,.
\end{equation}

\begin{proposition*}
\label{xb}
{\rm\,\,(i)} 
One can define a homomorphism\/ $\be_m:\X(\g_n)\to\B_m$ 
so that the series $S_{ij}(u)$ is mapped to
the series with coefficients in the algebra $\B_m\ts$,
$$
\de_{ij}\,+
\sum_{a,b=1}^m\,(\,
F_{-a,-b}\ts(u\pm{\textstyle\frac12}+m)
\ot\ts x_{ai}\,\d_{\ts bj}
\ts+\ts
F_{-a,b}\ts(u\pm{\textstyle\frac12}+m)
\ot\ts\th_j\,x_{ai}\,x_{\ts b\bj}
$$
\begin{equation}
\label{fhom}
-\,
F_{a,-b}\ts(u\pm{\textstyle\frac12}+m)
\ot\ts\th_i\,\d_{a\bi}\,\d_{\ts bj}
-\ts
F_{ab}\ts(u\pm{\textstyle\frac12}+m)
\ot\ts\th_i\,\th_j\,\d_{a\bi}\,x_{\ts b\bj}\,)\,.
\vspace{10pt}
\end{equation}
{\rm(ii)} 
The image of\/ $\X(\g_n)$ in\/ $\B_m$ relative to this homomorphism commutes
with the image of\/ $\U(\ts\f_m)$ in\/ $\B_m$ relative to the embedding
\eqref{xact}.
\end{proposition*}

Proposition \ref{xb} can be proved by direct calculation
using the relations \eqref{xrel}. That calculation
is tedious and is omitted here. In Section 6
we will give a more conceptual proof of this proposition.
Now let the indices $c$ and $d$ run through the sequence
$-\ts m\lcd-1\com1\lcd m\ts$. For $c<0$ put 
$p_{\ts ci}=x_{\ts-c,i}$ and $q_{\ts ci}=\d_{\ts-c,i}\ts$.
For $c>0$ put 
$p_{\ts ci}=-\,\th_i\,\d_{c\ts\bi}$
and 
$q_{\ts ci}=\th_i\,x_{c\ts\bi}\,$.
Then our definition of the homomorphism $\be_m$ can be written as
\begin{equation}
\label{fhompq}
\be_m:\,S_{ij}(u)\,\mapsto\,
\de_{ij}\,+
\sum_{|c|,|d|=1}^m\,
F_{cd}\ts(u\pm{\textstyle\frac12}+m)\ot p_{\ts ci}\,q_{\ts dj}\,,
\end{equation}
similarly to \eqref{ehom}. Moreover, then by the definition \eqref{gan}
\begin{equation}
\label{ganpq}
\zeta_{\ts n}:\ F_{cd}\,\mapsto\,
\de_{cd}\,\frac{n}2\,\,-\,
\sum_{k=1}^n\,q_{ck}\,p_{\ts dk}\,.
\end{equation}

Using \eqref{wu}, let us define a formal power series $\Wb(u)$ in $u^{-1}$
with coefficiens in the centre $\Z\ts(\f_m)$ of the algebra $\U(\ts\f_m)$
by the equation
$$
\bigl(\ts1\mp\frac1{2u}\ts\bigr)\,\,\Wb(u)\,=\,
\textstyle
W(u\pm\frac12+m)\,.
$$
By Corollary \ref{1.3},
$$
(\ts\Wb(u)-1\ts))\,(\ts\Wb(-u)-1\ts)=1\,.
$$
Hence there is a formal power series $\Wt(u)$ in $u^{-1}$
with coefficiens in $\Z\ts(\f_m)$ and leading term $1\ts$, such that
\begin{equation}
\label{wtu}
\Wt(-u)\,\ts\Wt(u)^{-1}=1-\Wb(u)\,.
\end{equation}
The series $\Wt(u)$ is not unique. But its coefficient 
at $u^{-1}$ is always $m\ts$, because the leading term
of the series $W(u)$ is $2m\ts u^{-1}$. 
Let $\bet_m$ be the homomorphism
$\X(\g_n)\to\B_m$ defined by assigning to
$S_{ij}(u)$ the series \eqref{fhom} multiplied by 
\begin{equation}
\label{wub}
\Wt(u)\ot1\ts\in\ts\B_m\ts[[u^{-1}]]\ts.
\end{equation}
The homomorphism property of $\bet_m$ follows from Part~(i) of 
Proposition \ref{xb}, see also the defining relations \eqref{xrel}.
Part (ii) implies that
the image of $\bet_m$ commutes with the image of\/ $\U(\ts\f_m)$ in
the algebra\/ $\B_m$ relative to the embedding
\eqref{xact}.

\begin{proposition*}
\label{tb}
The elements
$\,O^{(1)}\ns\com O^{(2)}\ns\com\,\ldots\,$ of\/ $\X(\g_n)$
belong to the kernel of\/ $\bet_m\ts$. 
\end{proposition*}

\begin{proof}
Let us denote by $\St_{ij}(u)$ the product of the 
series \eqref{fhom} and \eqref{wub}. Using the equivalent presentation
\eqref{spu} of the relation $O(u)=1$, we have to prove the equality
\begin{equation}
\label{stu}
\th_i\,\th_j\,\St_{\ts\bj\ts\bi\ts}(u)
\ts=\ts 
\St_{ij}(-u)\ts\pm\ts\frac{\St_{ij}(u)-\St_{ij}(-u)}{2u}
\end{equation}
for any $i\com j=1\lcd n\ts$.
By the definition of the series $\Wt(u)\ts$, we have the relation
\begin{equation}
\label{wtw}
\Wt(u)\,(\ts1-W(u\pm{\textstyle\frac12}+m)\ts)
\ts=\ts
\Wt(-u)\ts\pm\ts\frac{\Wt(u)-\Wt(-u)}{2u}\ .
\end{equation}
Further, let us introduce the $2m\times 2m$ matrix
\begin{equation}
\label{ftud}
\Ft(u)\ts=\ts\Wt(u)\,F(u\pm{\textstyle\frac12}+m)
\end{equation}
and its transpose $\Ftp(u)$ relative to our bilinear form 
on $\CC^{\ts2m}$. By Proposition~\ref{1.2},
\begin{equation}
\label{ftu}
\Ftp(u)
\ts=\ts 
-\ts\Ft(-u)\ts\mp\ts\frac{\Ft(u)-\Ft(-u)}{2u}\ .
\end{equation}

By changing the indices
$i\com j$ in \eqref{fhom} respectively 
to $\bj\ts\com\bi$ and multiplying
the resulting series by $\th_i\,\th_j$ we get
$$
\de_{ij}\,+
\sum_{a,b=1}^m\,(\,
F_{-a,-b}\ts(u\pm{\textstyle\frac12}+m)
\ot\ts\th_i\,\th_j\,x_{a\bj}\,\d_{\ts b\ts\bi}
\,\pm\,F_{-a,b}\ts(u\pm{\textstyle\frac12}+m)
\ot\ts\th_j\,x_{a\bj}\,x_{\ts bi}
$$
$$
\mp\,\,
F_{a,-b}\ts(u\pm{\textstyle\frac12}+m)
\ot\ts\th_i\,\d_{\ts aj}\,\d_{\ts b\ts\bi}
\ts-\ts
F_{ab}\ts(u\pm{\textstyle\frac12}+m)
\ot\ts\d_{\ts aj}\,x_{\ts bi}\,)
\vspace{12pt}
$$
$$
=\ (\ts1-W(u\pm{\textstyle\frac12}+m)\ts)\ot\de_{ij}\ +
$$
$$
\sum_{a,b=1}^m\,(\,
F_{\ts -b,-a}\ts(u\pm{\textstyle\frac12}+m)
\ot\ts\th_i\,\th_j\,\d_{\ts a\bi}\,x_{\ts b\bj}
\,\pm\,F_{\ts -b,a}\ts(u\pm{\textstyle\frac12}+m)
\ot\ts\th_j\,x_{ai}\,x_{\ts b\bj}
$$
$$
\mp\,\,
F_{\ts b,-a}\ts(u\pm{\textstyle\frac12}+m)
\ot\ts\th_i\,\d_{\ts a\bi}\,\d_{\ts bj}
\ts-\ts
F_{\ts ba}\ts(u\pm{\textstyle\frac12}+m)
\ot\ts x_{ai}\,\d_{\ts bj}\,)
\vspace{12pt}
$$
$$
=\ (\ts1-W(u\pm{\textstyle\frac12}+m)\ts)\ot\de_{ij}\ +
$$
$$
\sum_{a,b=1}^m\,(\,
\Fp_{ab}\ts(u\pm{\textstyle\frac12}+m)
\ot\ts\th_i\,\th_j\,\d_{\ts a\bi}\,\ts x_{\ts b\bj}
\ts-\ts
\Fp_{-a,b}\ts(u\pm{\textstyle\frac12}+m)
\ot\ts\th_j\,x_{ai}\,x_{\ts b\bj}\ +
$$
$$
\Fp_{a,-b}\ts(u\pm{\textstyle\frac12}+m)
\ot\ts\th_i\,\d_{\ts a\bi}\,\d_{\ts bj}
\ts-\ts
\Fp_{-a,-b}\ts(u\pm{\textstyle\frac12}+m)
\ot\ts x_{ai}\,\d_{\ts bj}\,)\,.
\vspace{8pt}
$$
Multiplying the expession in the last three lines by
$\Wt(u)\ot1$ and then using the definition \eqref{ftud}, we get
$$
\Wt(u)\,(\ts1-W(u\pm{\textstyle\frac12}+m)\ts)\ot\de_{ij}\ +
\vspace{4pt}
$$
$$
\sum_{a,b=1}^m\,(\,
\Ftp_{ab}\ts(u)
\ot\ts\th_i\,\th_j\,\d_{\ts a\bi}\,\ts x_{\ts b\bj}
\ts-\ts
\Ftp_{-a,b}\ts(u)
\ot\ts\th_j\,x_{ai}\,x_{\ts b\bj}\ +
$$
$$
\Ftp_{a,-b}\ts(u)
\ot\ts\th_i\,\d_{\ts a\bi}\,\d_{\ts bj}
\ts-\ts
\Ftp_{-a,-b}\ts(u)
\ot\ts x_{ai}\,\d_{\ts bj}\,)\,.
\vspace{8pt}
$$
The required equality \eqref{stu} now follows from \eqref{wtw} and \eqref{ftu}.
\end{proof}

Therefore the homomorphism $\bet_m:\X(\g_n)\to\B_m$
factors to a homomorphism $\Y(\g_n)\to\B_m\ts$.
This is an analogue of the homomorphism \eqref{ehom}
for the twisted Yangian $\Y(\g_n)$ instead of $\Y(\gl_n)$.
Note that
$$
\Wt(u)=1+m\ts u^{-1}+\ldots
$$
and so
$$
\bet_m:\,S_{ij}^{\ts(1)}\mapsto\,m\,\de_{ij}\,+\,\sum_{c=1}^m\,\,
(\ts
1\ot\ts x_{ci}\,\d_{\ts cj}
-
1\ot\ts\th_i\,\th_j\,\d_{\ts c\ts\bi}\,x_{c\ts\bj}\ts)
$$
$$
=\,\sum_{c=1}^m\,\,
(\ts
1\ot\ts x_{ci}\,\d_{\ts cj}
-
1\ot\ts\th_i\,\th_j\,x_{c\ts\bj}\,\d_{\ts c\ts\bi}\ts)\,.
$$
Thus for any formal power series $\Wt(u)$ 
in $u^{-1}$ which has its coefficiens in $\Z\ts(\f_m)$, 
has the leading term $1\ts$ and satisfies \eqref{wtu}, 
the restriction of the homomorphism $\Y(\g_n)\to\B_m\ts$
to the subalgebra $\U(\g_n)\subset\Y(\g_n)$ corresponds
to the natural action of the Lie algebra $\g_n$ on
the vector space $\P\ts(\CC^{\ts m}\ot\CC^{\ts n})$.

The series $\Wt(u)$ is not unique, and it will be more convenient
for us to work with the homomorphism $\be_m:\X(\g_n)\to\B_m$
defined in Proposition \ref{xb}. Using this homomorphism
and the action of the Lie algebra $\f_m$ on $\P\ts(\CC^{\ts m}\ot\CC^{\ts n})$
as defined by \eqref{gan}, for arbitrary $\f_m$-module $V$
we can turn the tensor product $V\ot\P\ts(\CC^{\ts m}\ot\CC^{\ts n})$
to a bimodule over $\f_m$ and $\X(\g_n)\ts$.
This bimodule will be denoted by $\F_m(V)\ts$.	

\vspace{1pt}
Consider again the triangular decomposition \eqref{tridec}
of the Lie algebra $\f_m\ts$.
For each $\f_m$-module $V$, denote by $V_{\ts\n}$ the vector space
$V\tts/\ts\n\cdot V$
of coinvariants of the action of the subalgebra $\n\subset\f_m$ on $V$.
The Cartan subalgebra $\h\subset\f_m$ acts on 
the vector space $V_{\ts\n}\ts$.
Now consider
the bimodule $\F_m(V)\ts$. 
The action of $\X(\g_n)$ on this bimodule commutes with the
action of the Lie algebra $\f_m\ts$, 
and hence with the action of the subalgebra $\n\subset\f_m\ts$.
Therefore the space $\F_m(V)_{\ts\n}$
of coinvariants of the action of $\n$ 
is a quotient of the $\X(\g_n)$-module $\F_m(V)\ts$.
Thus we get a functor from the category of all 
$\f_m$-modules to the category of bimodules over $\h$ and $\X(\g_n)\ts$,
\begin{equation}
\label{zelefun}
V\mapsto\ts\F_m(V)_{\ts\n}=
(\ts V\ot\P\ts(\CC^{\ts m}\ot\CC^{\ts n})\ts)_{\ts\n}\ts.
\end{equation}
It may be regarded
as a special case of a general construction
of A.\,Zelevinsky \cite{Z}.

The mappings $E_{ab}\mapsto F_{ab}$ for $a\com b=1\lcd m$
define a Lie algebra embedding $\gl_m\ns\to\f_m\ts$; see relations
\eqref{ufmrel}. Using this embedding, consider the decomposition
\begin{equation}
\label{pardec}
\f_m=\ts\r\ts\op\ts\gl_m\nns\op\rp
\end{equation}
where $\r$ and $\rp$ are the Abelian subalgebras of $\f_m$ 
spanned respectively by the elements $F_{a,-b}$ and
$F_{-a,b}$ for all $a\com b=1\lcd m\ts$.
For any $\gl_m$-module $U$, let $V$ be the $\f_m$-module 
\textit{parabolically induced\/} from the $\gl_m$-module $U$.
To define $V$, one first extends the action of the Lie algebra
$\gl_m$ on $U$ to the maximal parabolic
subalgebra $\gl_m\op\ts\rp\subset\f_m\ts$, so that
every element of the summand $\rp$ acts on $U$ as zero.
By definition, $V$ is the $\f_m$-module induced from the 
$\gl_m\op\ts\rp\ts$-module $U\ts$. Note that here we have
a canonical embedding $U\to V$ of 
$\gl_m\op\ts\rp\ts$-modules\ts; 
we will denote by $\uo$ the image of an element $u\in U$ under this
embedding. The $\f_m\ts$-module $V$ determines
the bimodule $\F_m(V)$ over $\f_m$ and $\X(\g_n)\ts$.
The space of $\r\ts$-coinvariants $\F_m(V)_{\ts\r}$ is
then a bimodule over $\gl_m$ and $\X(\g_n)\ts$.

On the other hand, for any $z\in\CC$ consider 
the bimodule $\E_m^{\ts z}\ts(U)$ over the Lie algebra
$\gl_m$ and over the Yangian $\Y(\gl_n)\ts$. 
By restricting the module $\E_m^{\ts z}\ts(U)$ from the algebra 
$\Y(\gl_n)$ to its subalgebra $\Y(\g_n)$ and then by using the homomorphism
$\X(\g_n)\to\Y(\g_n)$ defined by \eqref{xy}, we can regard
$\E_m^{\ts z}\ts(U)$ as a module over the algebra $\X(\g_n)$
instead of $\Y(\gl_n)\ts$. This module is 
determined by the homomorphism $\X(\g_n)\to\A_m$ such that
for any  $i\com j=1\lcd n$ the series $S_{ij}(u)$ is mapped to 
\begin{equation}
\label{rezhom}
\sum_{k=1}^n\,\th_i\,\th_k\,
\al_m\bigl(\ts T_{\ts\bk\ts\bi}(-\,u-z)\,T_{kj}(u-z)\bigr)\,;
\end{equation}
see \eqref{yser} and \eqref{ehom}.
Let us now map 
$S_{ij}(u)$ to the series \eqref{rezhom} multiplied~by
\begin{equation}
\label{difz}
(1-Z\ts(\ts u+z+m\ts))\ot1\ts\in\ts\A_m\ts[[u^{-1}]]\ts;
\end{equation}
see \eqref{zu} where the positive integer $l$ 
has to be replaced by $m\ts$. 
This mapping determines another homomorphism  
$\X(\g_n)\to\A_m\ts$. Using it, we turn the
vector space $U\ot\P\ts(\CC^{\ts m}\ot\CC^{\ts n})$
of the $\X(\g_n)\ts$-module $\E_m^{\ts z}\ts(U)$ to
another $\X(\g_n)\ts$-module,
to be denoted by $\Eb_m^{\ts z}\ts(U)\ts$.
Further, define an action of the Lie algebra 
$\gl_m$ on $\,\Eb_m^{\ts z}\ts(U)$ by pulling its action on
$\E_m^{\ts z}(U)$ back through the automorphism
\begin{equation}
\label{autom}
E_{ab}\,\mapsto\,\de_{ab}\,\frac{n}2\,+\,E_{ab}
\quad\text{for}\quad
a\com b=1\lcd m\ts.
\end{equation}
The following proposition is a
particular case of Theorem \ref{parind} from the next section.
 
\begin{proposition*}
\label{parres}
For the\/ $\f_m$-module $V$ parabolically induced from any\/
$\gl_m$-module $U$, the bimodule $\F_m(V)_{\ts\r}$
of\/ $\gl_m$ and\/ $\X(\g_n)$ is equivalent to\/ $\Eb_m^{\ts z}\ts(U)$
where $z=\pm\tts\frac12\ts$.
\end{proposition*}

Further, let $u$ and $f$ range over 
$U$ and $\P\ts(\CC^{\ts m}\ot\CC^{\ts n})$ respectively.
In the next section, we will show that the linear map
$$
U\ot\P\ts(\CC^{\ts m}\ot\CC^{\ts n})
\to
(\ts V\ot\P\ts(\CC^{\ts m}\ot\CC^{\ts n}))_{\ts\r}
$$
defined by mapping $u\ot f$ to the class of $\,\uo\ot f$ in 
the space of $\r\ts$-coinvariants, is an equivalence of bimodules
$\Eb_m^{\ts z}\ts(U)\to\F_m(V)_{\ts\r}$ over $\gl_m$ and~$\X(\g_n)$.

An element $\mu$ of the vector space $\h^\ast$ dual to
$\h$ is called a \textit{weight\/}.
A weight $\mu$ can be identified with
the sequence $(\ts\mu_1\lcd\mu_m\ts)$ of its \textit{labels},
where $$
\mu_a=\mu\ts(F_{\ts a-m-1,a-m-1})=-\ts\mu\ts(F_{\ts m-a+1,m-a+1})
\ \quad\textrm{for}\ \quad
a=1\lcd m\ts.
$$
The \textit{Verma module\/} $M_{\ts\mu}$ of the Lie algebra
$\f_m$ is the quotient of the algebra
$\U(\ts\f_m)$ by the left ideal generated by all elements
$X\in\np$, and by all elements $X-\mu\ts(X)$ where $X\in\h\ts$.
The elements of the Lie algebra $\f_m$ act on this quotient
via left multiplication. The image of the identity element $1\in\U(\ts\f_m)$
in this quotient is denoted by $1_\mu\ts$. Then
$X\cdot 1_\mu=0$ for all $X\in\np$, and
$X\cdot 1_\mu=\mu\ts(X)\cdot 1_\mu$
for all $X\in\h\ts$.
Denote by $L_\mu$ be the quotient of the Verma module $M_\mu$
relative to the maximal proper submodule. This quotient is a
simple $\f_m$-module of the highest weight $\mu\ts$. 

For $z\in\CC$ denote by $P_z$
the $\Y(\gl_n)\ts$-module obtained by pulling the standard action
of $\U(\gl_n)$ on $\P(\CC^{\ts n})$ back through 
the homomorphism $\Y(\gl_n)\to\U(\gl_n)$ defined by \eqref{eval}, 
and then back through the automorphism $\tau_z$ of $\Y(\gl_n)$
defined by \eqref{tauz}. Now let $x_1\lcd x_n$ be the coordinates
on $\CC^{\ts n}$ and let $\d_1\lcd\d_n$ be the corresponding
partial derivations. The action of $\Y(\gl_n)$ on $P_z$
is determined by the homomorphism $\Y(\gl_n)\to\PD\ts(\CC^{\ts n})\ts$,
\begin{equation}
\label{pz}
T_{ij}(u)\mapsto\de_{ij}+\frac{\,x_i\ts\d_j}{u-z}\ .
\end{equation}

Using the comultiplication \eqref{1.33},
for any $z_1\lcd z_m\in\CC$ define the tensor product of 
$\Y(\gl_n)\ts$-modules $P_{z_m}\ns\ot\ldots\ot P_{z_1}$. 
For $a=1\lcd m$ let $\deg{\nns}_a$ be the linear operator
on this tensor product, corresponding to evaluation of the polynomial 
degree in the tensor factor $P_{z_a}\ts$; 
that is the $a\ts$-th tensor factor when counting from right to left.
By restricting this tensor
product of $\Y(\gl_n)\ts$-modules to the subalgebra
$\Y(\g_n)\subset\Y(\gl_n)$ and then using the homomorphism
$\X(\g_n)\to\Y(\g_n)$ defined by \eqref{xy}, we can regard
$P_{z_m}\ns\ot\ldots\ot P_{z_1}$ as a module over the extended
twisted Yangian $\X(\g_n)\ts$.

\begin{corollary*}
\label{verma}
The bimodule\/ $\F_m(M_\mu)_{\ts\n}$ over\/ $\h$ and\/ $\X(\g_n)$ is
equivalent to the tensor product
\begin{equation}
\label{promuz}
P_{\mu_m+z}\ot P_{\mu_{m-1}+z+1}\ot\ldots\ot P_{\mu_1+z+m-1}
\vspace{6pt}
\end{equation}
pulled back through the automorphism of\/ $\X(\g_n)$ defined by \eqref{fus},
where $f(u)$ equals
\begin{equation}
\label{fuprod}
\prod_{a=1}^m\,\ts
\Bigl(1-\frac1{u+z+m-a+1+\mu_a}\ts\Bigr)\,;
\end{equation}
here $z=\pm\tts\frac12\ts$. 
The element\/ $F_{\ts m-a+1,m-a+1}\in\h$ acts on 
\eqref{promuz} as the operator
\begin{equation}
\label{faa}
\frac{n}2\ts+\ts\deg{\nns}_a-\mu_{\tts a}\,.
\end{equation}
\end{corollary*}

\begin{proof}
We have an embedding of $\gl_m$ to $\f_m$ such that
$E_{aa}\mapsto F_{aa}$ for any index $a=1\lcd m\ts$.
Then the Cartan subalgebra $\a$ of $\gl_m$ is identified
with the Cartan subalgebra $\h$ of $\f_m\ts$. Put
$\ab=m-a+1$ for short. If we regard the 
weight $\mu$ as an element of $\a^\ast$, then
$$
\mu\ts(E_{\ts\ab\ts\ab\ts})=-\,\mu_{\tts a}
\quad\text{for}\quad
a=1\lcd m\ts.
$$
Let $U$ be the Verma module of the Lie algebra
$\gl_m$ corresponding to $\mu\in\a^\ast\ts$.
It is defined as the quotient of the algebra
$\U(\gl_m)$ by the left ideal generated by all elements
$E_{ab}$ with $a<b\ts$, 
and by all elements $E_{aa}-\mu\ts(E_{aa})\ts$.
The Verma module $M_\mu$ of the Lie algebra $\f_m$ is then
equivalent to the module $V$ parabolically induced from
the $\gl_m$-module $U$. 
Here we use the decomposition \eqref{pardec}.

Let $\s$ denote the subalgebra of the Lie algebra $\gl_m$
spanned by all elements $E_{ab}$ with $a>b\,$.
Using our embedding of $\gl_m$ to $\f_m\ts$, we can also
regard $\s$ as a subalgebra of $\f_m\ts$.
The Lie algebra $\n$ of $\f_m$ is then spanned by
$\r$ and $\s\ts$. By Proposition \ref{parres}, 
the bimodule $\F_m(M_\mu)_{\ts\n}$ over $\h$ and $\X(\g_n)$
is equivalent to $\Eb_m^{\ts z}\ts(U)_{\ts\s}$ where
$z=\pm\tts\frac12\ts$.
To describe the latter bimodule, let us firstly
consider the bimodule $\E_m^{\ts z}\ts(U)_{\ts\s}$
over $\a$ and $\Y(\gl_n)\ts$. 
Using \cite[Corollary 2.4]{KN1}, the bimodule
$\E_m^{\ts z}\ts(U)_{\ts\s}$ is equivalent to the
tensor product of $\Y(\gl_n)\ts$-modules \eqref{promuz}
where the element $E_{\ts\ab\ts\ab\ts}\in\a$ acts as 
$\,\deg{\nns}_a-\mu_{\tts a}\ts$. After pulling the action
of Lie algebra $\gl_m$ on $\E_m^{\ts z}(U)$ back 
through the automorphism \eqref{autom},
the element $E_{\ts\ab\ab\ts}\in\a$ will act on the tensor product
of vector spaces \eqref{promuz} as \eqref{faa}. 

To complete the proof of Corollary \ref{verma},
recall that the action of $\X(\gl_n)$ on 
$\Eb_m^{\ts z}\ts(U)$ differs from that on 
$\E_m^{\ts z}\ts(U)$ by multiplying the series \eqref{rezhom} by 
\eqref{difz}. Using \eqref{hczu}, the series
$1-Z\ts(\ts u+z+m\ts)$ in $u^{-1}$ with the coefficients in $\Z(\gl_m)$
acts on the Verma module $U$ via scalar multiplication by 
the series \eqref{fuprod}.
\end{proof}

By definition, the vector spaces of two equivalent
bimodules in Corollary~\ref{verma} are 
$(M_\mu\ot\P\ts(\CC^{\ts m}\ot\CC^{\ts n}))_{\ts\n}$
and $\P\ts(\CC^{\ts n})^{\ts\ot m}$ respectively.
We can determine a linear map from the latter vector space
to the former one,
by mapping $f_1\ot\ldots\ot f_m$ to
the class of the element $1_\mu\ot f$ in
the space of $\n\ts$-coinvariants. Here
for any $f_1\lcd f_m\in\P\ts(\CC^{\ts n})$ the polynomial
$f\in\P\ts(\CC^{\ts m}\ot\CC^{\ts n})$ is defined by setting
\begin{equation}
\label{fff}
f(x_{11}\lcd 
x_{mn})\ts=
f_1(x_{11}\lcd x_{1n})\ts\ldots\ts f_m(x_{m1}\lcd x_{mn})\ts.
\end{equation}
This provides the bimodule equivalence in
Corollary \ref{verma},
see \cite[Corollary 2.4]{KN1}
and the remarks made immediately after stating Proposition \ref{parres} here.

For any $z\in\CC$ denote by $\Pp_z$ 
the $\Y(\gl_n)\ts$-module obtained by pulling $P_z$
back through the automorphism $T(u)\mapsto\Tp(-u)$ of $\Y(\gl_n)\ts$. 
According to \eqref{pz}, the action of $\Y(\gl_n)$ on $\Pp_z$
is determined by the homomorphism $\Y(\gl_n)\to\PD\ts(\CC^{\ts n})\ts$,
\begin{equation}
\label{ppz}
T_{ij}(u)\mapsto\de_{ij}-
\frac{\,\th_i\ts\th_j\ts x_{\ts\bj}\,\d_{\ts\bi}}{u+z}\ .
\end{equation}

\begin{lemma*}
\label{ppl}
The\/ $\Y(\gl_n)\ts$-module $\Pp_z$
can also be obtained by pushing the action of\/
$\Y(\gl_n)$ on $P_{\ts-z-1}$ forward through the automorphism 
of\/ $\PD\ts(\CC^{\ts n})$ such that for each $i=1\lcd n$ 
\begin{equation}
\label{onefour}
x_i\mapsto-\,\th_i\,\d_{\ts\bi}
\quad\text{and}\quad
\d_i\mapsto\th_i\,x_{\ts\bi}\,,
\end{equation}
and by pulling the resulting action back through the automorphism
\eqref{fut} of\/ $\Y(\gl_n)$ where 
\begin{equation}
\label{fuz}
g(u)=1+\frac1{u+z}\ .
\end{equation}
\end{lemma*}

\begin{proof} 
By applying the automorphism \eqref{onefour} to the right hand
side of \eqref{pz} and replacing the parameter $z$ there by $-z-1$ we get
$$
\de_{ij}-\frac{\,\th_i\ts\th_j\ts\d_{\ts\bi}\,x_{\ts\bj}}{u+z+1}
\,=\,
\de_{ij}-\frac{\,\de_{ij}+\th_i\ts\th_j\ts x_{\ts\bj}\,\d_{\ts\bi}}{u+z+1}
\,=\,
\frac{u+z}{u+z+1}
\,\Bigl(\de_{ij}-
\frac{\,\th_i\ts\th_j\ts x_{\ts\bj}\,\d_{\ts\bi}}{u+z}\ts\Bigr)
$$
which, after multiplying it by \eqref{fuz},
becomes the right hand side of \eqref{ppz}.
\end{proof}

We finish this section with 
the another simple lemma\ts; it will be used several times in
Sections 3 and 6. This lemma is well known,
and applies to matrices over an arbitrary unital ring.
Take any matrix of the size $(p+q)\times(p+q)$ over such a ring,  
and write it as
\begin{equation}
\label{blockmat}
\begin{bmatrix}\,A\,&B\,\\\,C\,&D\,\end{bmatrix}
\end{equation}
where 
$A,B,C,D$ are matrices of sizes
$p\times p$, $p\times q$, $q\times p$, $q\times q$ respectively. 
Here we assume that both $p$ and $q$ are any positive integers.

\begin{lemma*}
\label{lemma0}
Suppose that the matrix \eqref{blockmat} is invertible,
and that the matrices $A$ and $D$ are also invertible. Then the matrices
$A-B\ts D^{-1}\ts C$ and $D-C\ts A^{-1}B$ are invertible too, and
$$
\begin{bmatrix}A&B\\C&D\end{bmatrix}^{-1}\!=\ \  
\begin{bmatrix}
(A-B\ts D^{-1}\ts C)^{-1}
&
\,-A^{-1}B\ts(\ts D-C\ts A^{-1}B)^{-1}\,
\\
\,-D^{-1}C\ts(A-B\ts D^{-1}\ts C)^{-1}\,
&
(\ts D-C\ts A^{-1}B)^{-1}
\end{bmatrix}
\,.
$$
\end{lemma*}


\section{Parabolic induction}
\setcounter{section}{3}
\setcounter{equation}{0}
\setcounter{theorem*}{0}

The twisted Yangian $\Y(\g_n)$ is not just a subalgebra of
$\Y(\gl_n)\ts$, it is also a right coideal of the coalgebra 
$\Y(\gl_n)$ relative to the comultiplication \eqref{1.33}. Indeed,
let us apply this comultiplication to the 
$i\com j$ entry of the $n\times n$ matrix $\Tp(-u)\,T(u)\ts$.
We get the sum
$$
\sum_{k=1}^n\,\th_i\ts\th_k\,
\De\ts(\,T_{\,\bk\ts\bi\ts}(-u)\,T_{kj}(u))\,=
\vspace{-2pt}
$$
$$
\sum_{g,h,k=1}^n\,\th_i\ts\th_j\,
(\,T_{\ts\bk\ts\bg\ts}(-u)\ot T_{\ts\bg\ts\bi\ts}(-u))\,
(\,T_{kh}(u)\ot T_{hj}(u))\,=
$$
$$
\sum_{g,h,k=1}^n\,\th_g\ts\th_k\,
T_{\ts\bk\ts\bg\ts}(-u)\,T_{kh}(u)
\ts\ot\ts
\th_i\ts\th_g\,T_{\ts\bg\ts\bi\ts}(-u)\,T_{hj}(u)\ts.\ 
\vspace{2pt}
$$
In the last displayed line, by performing the summation
over $k=1\lcd n$ in the first tensor factor we get
the $g\com h$ entry of the matrix $\Tp(-u)\,T(u)\ts$. Therefore
$$
\De\ts(\ts\Y(\g_n))\subset\Y(\g_n)\ot\Y(\gl_n)\ts.
$$

Moreover, 
one can define a homomorphism of associative algebras 
$$
\X(\g_n)\to\X(\g_n)\ot\Y(\gl_n)
$$
by assigning
\begin{equation}
\label{comod}
S_{ij}(u)\,\mapsto
\sum_{g,h=1}^n\,S_{gh}(u)
\ts\ot\ts
\th_i\ts\th_g\,T_{\ts\bg\ts\bi\ts}(-u)\,T_{hj}(u)\ts.
\vspace{2pt}
\end{equation}
The homomorphism property can be verified directly,
by consecutively using the relations
\eqref{1.0.3},\eqref{rttp},\eqref{rsrs},\eqref{rtt} and \eqref{1.0.2}.
He we omit further details of this verification, see our proof 
of Proposition \ref{xyp}.
By using the homomorphism \eqref{comod}, the tensor product
of any modules over the algebras $\X(\g_n)$ and $\Y(\gl_n)$
becomes another module over $\X(\g_n)\ts$.

The homomorphism \eqref{comod} is a \textit{coaction}
of the Hopf algebra $\Y(\gl_n)$ on the algebra $\X(\g_n)\ts$.
Formally, one can define a homomorphism of associative algebras 
$$
\X(\g_n)\to\X(\g_n)\ot\Y(\gl_n)\ot\Y(\gl_n)
$$
in two different ways: either by using the assignment \eqref{comod} twice,
or by using \eqref{comod} and then \eqref{1.33}. Let us show that both ways
lead to the same result. Using the first,
\begin{equation}
\label{cone}
S_{ij}(u)\,\mapsto
\sum_{g,h,p,q=1}^n\,S_{pq}(u)
\ts\ot\ts
\th_g\ts\th_p\,T_{\ts\bp\ts\bg\ts}(-u)\,T_{qh}(u)
\ts\ot\ts
\th_i\ts\th_g\,T_{\ts\bg\ts\bi\ts}(-u)\,T_{hj}(u)\ts.
\vspace{2pt}
\end{equation}
Using the second definition,
\begin{equation}
\label{ctwo}
S_{ij}(u)\,\mapsto
\sum_{g,h,p,q=1}^n\,S_{gh}(u)
\ts\ot\ts
\th_i\ts\th_g\,T_{\ts\bg\ts\bp\ts}(-u)\,T_{hq}(u)
\ts\ot\ts
T_{\ts\bp\ts\bi\ts}(-u)\,T_{qj}(u)\ts.
\vspace{2pt}
\end{equation}
But the right hand side of \eqref{ctwo} can be obtained
from the right hand side of \eqref{cone} by exchanging the pairs
of running indices $g\com h$ and $p\com q\ts$. 

Now for any positive integer $l$ consider the general linear Lie algebra
$\gl_{\ts2m+2l}$ and its subalgebra $\f_{m+l}\ts$. This subalgebra
is spanned by the elements $F_{ab}$ where
\begin{equation}
\label{mliml}
a\com b=-\ts m-l\lcd-1\com1\lcd m+l\ts.
\end{equation}
We will extend the notation \eqref{eab} and \eqref{fab}
to all these indices $a\com b\ts$.
Let us now identify $\f_m$ with the subalgebra
of $\f_{m+l}$ spanned by the elements $F_{ab}$ where
$a\com b=-\ts m\lcd-1\com1\lcd m\ts$. Further, let us choose
the embedding of the Lie algebra $\gl_{\ts l}$ to $\f_{m+l}$
determined by the mappings
\begin{equation}
\label{gll}
E_{ab}\mapsto F_{m+a,m+b}
\ \quad\text{for}\ \quad
a\com b=1\lcd l\ts.
\end{equation}
Let $\q\com\qp$ be the subalgebras of $\f_{m+l}$
spanned respectively by the elements $F_{\ts ab}\com F_{\ts ba}$ where 
$$
a=m+1\lcd m+l
\ \quad\text{and}\ \quad
b=-\ts m-l\lcd-1\com1\lcd m\,;
$$
these two subalgebras of $\f_{m+l}$ are nilpotent.
Put $\p=\f_m\op\ts\gl_{\ts l}\op\qp\ts$.
Then $\p$ is a maximal parabolic subalgebra of the reductive 
Lie algebra $\f_{m+l}\ts$, and
$
\f_{m+l}=\q\op\p\ts.
$
We do not exclude the case $m=0$ here.
In this case the nilpotent subalgebras 
$\q$ and $\qp$ of $\f_{m+l}$ become the Abelian subalgebras 
$\r$ and $\rp$ of the Lie algebra $\f_{\ts l}\ts$; 
see the decomposition \eqref{pardec}
where the positive integer $m$ is to be replaced by $l\ts$.
Note that the meaning of the symbols $\p$ and $\q$ here is
different from that in Section 0.

Let $V$ and $U$ be any modules of the Lie algebras 
$\f_m$ and $\gl_{\ts l\ts}$ respectively.
Denote by $V\ns\bt U$ the $\f_{m+l}\ts$-module 
\textit{parabolically induced\/}
from the $\f_m\op\ts\gl_{\ts l}\ts$-module $V\ns\ot U$.
To define $V\ns\bt U$, one first extends the action of the Lie algebra
$\f_m\op\ts\gl_{\ts l}$ on $V\ns\ot U$ to the Lie algebra $\p\ts$, so that
every element of the subalgebra $\qp\subset\p$ acts on $V\ns\ot U$ as zero.
By definition, $V\ns\bt U$ is the $\f_{m+l}\ts$-module induced from the 
$\p\ts$-module $V\ns\ot U\ts$.  
Note that here we have
a canonical embedding $V\ot U\to V\ns\bt U$ of 
$\p\ts$-modules\ts; 
we will denote by $\vuo$ the image of an element $v\ot u\in V\ot U$ 
under this embedding. 

Consider the bimodule $\F_{m+l}\ts(\ts V\ns\bt U\ts)$
over $\f_{m+l}$ and $\X(\g_n)\ts$.
Here the action of $\,\X(\gl_n)$ commutes with the
action of the Lie algebra $\f_{m+l}\ts$, 
and hence with the action of the subalgebra $\q\subset\f_{m+l}\ts$.
Therefore the vector space 
$\F_{m+l}\ts(\ts V\ns\bt U\ts)_{\ts\q}$
of coinvariants of the action of the subalgebra $\q$ 
is a quotient of the $\X(\g_n)$-module 
$\F_{m+l}\ts(\ts V\ns\bt U\ts)\ts$.
Note that the subalgebra 
$\f_m\op\ts\gl_{\ts l}\subset\f_{m+l}$ also acts on this quotient.


For any $z\in\CC$ consider the bimodule $\E_{\ts l}^{\ts z}\ts(U)$
over $\gl_{\ts l}$ and $\Y(\gl_n)$ defined as in the end of Section~1.
Also consider the bimodule $\F_m(V)$ over $\f_m$ and $\X(\g_n)\ts$.
Using the homomorphism \eqref{comod}, the tensor product
of vector spaces $\F_m(V)\ot\ts\E_{\ts l}^{\ts z}\ts(U)$
becomes a module over $\X(\g_n)\ts$. This module is determined
by the homomorphism $\X(\g_n)\to\B_m\ot\A_l$ such that
for any $i\com j=1\lcd n$ the series $S_{ij}(u)$ is mapped to
\begin{equation}
\label{baser}
\sum_{g,h=1}^n\,\be_m\bigl(\ts S_{gh}(u)\bigr)
\ts\ot\ts
\th_i\ts\th_g\,
\al_{\ts l\ts}\bigr(\ts T_{\ts\bg\ts\bi\ts}(-\,u-z)\,T_{hj}(u-z)\bigl)\ts.
\end{equation}
Let us now map the series $S_{ij}(u)$ to
the series \eqref{baser} multiplied by
\begin{equation}
\label{zuser}
\bigl(\ts1\ot1\tts\bigr)
\ot
\bigl(\ts(1-Z\ts(\ts u+z+l\ts))\ot1\ts\bigr)
\ts\in\ts\B_m\ot\A_l\,[[u^{-1}]]\ts,
\end{equation}
see \eqref{zu}.
This mapping determines another homomorphism  
$\X(\g_n)\to\B_m\ot\A_l\ts$. 
Using it, we turn the vector space
of the $\X(\g_n)\ts$-module $\F_m(V)\ot\ts\E_{\ts l}^{\ts z}\ts(U)$
to yet another $\X(\g_n)\ts$-module, which will be denoted by 
$\F_m(V)\ts\,\widetilde{\ot}\,\ts\E_{\ts l}^{\ts z}\ts(U)\,$.
Further, define an action of the Lie algebra $\gl_{\ts l}$ on
the latter $\X(\g_n)\ts$-module by pulling the action of $\gl_{\ts l}$
on $\E_{\ts l}\ts(U)$ back through the automorphism
\begin{equation}
\label{autol}
E_{ab}\,\mapsto\,\de_{ab}\,\frac{n}2\,+\,E_{ab}
\quad\text{for}\quad
a\com b=1\lcd l\ts.
\end{equation}
The Lie algebra $\f_m$ acts on the $\X(\g_n)\ts$-module
$\F_m(V)\ts\,\widetilde{\ot}\,\ts\E_{\ts l}^{\ts z}\ts(U)$
via the tensor factor $\F_m(V)\ts$. Thus 
$\F_m(V)\ts\,\widetilde{\ot}\,\ts\E_{\ts l}^{\ts z}\ts(U)$
becomes a bimodule
over the direct sum of Lie algebras $\f_m\op\ts\gl_{\ts l}$
and over the extended twisted Yangian $\X(\g_n)\ts$.
The next theorem is the most difficult result of our article.
In the case $m=0$ this theorem becomes Proposition~\ref{parres}, 
where the positive integer $m$ has to be replaced by $l\ts$.
Here we assume that $\F_{\ts0}(V)=\CC$ so that
$\be_{\ts0}(S_{ij}(u))=\de_{ij}\ts$.

\begin{theorem*}
\label{parind}
The bimodule $\F_{m+l}\ts(\ts V\ns\bt U\ts)_{\ts\q}$
over\/ $\f_m\op\ts\gl_{\ts l}$ and\/ $\X(\g_n)$ is 
equivalent to
$\F_m(V)\,\ts\widetilde{\ot}\,\ts\E_{\ts l}^{\ts z}\ts(U)$
where $z=m\pm\tts\frac12\ts$.
\end{theorem*}

In the remainder of this section we shall prove Theorem \ref{parind}.
As vector spaces, 
\begin{gather*}
\F_{m+l}\ts(\ts V\ns\bt U\ts)_{\ts\q}=
(\ts V\ns\bt U\ot\P\ts(\CC^{\ts m+l}\ot\CC^{\ts n}))_{\ts\q}\,,
\\
\F_m(V)\,\ts\widetilde{\ot}\,\ts\E_{\ts l}^{\ts z}\ts(U)=
V\ot\P\ts(\CC^{\ts m}\ot\CC^{\ts n})\ot U\ot\P\ts(\CC^{\ts l}\ot\CC^{\ts n})
\ts.
\end{gather*}
We can determine a linear map from the latter vector space to the 
former one by mapping any element  
$v\ot f\ot u\ot g$
to the class of $\vuo\ot f\ot g$ in the space of $\q\ts$-coinvarians.
Here $v\in V$, $f\in\P\ts(\CC^{\ts m}\ot\CC^{\ts n})$ and
$u\in U$, $g\in\P\ts(\CC^{\ts l}\ot\CC^{\ts n})$ whereas 
the tensor product $f\ot g$
is identified with an element of $\P\ts(\CC^{\ts m+l}\ot\CC^{\ts n})$
in a natural way, which corresponds to the decomposition
\begin{equation}
\label{mnl}
\CC^{\ts m+l}\ot\CC^{\ts n}=
\CC^{\ts m}\ot\CC^{\ts n}\op\ts\CC^{\ts l}\ot\CC^{\ts n}\ts.
\end{equation}
We will show
that this map establishes an equivalence of bimodules
in Theorem~\ref{parind}.

The 
vector space of the $\f_{m+l}\ts$-module $V\ns\bt U$
can be identified with the tensor product $\U(\q)\ot V\ns\ot U$
where the Lie subalgebra $\q\subset\f_{m+l}$ acts via
left multiplication on the first tensor factor.
Then $\vuo=1\ot v\ot u\ts$, so that
the tensor product $V\ns\ot U$ gets identified with the subspace
\begin{equation}
\label{onesub}
1\ot V\ns\ot U\subset\ts\U(\q)\ot V\ns\ot U\ts.
\end{equation}
On this subspace,
every element of the subalgebra $\qp\subset\f_{m+l}$ acts as zero, 
while the two direct summands
of subalgebra $\f_m\op\ts\gl_{\ts l}\subset\f_{m+l}$
act non-trivially only on the tensor factors $V$ and $U$ respectively.
All this determines the action of Lie algebra $\f_{m+l}$ on
$\U(\q)\ot V\ns\ot U$. Now consider $\F_{m+l}\ts(\ts V\ns\bt U\ts)$
as a $\f_{m+l}\ts$-module, we will denote it by $M$ for short.
Then $M$ is the tensor product of two $\f_{m+l}\ts$-modules,
\begin{equation}
\label{mmod}
M=(\ts V\ns\bt U)\ot\P\ts(\CC^{\ts m+l}\ot\CC^{\ts n})=
\U(\q)\ot V\ns\ot U\ot\P\ts(\CC^{\ts m+l}\ot\CC^{\ts n})\ts.
\end{equation}

The vector spaces of the $\X(\g_n)\ts$-module $\F_m\ts(V)$ and of the 
$\Y(\gl_n)\ts$-module $\E_{\ts l}^{\ts z}\ts(U)$ are respectively
$V\ns\ot\ts\P\,(\CC^{\ts m}\ot\CC^{\ts n})$ and
$U\ns\ot\ts\P\,(\CC^{\ts l}\ot\CC^{\ts n})\ts$.
The action of the Lie algebra $\f_m$ on the first
vector space is defined by \eqref{gan}. By pulling back through the
automorphism \eqref{autol}, the action 
of the Lie algebra $\gl_{\ts l}$ on the second
vector space is defined by mapping
$$
E_{ab}\,\mapsto\,\de_{ab}\,\frac{n}2\,+\,
E_{ab}\ot1\,+\,
\sum_{k=1}^n\,
1\ot x_{ak}\,\d_{\ts bk}
\ \quad\text{for}\ \quad
a\com b=1\lcd l\ts.
$$
Identify the tensor product of these two vector spaces with
\begin{equation}
\label{vuprod}
V\ns\ot U\ot
\P\,(\CC^{\ts m}\ot\CC^{\ts n})\ot\P\,(\CC^{\ts l}\ot\CC^{\ts n})=
V\ns\ot U\ot
\P\,(\CC^{\ts m+l}\ot\CC^{\ts n})
\end{equation}
where we use the direct sum decomposition \eqref{mnl}.
We get an action of 
the direct sum of Lie algebras $\f_m\op\ts\gl_{\ts l}$ on 
the vector space \eqref{vuprod}.

Let us now define a linear map
$$
\ph:\,V\ns\ot U\ot\P\,(\CC^{\ts m+l}\ot\CC^{\ts n})\ts\to M\ts/\,\q\cdot M
$$
by the assignment
$$
\ph:\,y\ot x\ot t\,\mapsto\,1\ot y\ot x\ot t\,+\,\q\cdot M
$$
for any vectors $y\in V$, $x\in U$ and
$t\in\P\,(\CC^{\ts m+l}\ot\CC^{\ts n})\ts$.
The operator $\ph$ intertwines the actions of the Lie
algebra $\f_m\op\ts\gl_{\ts l}\ts$; see the definition \eqref{gan} 
where $m$ is to be replaced by $m+l\ts$. 
Let us demonstrate that the operator $\ph$ is bijective.

Firstly consider the
action of the Lie subalgebra $\q\subset\f_{m+l}$ on the vector space
$$
\P\ts(\CC^{\ts m+l})\ts=\ts
\P\ts(\CC^{\ts m})\ot\P\ts(\CC^{\ts l})\ts;
$$
the action is defined by \eqref{gan} where $n=1$, 
and the integer $m$ is replaced by $m+l\ts$.
This vector space admits a descending filtration by the subspaces
$$
\mathop{\op}\limits_{K=N}^{\infty}\,
\P\,(\CC^{\ts m})\ot
\P^{\ts K}(\CC^{\ts l})
\ \quad\textrm{where}\ \quad
N=0,1,2,\ts\ldots\,.
$$
Here $\P^{\ts K}(\CC^{\ts l})$ is
the space of polynomial functions on $\CC^{\ts l}$ of degree $K\ts$.
The action of the Lie algebra $\q$ on $\P\ts(\CC^{\ts m+l})$
preserves each of these subspaces, and becomes trivial
on the graded space associated with this filtration. 

Similarly, for any $n=1\com2\com\ts\ldots$
the vector space $\P\,(\CC^{\ts m+l}\ot\CC^{\ts n})$ admits a descending
filtration by $\q\ts$-submodules such that $\q$ acts trivially on
each of the corresponding graded subspaces.
The latter filtration induces a filtration of $M$
by $\q\ts$-submodules such that on
the corresponding graded quotient
$\operatorname{gr}M\ts$, the Lie algebra $\q$
acts via left multiplication
on the first tensor factor $\U(\q)$ in \eqref{mmod}.
The space $V\ns\ot U\ot\P\,(\CC^{\ts m+l}\ot\CC^{\ts n})$
is therefore isomorphic to the space
of coinvariants $(\ts\operatorname{gr}M)_{\ts\q}$ 
via the bijective linear map
$$
y\ot x\ot t\,\mapsto\,
1\ot y\ot x\ot t\,+\,\q\cdot(\ts\operatorname{gr}M\ts)\ts.
$$
Therefore the linear map $\ph$ is bijective as well.

It remains to show that the map $\ph$ intertwines the actions of
the algebra $\X(\g_n)\ts$. 
The definition of the $\X(\g_n)\ts$-module 
$M$ involves
the $(2m+2l)\times(2m+2l)$ matrix whose
$a\com b$ entry is $\de_{ab}\,v-F_{ab}$ where
\begin{equation}
\label{v}
\textstyle
v=u\pm\frac12+m+l\ts.
\end{equation}
Label the rows and columns of this matrix by the indices
\eqref{mliml}. Write this matrix and its inverse as
\begin{equation}
\label{abcd}
\begin{bmatrix}
\,A\,&B\,\\\,C\,&D\,
\end{bmatrix}
\quad\textrm{and}\quad
\begin{bmatrix}\ \At\,&\Bt\ \\ \ \Ct\,&\Dt\ 
\end{bmatrix}
\end{equation}
where the blocks $A,B,C,D$ and $\At,\Bt,\Ct,\Dt$ are matrices of sizes
$(\ts l+2m)\times(\ts l+2m)$, $(\ts l+2m)\times l$, $l\times(\ts l+2m)$,
$l\times l$ respectively. We will label the rows and columns of all these
blocks by the same indices as in the compound matrices \eqref{abcd}. For
instance, the rows and columns of the $l\times l$ matrix
$D$ will be labelled by the indices $m+1\lcd m+l\ts$.
Observe that $D=v-E$
where $E$ denotes the $l\times l$ matrix whose $a\com b$ entry is $F_{ab}$ 
for $a\com b=m+1\lcd m+l\ts$.
Using our embedding \eqref{gll}
of the Lie algebra $\gl_{\ts l}$ to $\f_{m+l}\ts$,
this notation agrees with the notation of Section 1.
However, now we use the indices $m+1\lcd m+l$
to label the rows and columns of the matrix $E$.
Let $E(v)$ the inverse to the matrix $v-E\ts$. 
Let $E_{ab}(v)$ be the $a\com b$ entry of the inverse matrix.

In  this section we will use the symbol $\,\equiv\,$ to indicate 
equalities in the algebra $\U(\ts\f_{m+l})$ modulo the left ideal
generated by the elements of the subalgebra $\qp\subset\f_{m+l}\ts$.
Any two elements of $\U(\ts\f_{m+l})$ related by $\,\equiv\,$
act on the subspace \eqref{onesub} in the same way.
We will extend the relation $\,\equiv\,$ to 
formal power series in $u^{-1}$ with
coefficients in $\U(\ts\f_{m+l})\ts$, and then
to matrices whose entries are these series.

\begin{proposition*}
\label{prop1}
We have the relations\/ $\Bt\ts\equiv\ts0$ and\/ 
$\Dt\ts\equiv\ts E(v)\ts$.
\end{proposition*}

\begin{proof}
All entries of the matrix $B$ belong to $\qp\ts$.
All entries of the matrix $E$
belong to the subalgebra $\gl_{\ts l}\subset\f_{m+l}\ts$, and
adjoint action of $\gl_{\ts l}$ on $\f_{m+l}$
preserves $\qp$. Therefore
$$
\Dt=(\ts D-C\ts A^{-1}B\ts)^{-1}=(\ts v-E-C\ts A^{-1}B\ts)^{-1}
\ts\equiv\ts(\ts v-E\ts)^{-1}=\ts E(v)
$$
where we used Lemma~\ref{lemma0}. Similar arguments show that
$$
\Bt=-\ts A^{-1}B\ts(\ts D-C\ts A^{-1}B)^{-1}=-\ts A^{-1}B\ts\Dt
\ts\equiv\ts0\ts.
\eqno\qedhere
$$
\end{proof}

Consider the product of matrices $D^{-1}C$ appearing 
in the expressions for $\At$ and $\Ct$ from Lemma \ref{lemma0}.
This product is a $l\times(2m+l)$ matrix whose $a\com d$ entry
for the indices $a=m+1\lcd m+l$ and $d=-\ts m-l\lcd-1\com1\lcd m$ is the sum
over the index $b\ts$,
$$
-\,\sum_{b>m}\, 
E_{ab}(v)\,F_{\ts bd}\,.
$$
Let $Z(v)$ be the trace of the matrix $E(v)\ts$.
The coefficients of the formal power series $Z(v)$ in $v^{-1}$ belong to the
centre of the algebra $\U(\gl_{\ts l})\ts$, which is
now regarded as a subalgebra of $\U(\ts\f_{m+l})\ts$.
For any $c\com d=-\ts m-l\lcd-1\com1\lcd m$ put
$$
X_{cd\ts}(v)\,=\,
\de_{cd}\,(1-Z(v))\,\mp\,
\left\{
\begin{array}{ll}
E_{\ts-d,-c\ts}(v)
&\ \textrm{if}\ \,c\com d<-\,m\ts;\\[2pt]
\,0
&\ \textrm{otherwise.}
\end{array}
\right.
$$

\begin{lemma*}
\label{lemma1}
For any indices\/ $a=m+1\lcd m+l$ and\/ $d=-\ts m-l\lcd-1\com1\lcd m$
\begin{equation}
\label{3.3.1}
\sum_{b>m\ge c} 
E_{ab}(v)\,F_{\ts bc}\,X_{cd\ts}(v)\,\,=\,
\sum_{b>m}\,
F_{\ts bd}\,E_{ab}(v)\ts.
\end{equation}
\end{lemma*}

\begin{proof}
For any index $e=m+1\lcd m+l$ let us multiply both sides of
relation \eqref{3.3.1} by $\de_{ea}\ts v-E_{ea}$ on the left,
and then perform the summation over the index $a$. We get 
\begin{equation}
\label{3.3.2}
\sum_{c\le m} 
F_{ec}\,X_{cd\ts}(v)\,\,=\,
\sum_{a,b>m}\,
(\ts\de_{ea}\ts v-E_{ea})\,
F_{\ts bd}\,E_{ab}(v)\ts.
\end{equation}
By the commutation relations \eqref{ufmrel} in $\U(\ts\f_{m+l})\ts$,
the right hand side of \eqref{3.3.2} equals
$$
\sum_{a,b>m}\,
(\ts F_{\ts bd}\,(\ts\de_{ea}\ts v-E_{ea})
-\de_{ab}\,F_{ed}-\de_{-a,d}\,F_{\ts b,-e}\ts)\,
E_{ab}(v)
$$
$$
=\,F_{ed}\,(1-Z(v))\,\mp\,
\sum_{a,b>m}\de_{-a,d}\,F_{e,-b}\,E_{ab}(v)\,.
$$
By the definition of the series $X_{cd\ts}(v)\ts$,
the last displayed expression coincides with the left hand side
of \eqref{3.3.2}.
\end{proof}

Let us now write the matrix $A-B\ts D^{-1}\ts C$ and its inverse $\At$ as
\begin{equation}
\label{ghij}
\begin{bmatrix}
\,G\,&H\,\\\,I\,&J\,
\end{bmatrix}
\quad\textrm{and}\quad
\begin{bmatrix}\ \Gt\,&\Ht\ \\ \ \It\,&\Jt\ 
\end{bmatrix}
\end{equation}
where the blocks $G,H,I,J$ and $\Gt,\Ht,\It,\Jt$ are matrices of sizes
$l\times l$, $l\times2m$, $2m\times l$, $2m\times2m$ 
respectively. We label the rows and columns of all these blocks
by the same indices as in the compound matrices \eqref{ghij}. For instance,
the rows and columns of the $l\times l$ matrix
$G$ will be labelled by the indices $-\ts m-l\lcd-m-1\ts$.

Keeping to the notation of Section 2, let $F$ be the $2m\times2m$
matrix whose $c\com d$ entry is $F_{cd}$ for
$c\com d=-\ts m\lcd-1\com1\lcd m\ts$. Let $F(u)$ be the inverse 
to the matrix $u-F\ts$. The entries of this matrix are
formal power series in $u^{-1}$ with coefficients in the
algebra $\U(\ts\f_m)\ts$, see \eqref{fabu}.
But now the algebra $\U(\ts\f_m)$ is regarded as a as subalgebra of
$\U(\ts\f_{m+l})\ts$.

\begin{proposition*}
\label{prop2}
We have the relations\/ $\Ht\ts\equiv\ts0$ and\/ 
$\Jt\ts\equiv\ts(1-Z(v))\,F(v-l\ts)\ts$.
\end{proposition*}

\begin{proof}
The matrices $H$ and $J$ are the blocks of 
$\,A-B\ts D^{-1}\ts C\,$ corresponding to the rows
$-\ts m-l\lcd-\ts m-1$ and $-\ts m\lcd-1\com1\lcd m\ts$.
Both blocks correspond to the same columns $-\ts m\lcd-1\com1\lcd m\ts$.
For any indices $c=-\ts m-l\lcd-1\com1\lcd m$ and
$d=-\ts m\lcd-1\com1\lcd m$ the $c\com d$ entry of the matrix
$A$ is $\de_{cd}\,v-F_{cd}\ts$.
For the same indices the $c\com d$ entry
of the matrix $\,B\ts D^{-1}\ts C\,(1-Z(v))$ is
$$
\sum_{a,b>m}
F_{ca}\,
E_{ab}(v)\,F_{\ts bd}\,(1-Z(v))
\,=
\sum_{a,b>m} 
F_{ca}\,F_{\ts bd}\,E_{ab}(v)\,=
$$
$$
\sum_{a,b>m}
(\ts F_{\ts bd}\,F_{ca}+\de_{ab}\,F_{cd}-\de_{cd}\,F_{\ts ba}
\,\pm\,\de_{-b,c}\,F_{-a,d}\ts)
\,E_{ab}(v)
$$
where we used Lemma \ref{lemma1} and 
the commutation relations \eqref{ufmrel} in $\U(\ts\f_{m+l})\ts$.
The factors $F_{ca}$ and $F_{-a,d}$
belong to $\qp$. So the last displayed sum 
is related by $\,\equiv\,$ to
$$
F_{cd}\,Z(v)
\,-\sum_{a,b>m}\de_{cd}\,F_{\ts ba}\,E_{ab}(v)
\,=\,F_{cd}\,Z(v)+\de_{cd}\,(\ts l-\,v\,Z(v))\,.
$$
Therefore $c\com d$ entry of the matrix $\,A-B\ts D^{-1}\ts C\,$
is related by $\,\equiv\,$ to
$$
\de_{cd}\ts v-F_{cd}-F_{cd}\,\frac{Z(v)}{\,1-Z(v)}\,-\,
\de_{cd}\,\,\frac{l-\,v\,Z(v)}{\,1-Z(v)}
\,=\,(\ts\de_{cd}\,(v-l\ts)-F_{cd}\ts)\,\frac1{\,1-Z(v)}\ .
$$
This relation implies that $H\ts\equiv\ts0\ts$, and that
$$
J\ts\equiv\ts(v-l-F\ts)\,\frac1{\,1-Z(v)}
\ts.
$$
The adjoint action of the subalgebra $\f_m\subset\f_{m+l}$
on $\f_{m+l}$ preserves $\qp$. Hence
$$
\Jt=(\ts J-I\ts G^{-1}H\ts)^{-1}\ts\equiv\ts(1-Z(v))\,(\ts v-l-F\ts)^{-1}
\textstyle
=(1-Z(v))\,F\ts(v-l)
$$
where we used Lemma~\ref{lemma0}.
Similar arguments show that
$$
\Ht=-\ts G^{-1}H\ts(\ts J-I\ts G^{-1}H)^{-1}=-\ts G^{-1}H\ts\Jt
\ts\equiv\ts0\ts.
\eqno\qedhere
$$
\end{proof}

We have
\begin{equation}
\label{ct}
\Ct=-\ts D^{-1}C\ts(A-B\ts D^{-1}\ts C)^{-1}=-\ts D^{-1}C\ts\At
\end{equation}
by Lemma \ref{lemma0}. 
Let us now write the $l\times(\ts l+2m)$ matrices $C$ and $\Ct$ as
\begin{equation}
\label{pq}
\begin{bmatrix}
\,P\,&Q\,
\end{bmatrix}
\quad\textrm{and}\quad
\begin{bmatrix}
\,\Pt\,&\Qt\,
\end{bmatrix}
\end{equation}
where the blocks $P\com Q$ and $\Pt\com\Qt$ are matrices of sizes
$l\times l$ and $l\times2m$ respectively. 
We label the rows and columns of these blocks
by the same indices as in the compound matrices \eqref{pq}. 
For instance, the rows and columns of the matrix $P$ are labelled 
by the indices $m+1\lcd m+l$ and $-\ts m-l\lcd-m-1$ respectively.

\begin{proposition*}
\label{prop3}
For any\/ $a=m+1\lcd m+l$ and\/ $d=-\ts m\lcd-1\com1\lcd m$
$$
\Qt_{ad}\,\equiv
\sum_{b>m\ge c\ge-m}
F_{\ts bc}\,E_{ab}(v)\,F_{cd\ts}(v-l\ts)\ts.
$$
\end{proposition*}

\begin{proof}
Using \eqref{ct} and Proposition \ref{prop2}, 
$$
\Qt=-\,D^{-1}P\,\Ht-D^{-1}Q\,\Jt
\,\equiv\ts-\ts D^{-1}Q\,(1-Z(v))\,F(v-l\ts)\ts.
$$
Therefore the $a\com d$ entry of the matrix $\Qt$ is related by
$\,\equiv\,$ to the sum
$$
\sum_{b>m\ge c\ge-m}
E_{ab}(v)\,F_{\ts bc}\,(1-Z(v))\,F_{cd\ts}(v-l\ts)
$$
which is equal to the required sum by Lemma \ref{lemma1}.
\end{proof}

Denote by $Y(v)$ the trace of the second matrix in \eqref{abcd}\ts;
then $Y(v)$ plays the role of series $W(v)$ for the Lie algebra $\f_{m+l}$
instead of $\f_m\ts$. In particular, all coefficients of the series
$Y(v)$ belong to $\Z\ts(\ts\f_{m+l})\ts$.
Recall that the parameter $v$ is defined by \eqref{v}. Put
\begin{equation}
\label{w}
\textstyle
w=-\ts u\pm\frac12+m+l\ts.
\vspace{4pt}
\end{equation}

\begin{lemma*}
\label{lemma2}
We have the relation
\begin{equation}
\label{yu}
Y(v)\ts\pm\ts\frac1{2u}-1\,\equiv\ts
\bigl(\ts W(v-l\ts)\ts\pm\ts\frac1{2u}-1\ts\bigr)\,
\frac{1-Z(v)}{\,1-Z(w)}\ .
\end{equation}
\end{lemma*}

\begin{proof}
There is a unique power series in $u^{-1}$ with
coefficients in the subalgebra 
$\Z\ts(\ts\f_m\op\ts\gl_{\ts l})\subset\Z\ts(\ts\f_{m+l})$
related by $\,\equiv\,$ to the left hand side of \eqref{yu}.
To show that this unique series is the right hand
side of \eqref{yu}, it suffices to show that 
both sides have the same image under the Harish-Chandra
homomorphism $\psi_{m+l}\ts$.
This homomorphism corresponds to the triangular decomposition of
the Lie algebra $\f_{m+l}$ similar to \eqref{tridec}. 
Using the equality \eqref{hcwu} for $\f_{m+l}$ instead of
$\f_m\ts$, the image of the left hand side equals
$$
\Bigl(\ts\pm\,\frac1{2u}-1\Bigr)\,
\prod_{a=1}^{m+l}\,\ts
\Bigl(1-\frac1{v-m-l+a-F_{aa}}\ts\Bigr)
\Bigl(1-\frac1{m+l\pm1-v+a-F_{aa}}\ts\Bigr)^{\ns-1}\,.
$$
Applying $\psi_{m+l}$ to the right hand side of \eqref{yu}
amounts to applying the Harish-Chandra homomorphisms 
$\psi_m$ and $\varphi_{\ts l}$ to $W(v-l\ts)$ and to 
$Z(v)\com Z(w)$ respectively. Employing our embeddings of the Lie algebras 
$\f_m$ and $\gl_{\ts l}$ to $\f_{m+l}\ts$, we have 
$$
\psi_m(\ts W(v-l\ts))\ts\pm\ts\frac1{2u}-1\ts\,=\,
\Bigl(\ts\pm\,\frac1{2u}-1\Bigr)\ \times
$$
\begin{equation}
\label{yu1}
\prod_{a=1}^m\,\ts
\Bigl(1-\frac1{v-m-l+a-F_{aa}}\ts\Bigr)
\Bigl(1-\frac1{m+l\pm1-v+a-F_{aa}}\ts\Bigr)^{\ns-1},
\end{equation}
\begin{align}
\label{yu2}
1-\varphi_{\ts l}\ts(Z(v))
&\ =
\prod_{a=m+1}^{m+l}\,\ts
\Bigl(1-\frac1{v-m-l+a-F_{aa}}\ts\Bigr)\,,
\\
\label{yu3}
1-\varphi_{\ts l}\ts(Z(w))
&\ =
\prod_{a=m+1}^{m+l}\,\ts
\Bigl(1-\frac1{m+l\pm1-v+a-F_{aa}}\ts\Bigr)
\end{align}
by \eqref{hcwu} and \eqref{hczu}.
Multiplying \eqref{yu1} by \eqref{yu2} and dividing by \eqref{yu3}
we complete the proof of Lemma \ref{lemma2}.
\end{proof}

For any indices $a\com b=m+1\lcd m+l\ts$ denote
$\Et_{ab}(v)=(v-l-\Ep\ts)_{\ts ba}^{-1}\ts$.
Then by Lemma \ref{eep},
\begin{equation}
\label{eet}
(1-Z(v))\,\Et_{ab}(v)=E_{ab}(v)\ts.
\vspace{8pt}
\end{equation}

\begin{proposition*}
\label{prop4}
For any\/ $a\com b=m+1\lcd m+l$ and\/ $d=-\ts m\lcd-1\com1\lcd m$
\begin{align*}
\Gt_{-b,-a}
&\,\equiv\,
(1-Z(v))\,\Bigl(
\bigl(\ts W(v-l\ts)\ts\pm\ts\frac1{2u}-1\ts\bigr)\,
\Et_{ab}(w)\,\mp\,\frac1{2u}\ts\,\Et_{ab}(v)\Bigr)\,,
\\[8pt]
\It_{-d,-a}
&\,\equiv\,
\sum_{b>m\ge c\ge-m}\,
\ep_{ad}\,F_{\ts bc}\,(1-Z(v))\,\Bigl(\ts
\ep_{cd}\,\Et_{ab}(w)\,
F_{-d,-c\ts}(v-l\ts)
\\
&\hspace{48pt}
\pm\,\ts\frac{\Et_{ab}(w)-\Et_{ab}(v)}{2u}\,
F_{cd\ts}(v-l\ts)
\Bigr)\ts.
\end{align*}
\end{proposition*}

\begin{proof}
By applying Proposition \ref{1.2} to the Lie algebra 
$\f_{m+l}$ instead of $\f_m$ and 
then using Proposition~\ref{prop1} along with Lemma \ref{lemma2}, 
we prove that
$$
\Gt_{-b,-a}\ts\equiv\ts
\bigl(\ts W(v-l\ts)\ts\pm\ts\frac1{2u}-1\ts\bigr)
\,\ts\frac{1-Z(v)}{\,1-Z(w)}\ts\,
E_{ab}(w)\,\mp\,\frac1{2u}\ts\,E_{ab}(v)\,.
$$
This can be rewritten as the first equivalence relation 
of Proposition \ref{prop4} by using \eqref{eet}.
By applying Proposition \ref{1.2} to the Lie algebra
$\f_{m+l}$ instead of $\f_m$
and then using Proposition~\ref{prop3} along with Lemma \ref{lemma2}, 
we similarly prove that
$$
\It_{-d,-a}\ts\equiv\ts
\!\sum_{b>m\ge c\ge-m}
\ep_{ad}\,F_{\ts bc}\,\Bigl(
\bigl(\ts W(v-l\ts)\ts\pm\ts\frac1{2u}-1\ts\bigr)\,
\frac{1-Z(v)}{\,1-Z(w)}\ts\,E_{ab}(w)\,F_{cd\ts}(w-l\ts)
\vspace{-4pt}
$$
$$
\mp\,\,\frac1{2u}\,E_{ab}(v)\,F_{cd\ts}(v-l\ts)\Bigr)\ts.
\vspace{6pt}
$$
Now by applying Proposition \ref{1.2} to the Lie algebra $\f_m$ itself 
and using \eqref{eet}, we 
obtain the second equivalence relation of Proposition \ref{prop4}.
\end{proof}

\begin{lemma*}
\label{lemma3}
For any\/ $b\com e\com f=m+1\lcd m+l$ and\/ $c=-\ts m\lcd-1\com1\lcd m$
$$
E_{ef}(v)\,\frac1{1-Z(v)}\,F_{\ts bc}\,(1-Z(v))
\,=\,F_{\ts bc}\,E_{ef}(v)\ +
$$
$$
\sum_{a>m}\,
F_{ac}\,\bigl(\ts 
E_{ef}(v)\,
(\ts\Et_{\ts ba}(v-1)-\Et_{\ts ba}(v+1))\ts+\ts
E_{ea}(v)\,
(\ts\Et_{\ts bf}(v-1)+\Et_{\ts bf}(v+1))
\bigr)\tts\big/\tts2\,.
$$
\end{lemma*}

\begin{proof}
Due to \eqref{ufmrel},
for any $a\com p\com q=m+1\lcd m+l$ we have the relation
$$
[\ts F_{pq}\com F_{ac}\ts]=\de_{aq}\,F_{pc}\ts.
$$
Multiplying this relation by $E_{ep}(v)$ on the left, by
$E_{\ts qf}(v)$ on the right, and then taking the sum over the indices
$p\com q$ we get 
$$
[\ts E_{ef}(v)\com F_{ac}\ts]\,=
\sum_{p>m}\,E_{ep}(v)\,F_{pc}\,E_{af}(v)\ts.
$$
Applying now Lemma \ref{lemma1} we get the commutation relation
\begin{equation}
\label{efbc}
[\ts E_{ef}(v)\com F_{ac}\ts]\,=
\sum_{p>m}\,F_{pc}\,E_{ep}(v)\,E_{af}(v)\,
\frac1{1-Z(v)}\ .
\end{equation}

By setting $e=f$ \eqref{efbc} and then performing
summation over the index $f\ts$, we get 
$$
[\ts Z(v)\com F_{ac}\ts]\,=
\sum_{f,p>m}\,F_{pc}\,E_{fp}(v)\,E_{af}(v)\,
\frac1{1-Z(v)}\ =
$$
$$
\sum_{f,p>m}\,F_{pc}\,\Et_{fp}(v)\,\Et_{af}(v)\,
(1-Z(v))\ =
$$
$$
\sum_{f,p>m}\,F_{pc}\,
(v-l-\Ep\ts)^{-1}_{\ts pf}\,
(v-l-\Ep\ts)^{-1}_{fa}\,
(1-Z(v))\ =
$$
$$
\sum_{p>m}\,F_{pc}\,
(v-l-\Ep\ts)^{-2}_{pa}\,
(1-Z(v))
$$
where also we used the relation \eqref{eet} and
the definition of the series $\Et_{ab}(v)\ts$.
The last displayed expression for the
commutator $[\ts Z(v)\com F_{\ts bc}\ts]$ shows that
$$
\sum_{p>m}\,
\frac1{1-Z(v)}\,F_{pc}\,
(\ts1-(v-l-\Ep\ts)^{-2}\ts)_{pa}\,=\,
F_{\ts ac}\,\frac1{1-Z(v)}\ .
$$
By multiplying this relation on the right by the matrix entry
$(\ts1-(v-l-\Ep\ts)^{-2}\ts)_{\ts ab}^{-1}$ and then taking the sum
over the index $a$ we finally get
\begin{equation}
\label{zbc}
\frac1{1-Z(v)}\,F_{\ts bc}\,=
\sum_{a>m}\,
F_{ac}\,\frac1{1-Z(v)}\,(\ts1-(v-l-\Ep\ts)^{-2}\ts)_{\ts ab}^{-1}\,.
\end{equation}

By using the relations \eqref{efbc} and \eqref{zbc},
we obtain the equalities
$$
E_{ef}(v)\,\frac1{1-Z(v)}\,F_{\ts bc}\,(1-Z(v))\,=
\sum_{a>m}\,
E_{ef}(v)\,
F_{ac}\,(\ts1-(v-l-\Ep\ts)^{-2}\ts)_{\ts ab}^{-1}\ =
$$
$$
\sum_{a>m}\,
F_{ac}\,E_{ef}(v)\,
(\ts1-(v-l-\Ep\ts)^{-2}\ts)_{\ts ab}^{-1}
\ +
\vspace{4pt}
$$
$$
\sum_{a,p>m}
F_{pc}\,E_{ep}(v)\,(v-l-\Ep\ts)^{-1}_{fa}\,
(\ts1-(v-l-\Ep\ts)^{-2}\ts)_{\ts ab}^{-1}
$$
where we also used the relation \eqref{eet}. Observe that here
$$
\sum_{a>m}\,
(v-l-\Ep\ts)^{-1}_{fa}\,
(\ts1-(v-l-\Ep\ts)^{-2}\ts)_{ab}^{-1}=
(\ts(v-l-\Ep\ts)-(v-l-\Ep\ts)^{-1}\ts)^{-1}_{fb}\,.
$$
Using the matrix identities
$$
(\ts1-(v-l-\Ep\ts)^{-2}\ts)^{-1}\ts=\,
1\,+\,\frac{(v-l-1-\Ep\ts)^{-1}-(v-l+1-\Ep\ts)^{-1}}2\,,
$$
$$
(\ts(v-l-\Ep\ts)-(v-l-\Ep\ts)^{-1}\ts)^{-1}\ts=\,
\frac{(v-l-1-\Ep\ts)^{-1}+(v-l+1-\Ep\ts)^{-1}}2
$$
and the definition of the series $\Et_{ab}(v)\ts$, we
complete the proof of Lemma \ref{lemma3}.
\end{proof}

\begin{proposition*}
\label{prop5}
For any indices\/ $a\com e=m+1\lcd m+l$ we have the relation
$$
\Pt_{e,-a}\,\equiv\ 
-\sum_{b,f>m}
F_{\ts f,-b}\,E_{ef}(v)\,\Et_{ab}(w)
$$
$$
\ \mp\hspace{-5pt}
\sum\limits_{\substack{b,f>m\\m\ge c,d\ge-m}}\hspace{-5pt}
\ep_{ad}\,\ts
F_{\ts f,-d}\,F_{\ts bc}\,E_{eb}(v)\,
\Et_{af}(w)\,F_{cd\ts}(v-l\ts)\,.
$$
\end{proposition*}

\begin{proof}
By \eqref{ct} we have
$$
\Pt=-\,D^{-1}P\,\Gt-D^{-1}Q\,\It\ts.
$$
So for any indices $a\com e=m+1\lcd m+l$
the $e\com\ns-\ts a$ entry of the matrix $\Pt$ equals
$$
\sum_{b,f>m}
E_{ef}(v)\,F_{\ts f,-b}\,\Gt_{-b,-a}
\ +
\sum_{f>m\ge d\ge-m}
E_{ef}(v)\,F_{\ts f,-d}\,\It_{-d,-a}\,.
$$
Using Proposition \ref{prop4}, 
the latter sum is related by $\,\equiv\,$ to
$$
\sum_{b,f>m}
E_{ef}(v)\,F_{\ts f,-b}\,(1-Z(v))\,
\Bigl(
\bigl(\ts W(v-l\ts)\ts\pm\ts\frac1{2u}-1\ts\bigr)\,
\Et_{ab}(w)\,\mp\,\frac1{2u}\ts\,\Et_{ab}(v)
\Bigr)\ +
\vspace{4pt}
$$
$$
\sum\limits_{\substack{b,f>m\\m\ge c,d\ge-m}}\hspace{-5pt}
\ep_{ad}\,\ts
E_{ef}(v)\,F_{\ts f,-d}\,
F_{\ts bc}\,(1-Z(v))\ \times
$$
\begin{equation}
\label{thrlin}
\Bigl(\ts
\ep_{cd}\,\Et_{ab}(w)\,F_{-d,-c\ts}(v-l\ts)
\ts\pm\ts
\frac{\ts\Et_{ab}(w)-\Et_{ab}(v)}{2u}\,
F_{cd\ts}(v-l\ts)\Bigr)\,.
\vspace{4pt}
\end{equation}

Using \eqref{eet} and the definition of $\Et_{ab}(v)\ts$,
it follows from Lemma \ref{lemma1} that for any $b=m+1\lcd m+l$ we have
$$
\sum_{f>m}
E_{ef}(v)\,F_{\ts f,-b}\,(1-Z(v))
\,=\,
\sum_{f,p>m}
F_{\ts f,-p}\,E_{ef}(v)\,
\bigl(\,\de_{\ts bp}\ts\pm\ts\Et_{\ts bp}(\ts v\mp1\ts)\bigr)\ts.
$$
Here the index $p$ is ranging over
$m+1\lcd m+l\ts$ as well as the index $f$ does. Further,
by the definition of the series $\Et_{\ts ab}(v)$ we have
\begin{equation}
\label{can1}
\sum_{b>m}\,\Et_{\ts bp}(\ts v\mp1\ts)\,\Et_{ab}(w)\,=\,
\frac{\ts\Et_{ap}(w)-\Et_{ap}(\ts v\mp1\ts)\ts}{2u\mp1}
\end{equation}
and
\begin{equation}
\label{can2}
\sum_{b>m}\,\Et_{\ts bp}(\ts v\mp1\ts)\,\Et_{ab}(v)\,=\,
\mp\,\ts\Et_{ap}(v)\tts\pm\tts\Et_{ap}(\ts v\mp1\ts)\ts.
\end{equation}
So after cancellations the sum over $b\com f$ in the first of the three
lines \eqref{thrlin} equals
\begin{equation}
\label{toadd}
\sum_{f,p>m}
F_{\ts f,-p}\,E_{ef}(v)\,\Bigl(\ts 
W(v-l)\,
\frac{\ts2u\,\Et_{ap}(w)\mp\Et_{ap}(\ts v\mp1\ts)}{2u\mp1}
\ts-\ts\Et_{ap}(w)\Bigr)\,.
\end{equation}

Also by Lemma \ref{lemma1}, for any index $d=-\ts m\lcd-1\com1\lcd m$
we have the relation
$$
\sum_{f>m}
E_{ef}(v)\,F_{\ts f,-d}
\,=\,
\sum_{f>m}
F_{\ts f,-d}\,E_{ef}(v)\,\frac1{1-Z(v)}\ .
$$
Using it with Lemma \ref{lemma3}, 
the sum over $b\com c\com d\com f$ in the last two lines
of \eqref{thrlin} equals
$$
\sum\limits_{\substack{b,f>m\\m\ge c,d\ge-m}}\hspace{-5pt}
\ep_{ad}\,\ts
F_{\ts f,-d}\,\Bigl(
F_{\ts bc}\,E_{ef}(v)\,+
\sum_{p>m}\,
F_{pc}\,\bigl(\ts
E_{ef}(v)\,
(\ts\Et_{\ts bp}(v-1)-\Et_{\ts bp}(v+1))\ \ts+
\vspace{-8pt}
$$
$$
E_{ep}(v)\,
(\ts\Et_{\ts bf}(v-1)+\Et_{\ts bf}(v+1))
\bigr)\tts\big/\tts2\Bigr)\,\ts
\ep_{cd}\,\Et_{ab}(w)\,F_{-d,-c\ts}(v-l\ts)
\vspace{12pt}
$$
$$
\pm\hspace{-5pt}
\sum\limits_{\substack{b,f>m\\m\ge c,d\ge-m}}\hspace{-5pt}
\ep_{ad}\,\ts
F_{\ts f,-d}\,\Bigl(
F_{\ts bc}\,E_{ef}(v)\,+
\sum_{p>m}\,
F_{pc}\,\bigl(\ts
E_{ef}(v)\,
(\ts\Et_{\ts bp}(v-1)-\Et_{\ts bp}(v+1))\ \ts+
\vspace{-2pt}
$$
$$
E_{ep}(v)\,
(\ts\Et_{\ts bf}(v-1)+\Et_{\ts bf}(v+1))
\bigr)\tts\big/\tts2\Bigr)\,\ts
\bigl(\ts\Et_{ab}(w)-\Et_{ab}(v)\bigr)\,
F_{cd\ts}(v-l\ts)\ts\big/\ts{2u}\Bigr)\,.
\vspace{2pt}
$$
Note by that by \eqref{ufmrel}
here we have the commutation relation
\begin{equation}
\label{nonab}
[\ts F_{\ts f,-d}\,,F_{\ts bc}\ts]=\de_{cd}\,\ts\ep_{f,-d}\,F_{b,-f}\ts.
\end{equation}
In the first two of the four lines displayed above of \eqref{nonab},
let us open the outer brackets and then
permute the factor $F_{\ts f,-d}$ with each of the factors
$F_{\ts bc}$ and $F_{\ts pc}$ in the summands.
Let us then also change the running indices $c\com d$ to $-\ts d\com-\ts c$
in the resulting sum. By adding to this sum the sum in the last
two of the four diplayed lines, we get
$$
\sum\limits_{\substack{b,f>m\\m\ge c,d\ge-m}}\hspace{-5pt}
\ep_{a,-d}\,\ts
F_{\ts f,-d}\,F_{\ts bc}\,E_{eb}(v)\,
\Et_{af}(w)\,F_{cd\ts}(v-l\ts)\ \ts+
\vspace{-4pt}
$$
$$
\sum\limits_{\substack{b,f,p>m\\m\ge c,d\ge-m}}\hspace{-5pt}
\ep_{a,-d}\,\ts
F_{\ts f,-d}\,
F_{pc}\,\bigl(\ts
E_{ep}(v)\,
(\ts\Et_{\ts bf}(v-1)-\Et_{\ts bf}(v+1))\ \ts+
\vspace{-1pt}
$$
$$
E_{ef}(v)\,
(\ts\Et_{\ts bp}(v-1)+\Et_{\ts bp}(v+1))
\bigr)\,
\Et_{ab}(w)\,F_{cd\ts}(v-l\ts)\ts/\ts2
\vspace{15pt}
$$
$$
\pm\hspace{-5pt}
\sum\limits_{\substack{b,f>m\\m\ge c,d\ge-m}}\hspace{-5pt}
\ep_{ad}\,\ts
F_{\ts f,-d}\,\Bigl(
F_{\ts bc}\,E_{ef}(v)\,+
\sum_{p>m}\,
F_{pc}\,\bigl(\ts
E_{ef}(v)\,
(\ts\Et_{\ts bp}(v-1)-\Et_{\ts bp}(v+1))\ \ts+
\vspace{-4pt}
$$
$$
E_{ep}(v)\,
(\ts\Et_{\ts bf}(v-1)+\Et_{\ts bf}(v+1))
\bigr)\tts\big/\tts2\Bigr)\,
(\ts\Et_{ab}(w)-\Et_{ab}(v))
\,F_{cd\ts}(v-l\ts)\ts/\ts2u\,.
\vspace{3pt}
$$
Note that here $\ep_{a,-d}\ts=\ts\mp\,\ep_{ad}\,$.
Using relations similar to \eqref{can1} and \eqref{can2},
all summands in the last four of the five lines displayed above cancel.
Hence by using the commutation relation \eqref{nonab},
the sum over the indices $b\com c\com d\com f$ in the last two lines 
of \eqref{thrlin} equals

$$
\sum\limits_{\substack{b,f>m\\m\ge c,d\ge-m}}\hspace{-5pt}
\ep_{a,-d}\,\ts
F_{\ts f,-d}\,F_{\ts bc}\,E_{eb}(v)\,
\Et_{af}(w)\,F_{cd\ts}(v-l\ts)\ \ts+
$$
$$
\sum\limits_{\substack{b,f>m\\m\ge c,d\ge-m}}\hspace{-5pt}
\ep_{ad}\,\ts
\de_{cd}\,\ts
\ep_{f,-d}\,
\Bigl(
F_{\ts b,-f}\,E_{ef}(v)\,+
\sum_{p>m}\,
F_{p,-f}\,\bigl(\ts
E_{ef}(v)\,
(\ts\Et_{\ts bp}(v-1)-\Et_{\ts bp}(v+1))
$$
$$
+\ E_{ep}(v)\,
(\ts\Et_{\ts bf}(v-1)+\Et_{\ts bf}(v+1))
\bigr)\tts\big/\tts2\Bigr)\,\ts
\ep_{cd}\,\Et_{ab}(w)\,F_{-d,-c\ts}(v-l\ts)\ =
\vspace{14pt}
$$
$$
\mp\hspace{-5pt}
\sum\limits_{\substack{b,f>m\\m\ge c,d\ge-m}}\hspace{-5pt}
\ep_{ad}\,\ts
F_{\ts f,-d}\,F_{\ts bc}\,E_{eb}(v)\,
\Et_{af}(w)\,F_{cd\ts}(v-l\ts)
$$
$$
\mp
\sum_{b,f>m}
\Bigl(
F_{\ts b,-f}\,E_{ef}(v)\,+
\sum_{p>m}\,
F_{p,-f}\,\bigl(\ts
E_{ef}(v)\,
(\ts\Et_{\ts bp}(v-1)-\Et_{\ts bp}(v+1))\ \ts+
\vspace{6pt}
$$
$$
E_{ep}(v)\,
(\ts\Et_{\ts bf}(v-1)+\Et_{\ts bf}(v+1))
\bigr)\tts\big/\tts2\Bigr)\,
\Et_{ab}(w)\,W(v-l\ts)\ =
\vspace{14pt}
$$
$$
\mp\hspace{-5pt}
\sum\limits_{\substack{b,f>m\\m\ge c,d\ge-m}}\hspace{-5pt}
\ep_{ad}\,\ts
F_{\ts f,-d}\,F_{\ts bc}\,E_{eb}(v)\,
\Et_{af}(w)\,F_{cd\ts}(v-l\ts)
$$
$$
\mp
\sum_{b,f>m}
F_{\ts b,-f}\,E_{ef}(v)\,\Et_{ab}(w)\,W(v-l\ts)
\vspace{8pt}
$$
\begin{equation}
\label{mulin}
\mp\,
\sum_{b,f,p>m}
F_{f,-p}\,E_{ef}(v)\,\Et_{\ts bp}(v\mp1)\,\Et_{ab}(w)\,W(v-l\ts)\,.
\vspace{4pt}
\end{equation}
To get the last equality, we used the relation 
$F_{p,-f}=\pm\,F_{f,-p}\ts$. Using the relation \eqref{can1},
the expression in the last line of \eqref{mulin} equals
$$
\mp\,
\sum_{f,p>m}
F_{f,-p}\,E_{ef}(v)\,
\frac{\ts\Et_{ap}(w)-\Et_{ap}(\ts v\mp1\ts)\ts}{2u\mp1}
\,W(v-l\ts)\,.
$$
By adding \eqref{toadd} to the last displayed expression, we get the sum
$$
\sum_{f,p>m}
F_{f,-p}\,E_{ef}(v)\,\Et_{ap}(w)
\,(\ts W(v-l\ts)-1\ts)\,.
$$
The sum \eqref{thrlin} is equal to \eqref{toadd} plus \eqref{mulin}. 
Thus we have shown that \eqref{thrlin} equals
$$
\mp\hspace{-5pt}
\sum\limits_{\substack{b,f>m\\m\ge c,d\ge-m}}\hspace{-5pt}
\ep_{ad}\,\ts
F_{\ts f,-d}\,F_{\ts bc}\,E_{eb}(v)\,
\Et_{af}(w)\,F_{cd\ts}(v-l\ts)
$$
$$
\mp
\sum_{b,f>m}
F_{\ts b,-f}\,E_{ef}(v)\,\Et_{ab}(w)\,W(v-l\ts)
$$
\vspace{4pt}
$$
\,+\!
\sum_{f,p>m}
F_{f,-p}\,E_{ef}(v)\,\Et_{ap}(w)
\,(\ts W(v-l\ts)-1\ts)\,=
\vspace{6pt}
$$
$$
\mp\hspace{-5pt}
\sum\limits_{\substack{b,f>m\\m\ge c,d\ge-m}}\hspace{-5pt}
\ep_{ad}\,\ts
F_{\ts f,-d}\,F_{\ts bc}\,E_{eb}(v)\,
\Et_{af}(w)\,F_{cd\ts}(v-l\ts)
$$
$$
\,-\!
\sum_{f,p>m}
F_{\ts f,-p}\,E_{ef}(v)\,\Et_{ap}(w)\,.
$$
By changing here the index $p$ to $b$ we complete the proof of
Proposition \ref{prop5}.
\end{proof}

Let us now summarize what we did so far 
to prove that $\chi$ is an intertwining operator.
We considered the $(2m+2l)\times(2m+2l)$ matrix whose
$a\com b$ entry is $\de_{ab}\,v-F_{ab}$ 
for any indices \eqref{mliml}. We have shown that its
inverse 
is related by $\,\equiv\,$ to the block matrix
$$
\begin{bmatrix}
\,\Gt\,&0\,&0\,
\\
\,\It\,&\Jt\,&0\,
\\
\,\Pt\,&\Qt\,&\Dt\,
\end{bmatrix}
$$
where the blocks $\Dt\com\Gt\com\It\com\Jt\com\Pt\com\Qt$
are described by the five propositions above. By the definition
of the $\X(\g_n)\ts$-module $M=\F_{m+l}\ts(\ts V\ns\bt U\ts)\ts$,
the action of $\X(\g_n)$ on the elements of the subspace
\begin{equation}
\label{1vu}
1\ot V\ns\ot U\ot\P\ts(\CC^{\ts m+l}\ot\CC^{\ts n})\subset M
\end{equation}
can be described by assigning to every series $S_{ij}(u)$ 
the following sum of series with coefficients in the algebra
$\B_{m+l}=\U(\ts\f_{m+l})\ot\PD\ts(\CC^{\ts m+l}\ot\CC^{\ts n})\ts$:
$$
\de_{ij}\,-\sum_{a,b>m}\,
\Dt_{ab}\ts\ot\ts\th_i\,\th_j\,\d_{a\bi}\,x_{\ts b\bj}\ \,+
$$
$$
\sum_{a,b>m}\,
\Gt_{-b,-a}\ts\ot\ts x_{\ts bi}\,\d_{aj}\ \,+
$$
$$
\sum_{a>m\ge d>0}\,(\,
\It_{-d,-a}\ts\ot\ts x_{di}\,\d_{aj}
\ts-\ts
\It_{d,-a}\ts\ot\ts\th_i\,\d_{d\bi}\,\d_{aj}\,)\ \,+
$$
$$
\sum_{m\ge c,d>0}\,(\,
\Jt_{-c,-d}\ts\ot\ts x_{ci}\,\d_{dj}
+
\Jt_{-c,d}\ts\ot\ts\th_j\,x_{ci}\,x_{d\bj}
-
\Jt_{c,-d}\ts\ot\ts\th_i\,\d_{c\bi}\,\d_{dj}
-
\Jt_{cd}\ts\ot\ts\th_i\,\th_j\,\d_{c\bi}\,x_{d\bj}\,)
$$
$$
-\,
\sum_{a,e>m}\,
\Pt_{e,-a}\ts\ot\ts\th_i\,\d_{e\bi}\,\d_{aj}
$$
\begin{equation}
\label{fivlin}
-\,
\sum_{a>m\ge d>0}\,(\ts
\Qt_{a,-d}\ts\ot\ts\th_i\,\d_{a\bi}\,\d_{dj}
\ts+\ts
\Qt_{ad}\ts\ot\ts\th_i\,\th_j\,\d_{a\bi}\,x_{d\bj}\,)\,.
\end{equation}
Here we use the standard coordinate functions
$x_{ai}$ on the vector space $\CC^{\ts m+l}\ot\CC^{\ts n}$
with $a=1\lcd m+l$ and $i=1\lcd n\ts$.
Then $\d_{ai}$ is the partial derivation on 
$\P\ts(\CC^{\ts m+l}\ot\CC^{\ts n})$ relative to $x_{ai}\ts$.
The functions $x_{ai}$ with $a\leqslant m$ and $a>m$
correspond to the first and the second
direct summands in \eqref{mnl}.

Consider the action of $\X(\g_n)$
on the elements of the subspace \eqref{1vu} modulo $\q\cdot M\ts$, using
the definition \eqref{gan} where $m$ is to be replaced by $m+l\ts$.
From now till the end of this section,
we will be assuming that
$a\com b\com e\com f=m+1\lcd m+l$ while $c\com d=1\lcd m\ts$.
The indices $g\com h$ and $k$ will run through $1\lcd n\ts$. 

By Proposition \ref{prop1},
the action of the sum displayed in the first of six lines \eqref{fivlin}
on the elements of \eqref{1vu}
coincides with action of
$$
\de_{ij}\,\ts-\ts\sum_{a,b}
E_{ab}(v)\ts\ot\ts\th_i\,\th_j\,\d_{a\bi}\,x_{\ts b\bj}\ =
$$
\begin{equation}
\label{term0}
\de_{ij}\,(1-Z(v))\ts-\ts\sum_{a,b}
E_{ab}(v)\ts\ot\ts\th_i\,\th_j\,x_{\ts b\bj}\,\d_{a\bi}\,.
\end{equation}

By Proposition \ref{prop4},
the action of the sum displayed in the second of six lines \eqref{fivlin}
on the elements of \eqref{1vu}
coincides with action of the sum over 
$a\com b$ of
\begin{equation}
\label{term1}
(1-Z(v))\ts\Bigl(
\bigl(\ts W(v-l\ts)\ts\pm\ts\frac1{2u}-1\ts\bigr)\,
\Et_{ab}(w)\,\mp\,\frac1{2u}\ts\,\Et_{ab}(v)\Bigr)
\ts\ot\ts x_{\ts bi}\,\d_{aj}\,.
\vspace{4pt}
\end{equation}

By Proposition \ref{prop4}, the action of the 
sum in the third of six lines \eqref{fivlin} 
on the elements of \eqref{1vu}
coincides with action of the sum over the indices 
$a\com b\com c\com d$ of
\begin{align*}
\mp\,& 
F_{\ts b,-c}(1-Z(v))
\Bigl(\ns
\Et_{ab}(w) F_{-d,c\ts}(v-l)
-
\frac{\Et_{ab}(w)\ns-\ns\Et_{ab}(v)}{2u}
F_{-c,d}(v-l\ts)
\!\Bigr)
\ot x_{di}\ts\d_{aj}
\\[6pt]
+\,& 
F_{\ts bc}\,(1-Z(v))
\Bigl(\ns
\Et_{ab}(w)F_{-d,-c\ts}(v-l)
\pm
\frac{\Et_{ab}(w)\ns-\ns\Et_{ab}(v)}{2u}
F_{cd\ts}(v-l)
\!\Bigr)
\ot x_{di}\ts\d_{aj}
\\[6pt]
\pm\,& 
F_{\ts b-c}(1-Z(v))\Bigl(\ns
\Et_{ab}(w)
F_{dc\ts}(v-l)
\pm\frac{\Et_{ab}(w)\ns-\ns\Et_{ab}(v)}{2u}
F_{-c,-d\ts}(v-l)
\!\Bigr)
\ot\th_i\ts\d_{d\bi}\ts\d_{aj}
\\[6pt]
-\,&
F_{\ts bc}(1-Z(v))\,\Bigl(\ns
\Et_{ab}(w)
F_{d,-c\ts}(v-l)
-\frac{\Et_{ab}(w)\ns-\ns\Et_{ab}(v)}{2u}
F_{c,-d\ts}(v-l)
\!\Bigr)
\ot\th_i\ts\d_{d\bi}\ts\d_{aj}\,.
\end{align*}
Here $F_{\ts b,-c}\in\q$ and $F_{\ts bc}\in\q\ts$. Hence
modulo $\q\cdot M\ts$, the expression displayed in the latter eight
lines acts on the elements of \eqref{1vu} as the sum over the index $k$ of
\begin{align*}
\Bigl(\, 
\mp\,&
\Bigl(\ts
\Et_{ab}(w)\,F_{-d,c\ts}(v-l\ts)
-
\frac{\Et_{ab}(w)-\Et_{ab}(v)}{2u}\,
F_{-c,d\ts}(v-l\ts)
\Bigr)
\ot
\th_k\,x_{\ts b\bk}\,x_{ck}\,
x_{di}\,\d_{aj}
\\[6pt]
-\,&
\Bigl(\ts
\Et_{ab}(w)\,F_{-d,-c\ts}(v-l\ts)
\pm
\frac{\Et_{ab}(w)-\Et_{ab}(v)}{2u}\,
F_{cd\ts}(v-l\ts)
\Bigr)
\ot x_{\ts bk}\,\d_{ck}\,x_{di}\,\d_{aj}
\\[6pt]
\pm\,&
\Bigl(\ts
\Et_{ab}(w)\,
F_{dc\ts}(v-l\ts)
\pm\ts\frac{\Et_{ab}(w)-\Et_{ab}(v)}{2u}\,
F_{-c,-d\ts}(v-l\ts)
\Bigr)
\ot\th_i\,
\th_k\,x_{\ts b\bk}\,x_{ck}\,
\d_{d\bi}\,\d_{aj}
\\[6pt]
+\,&
\Bigl(\ts
\Et_{ab}(w)\,
F_{d,-c\ts}(v-l\ts)
-\frac{\Et_{ab}(w)-\Et_{ab}(v)}{2u}\,
F_{c,-d\ts}(v-l\ts)
\Bigr)
\ot
\th_i\,
x_{\ts bk}\,\d_{ck}\,
\d_{d\bi}\,\d_{aj}\,\Bigr)
\end{align*}
\vspace{-16pt}
\begin{equation}
\label{term2}
\times\ \bigl(\ts(1-Z(v))\ot1\ts\bigr)\,.
\vspace{4pt}
\end{equation}

By Proposition \ref{prop2}, the action of the 
sum displayed in the fourth of six lines \eqref{fivlin} 
on the elements of \eqref{1vu}
coincides with the action of the sum over $c\com d$ of
$$
\bigl(\ts(1-Z(v))\ot1\ts\bigr)\,\bigl(\,
F_{-c,-d\ts}(v-l\ts)\ts\ot\ts x_{ci}\,\d_{dj}
\ts+\ts
F_{-c,d\ts}(v-l\ts)\ts\ot\ts\th_j\,x_{ci}\,x_{d\bj}
$$
\begin{equation}
\label{term3}
-\,\ts
F_{c,-d\ts}(v-l\ts)\ts\ot\ts\th_i\,\d_{c\bi}\,\d_{dj}
\ts-\ts
F_{cd\ts}(v-l\ts)\ts\ot\ts\th_i\,\th_j\,\d_{c\bi}\,x_{d\bj}\,\bigr)\,.
\vspace{4pt}
\end{equation}

By Proposition \ref{prop5}, the action of the 
sum in the fifth of six lines \eqref{fivlin} 
on the elements of \eqref{1vu} coincides with action 
of the sum over the indices $a\com b\com e\com f$~of
$$
F_{\ts f,-b}\,E_{ef}(v)\,\Et_{ab}(w)
\ts\ot\ts\th_i\,\d_{e\bi}\,\d_{aj}
$$
plus the action of the sum over the indices 
$a\com b\com c\com d\com e\com f$ of
\begin{align*}
&-\,
F_{fd}\,
F_{\ts b,-c}\,
E_{eb}(v)\,\Et_{af}(w)\,F_{-c,-d\ts}(v-l)
\ts\ot\ts\th_i\,\d_{e\bi}\,\d_{aj}
\\[6pt]
&\pm\,
F_{f,-d}\,
F_{\ts b,-c}\,
E_{eb}(v)\,\Et_{af}(w)\,F_{-c,d}(v-l)
\ts\ot\ts\th_i\,\d_{e\bi}\,\d_{aj}
\\[6pt]
&-\,
F_{fd}\,
F_{\ts bc}\,
E_{eb}(v)\,\Et_{af}(w)\,F_{c,-d\ts}(v-l)
\ts\ot\ts\th_i\,\d_{e\bi}\,\d_{aj}
\\[6pt]
&\pm\,
F_{f,-d}\,
F_{\ts bc}\,
E_{eb}(v)\,\Et_{af}(w)\,F_{cd\ts}(v-l)
\ts\ot\ts\th_i\,\d_{e\bi}\,\d_{aj}\,.
\end{align*}
Modulo $\q\cdot M\ts$, here the expression to be
summed over 
the indices $a\com b\com e\com f$ acts
on the elements of the subspace \eqref{1vu} as the sum over the index $k$ of
$$
E_{ef}(v)\,\Et_{ab}(w)
\ts\ot\ts
\th_i\,\th_k\,
x_{f\bk}\,x_{\ts bk}\,
\d_{e\bi}\,\d_{aj}
$$
while the expression to be summed over $a\com b\com c\com d\com e\com f$ acts 
as the sum over $g\com h$ of
\begin{align*}
\bigl(\ts E_{eb}(v)\,\Et_{af}(w)\ot1\ts\bigr)\,
\bigl(
&\,F_{-c,-d\ts}(v-l)
\ts\ot\ts
\th_i\,
\th_g\,
x_{\ts b\bg}\,x_{cg}\,
x_{fh}\,\d_{dh}\,
\d_{e\bi}\,\d_{aj}
\\[6pt]
\pm&\ 
F_{-c,d}(v-l)
\ts\ot\ts
\th_i\,\th_g\,\th_h\,
x_{\ts b\bg}\,x_{cg}\,
x_{f\bh}\,x_{dh}\,
\d_{e\bi}\,\d_{aj}
\\[6pt]
-&\ 
F_{c,-d\ts}(v-l)
\ts\ot\ts
\th_i\,
x_{\ts bg}\,\d_{cg}\,
x_{\ts fh}\,\d_{dh}\,
\d_{e\bi}\,\d_{aj}
\\[6pt]
\mp&\ 
F_{cd\ts}(v-l)
\ts\ot\ts
\th_i\,\th_h\,
x_{\ts bg}\,\d_{cg}\,
x_{f\bh}\,x_{dh}\,
\d_{e\bi}\,\d_{aj}\ts\bigr)\,.
\end{align*}

We have $\th_{\bk}=\pm\,\th_k$ for $k=1\lcd n\ts$.
Using the commutation relations in the ring
$\PD\ts(\CC^{\ts m+l}\ot\CC^{\ts n})\ts$,
the sum over the index $k$ above equals the sum over $k$ of
\begin{equation}
\label{term41}
E_{ef}(v)\,\Et_{ab}(w)
\ts\ot\ts
\th_i\,\th_k\,
x_{f\bk}\,\d_{e\bi}\,
x_{\ts bk}\,\d_{aj}
\end{equation}
plus
\begin{equation}
\label{term410}
\mp\,\ts\de_{\ts be}\,
E_{ef}(v)\,\Et_{ab}(w)
\ts\ot\ts
x_{f\bi}\,\d_{aj}\ts.
\vspace{4pt}
\end{equation}
Similarly, the sum over the indices
$g\com h$ equals the sum over $g\com h$ of
\begin{gather}
\nonumber
\bigl(\,
F_{-c,-d\ts}(v-l\ts)\ts\ot\ts x_{cg}\,\d_{dh}
\ts+\ts
F_{-c,d\ts}(v-l\ts)\ts\ot\ts\th_h\,x_{cg}\,x_{d\bh}
\\[4pt]
\nonumber
-\,\ts
F_{c,-d\ts}(v-l\ts)\ts\ot\ts\th_g\,\d_{c\bg}\,\d_{dh}
\ts-\ts
F_{cd\ts}(v-l\ts)\ts\ot\ts\th_g\,\th_h\,\d_{c\bg}\,x_{d\bh}\,\bigr)\ \times
\\[4pt]
\label{term421}
\bigl(\ts E_{eb}(v)\ot\th_i\,\th_g\,x_{\ts b\bg}\,\d_{e\bi}\ts\bigr)\,
\bigl(\ts\Et_{af}(w)\ot x_{\ts fh}\,\d_{aj}\ts\bigr)
\end{gather}
plus the sum over $k$ of
\begin{gather}
\nonumber
\bigl(\,\de_{ef}\,E_{eb}(v)\,\Et_{af}(w)\ot1\ts\bigr)\ \times
\\[4pt]
\nonumber
\bigl(\ts
-\,
F_{-c,-d\ts}(v-l\ts)
\ts\ot\ts\th_i\,
\th_k\,x_{\ts b\bk}\,x_{ck}\,
\d_{d\bi}\,\d_{aj}
\ts\,\mp\,\ts
F_{-c,d\ts}(v-l\ts)
\ts\ot\ts 
\th_k\,x_{\ts b\bk}\,x_{ck}\,
x_{di}\,\d_{aj}
\\[4pt]
\label{term422}
+\,
F_{c,-d\ts}(v-l\ts)
\ts\ot\ts
\th_i\,
x_{\ts bk}\,\d_{ck}\,
\d_{d\bi}\,\d_{aj}
\ts\,\pm\,\ts
F_{cd\ts}(v-l\ts)
\ts\ot\ts 
x_{\ts bk}\,\d_{ck}\,x_{di}\,\d_{aj}\ts\bigl)\,.
\end{gather}

By Proposition \ref{prop3}, the action of the 
sum in the last of the six lines \eqref{fivlin}
on the elements of \eqref{1vu} 
coincides with the action of the sum over $a\com b\com c\com d$ of
$$
\begin{aligned}
&-\,(\ts 
F_{\ts b,-c}\,E_{ab}(v)\,F_{-c,-d\ts}(v-l\ts)
\ts+\ts
F_{\ts bc}\,E_{ab}(v)\,F_{c,-d\ts}(v-l\ts))
\ts\ot\ts\th_i\,\d_{a\bi}\,\d_{dj}
\\[6pt]
&-\,(\ts 
F_{\ts b,-c}\,E_{ab}(v)\,F_{-c,d\ts}(v-l\ts)
\ts+\ts
F_{\ts bc}\,E_{ab}(v)\,F_{cd\ts}(v-l\ts))
\ts\ot\ts\th_i\,\th_j\,\d_{a\bi}\,x_{d\bj}\,.
\end{aligned}
$$
Modulo $\q\cdot M\ts$, the expression in the above two lines acts
on the elements of \eqref{1vu} as the sum over $k$ of
$$
\bigl(\ts E_{ab}(v)\ot1\ts\bigr)\ \times
$$
\begin{gather*}
\bigl(\ts-\,
F_{-c,-d\ts}(v-l\ts)
\ts\ot\ts
\th_i\,\th_k\,x_{\ts b\bk}\,x_{ck}\,
\d_{a\bi}\,\d_{dj}
\,+\,
F_{c,-d\ts}(v-l\ts)
\ts\ot\ts
\th_i\,
x_{\ts bk}\,\d_{ck}\,
\d_{a\bi}\,\d_{dj}
\\[4pt]
-\,
F_{-c,d\ts}(v-l\ts)
\ts\ot\ts
\th_k\,x_{\ts b\bk}\,x_{ck}\,
\th_i\,\th_j\,\d_{a\bi}\,x_{d\bj}
\,+\,
F_{cd\ts}(v-l\ts)
\ts\ot\ts
x_{\ts bk}\,\d_{ck}\,\th_i\,\th_j\,\d_{a\bi}\,x_{d\bj}\ts\bigr)\ts.
\end{gather*}
Note that this sum over the index $k$ 
can be rewritten as the sum over $k$ of
\begin{gather}
\nonumber
\bigl(\ts
F_{-c,-d\ts}(v-l\ts)
\ts\ot\ts
x_{ck}\,\d_{dj}
\,+\,
F_{-c,d\ts}(v-l\ts)
\ts\ot\ts
\th_j\,
x_{ck}\,\,x_{d\bj}
\\[4pt]
\nonumber
-\,
F_{c,-d\ts}(v-l\ts)
\ts\ot\ts
\th_k\,
\d_{c\bk}\,\,\d_{dj}
\,-\,
F_{cd\ts}(v-l\ts)
\ts\ot\ts
\th_k\,\th_j
\d_{c\bk}\,x_{d\bj}\ts\bigr)\ \times
\\[4pt]
\label{term5}
\bigl(\ts-\,E_{ab}(v)
\ts\ot\ts 
\th_i\,\th_k\,
x_{\ts b\bk}\,\d_{a\bi}\ts\bigr)\,.
\end{gather}

Consider the sum of the expressions \eqref{term422} 
over the running indices $e\com f\ts$. Add this sum to
the expression displayed in the five lines \eqref{term2}. Using the relation
\begin{equation}
\label{can3}
\sum_{e}\,\Et_{\ts eb}(v)\,\Et_{ae}(w)\,=\,
\frac{\ts\Et_{ab}(w)-\Et_{ab}(v)\ts}{2u}
\end{equation}
together with \eqref{eet} and performing cancellations, 
we get the expression
\begin{gather*}
\bigl(\,
\mp\,\ts
F_{-d,c\ts}(v-l\ts)
\ts\ot\ts 
\th_k\,x_{\ts b\bk}\,x_{ck}\,
x_{di}\,\d_{aj}
\,-\,
F_{-d,-c\ts}(v-l\ts)
\ts\ot\ts x_{\ts bk}\,\d_{ck}\,x_{di}\,\d_{aj}
\\[4pt]
F_{d,-c\ts}(v-l\ts)
\ts\ot\ts
\th_i\,
x_{\ts bk}\,\d_{ck}\,
\d_{d\bi}\,\d_{aj}
\,\pm\,
F_{dc\ts}(v-l\ts)
\ts\ot\ts\th_i\,
\th_k\,x_{\ts b\bk}\,x_{ck}\,
\d_{d\bi}\,\d_{aj}\ts\bigr)\ \times
\\[4pt]
\bigl(\ts(1-Z(v))\,\Et_{ab}(w)\ot1\ts\bigr)\,.
\end{gather*}

After exchanging the running indices $c$ and $d\ts$,
the sum over the index $k$ of the expressions in the last three 
displayed lines 
can be rewritten as 
\begin{equation}
\label{term24221}
-\,\,\de_{cd}\,(1-Z(v))\,
\bigl(\ts F_{-c,-d\ts}(v-l\ts)+F_{cd\ts}(v-l\ts)\bigr)\,
\Et_{ab}(w)
\ts\ot\ts 
x_{\ts bi}\,\d_{aj}
\end{equation}
plus the sum over $k$ of
\begin{gather}
\nonumber
\bigl(\ts
F_{-c,-d\ts}(v-l\ts)
\ts\ot\ts
x_{ci}\,\d_{dk}
\,+\,
F_{-c,d\ts}(v-l\ts)
\ts\ot\ts
\th_k\,
x_{ci}\,x_{d\bk}
\\[4pt]
\nonumber
-\,
F_{c,-d\ts}(v-l\ts)
\ts\ot\ts
\th_i\,
\d_{c\bi}\,\d_{dk}
\,-\,
F_{cd\ts}(v-l\ts)
\ts\ot\ts
\th_i\,\th_k
\d_{c\bi}\,x_{d\bk}
\ts\bigr)\ \times
\\[4pt]
\label{term24222}
\bigl(\ts-\,(1-Z(v))\,\Et_{ab}(w)
\ts\ot\ts 
x_{\ts bk}\,\d_{aj}\ts\bigr)\,.
\end{gather}
Here we again used the commutation relations in the ring
$\PD\ts(\CC^{\ts m+l}\ot\CC^{\ts n})\ts$.

Let us now perform the summation over all running indices
in the 
expressions 
\eqref{term3},\eqref{term421},\eqref{term5},\eqref{term24222}
and then take their total. By exchanging the running indices
$b$ and $f$ in \eqref{term421}, and by replacing the running index $k$
in \eqref{term5},\eqref{term24222} by $g\com h$ respectively, 
the total can be written as the sum over indices $c\com d$ and
$g\com h$ of
\begin{gather}
\nonumber
\bigl(\ts(1-Z(v))\ot1\ts\bigr)\ \times
\\[4pt]
\nonumber
\bigl(\,
F_{-c,-d\ts}(v-l\ts)\ts\ot\ts x_{cg}\,\d_{dh}
\,+\,
F_{-c,d\ts}(v-l\ts)\ts\ot\ts\th_h\,x_{cg}\,x_{d\bh}
\\[4pt]
\nonumber
-\,\ts
F_{c,-d\ts}(v-l\ts)\ts\ot\ts\th_g\,\d_{c\bg}\,\d_{dh}
\,-\,
F_{cd\ts}(v-l\ts)\ts\ot\ts\th_g\,\th_h\,\d_{c\bg}\,x_{d\bh}\,\bigr)\ \times
\\[4pt]
\label{cdgh}
\bigl(\ts\de_{ig}-
\sum_{e,f}\,
\Et_{ef}(v)\ot\th_i\,\th_g\,x_{\ts f\bg}\,\d_{e\bi}\,\bigr)\,
\bigl(\ts\de_{hj}-
\sum_{a,b}\,
\Et_{ab}(w)\ts\ot\ts x_{\ts bh}\,\d_{aj}\ts\bigr)\,.
\end{gather}

Let us perform the summation in \eqref{term410} 
over the running indices $b\com e\ts$. Then let us replace
the running index $f$ by the index $b\ts$, which becomes free 
after the summation. By adding the resulting sum to 
the expession \eqref{term1} we get
$$
(1-Z(v))\ts
\bigl(\ts W(v-l\ts)-1\ts\bigr)\,\Et_{ab}(w)
\ts\ot\ts x_{\ts bi}\,\d_{aj}
$$
due to \eqref{eet} and \eqref{can3}.
By performing the summation in \eqref{term24221}
over the running indices $c\com d$ and then adding the result
to the last displayed expression, we get
\begin{equation}
\label{term000}
-\,(1-Z(v))\,\Et_{ab}(w)
\ts\ot\ts x_{\ts bi}\,\d_{aj}\,.
\end{equation}

Now perform the summation over all running indices
in the two expressions \eqref{term41},\eqref{term000} and then
add the two resulting sums to \eqref{term0}. By using the relation
\eqref{eet} once again, the total can be written as
the sum over the index $k$ of
\begin{gather}
\nonumber
\bigl(\ts(1-Z(v))\ot1\ts\bigr)\ \times
\\[4pt]
\label{k}
\bigl(\ts\de_{ik}-
\sum_{e,f}\,
\Et_{ef}(v)\ot\th_i\,\th_k\,x_{\ts f\bk}\,\d_{e\bi}\,\bigr)\,
\bigl(\ts\de_{kj}-
\sum_{a,b}\,
\Et_{ab}(w)\ts\ot\ts x_{\ts bk}\,\d_{aj}\ts\bigr)\,.
\end{gather}

By using the definition of the series $\Et_{ab}(v)$ 
as given before the relation \eqref{eet},
\begin{gather*}
\textstyle
\Et_{ef}(v)=(\ts v-l-\Ep)^{-1}_{fe}=-\,(-\ts u\mp\frac12-m+\Ep)^{-1}_{fe}\,,
\\[4pt]
\textstyle
\Et_{ab}(w)=(\ts w-l-\Ep)^{-1}_{ba}=-\,(\ts u\mp\frac12-m+\Ep)^{-1}_{ba}\,.
\end{gather*}
We also used the definitions \eqref{v} and \eqref{w}.
Hence the sum of the expressions \eqref{cdgh} over 
the indices $c\com d$ and $g\com h$ plus the sum of the expessions 
\eqref{k} over the index $k$ can be rewritten
as the sum over the indices $g\com h$ of the series in $u^{-1}$,
\begin{gather*}
\textstyle
\bigl(\ts(1-Z(\ts u\pm\frac12+m+l\ts))\ot1\ts\bigr)\ \times
\\[8pt]
\bigl(\,\de_{gh}\ \ts+
\sum_{c,d}\ 
\bigl(\,
\textstyle
F_{-c,-d\ts}(\ts u\pm\frac12+m\ts)
\ts\ot\ts 
x_{cg}\,\d_{dh}
\,+\,
\textstyle
F_{-c,d\ts}(\ts u\pm\frac12+m\ts)
\ts\ot\ts
\th_h\,x_{cg}\,x_{d\bh}
\\
-\,\ts
\textstyle
F_{c,-d\ts}(\ts u\pm\frac12+m\ts)
\ts\ot\ts
\th_g\,\d_{c\bg}\,\d_{dh}
\,-\,
\textstyle
F_{cd\ts}(\ts u\pm\frac12+m\ts)
\ts\ot\ts
\th_g\,\th_h\,\d_{c\bg}\,x_{d\bh}\,\bigr)\bigr)\ \times
\\[8pt]
\bigl(\,\de_{ig}+
\sum_{e,f}\,
\textstyle
(-\ts u\mp\frac12-m+\Ep)^{-1}_{fe}
\displaystyle
\ot\th_i\,\th_g\,x_{\ts f\bg}\,\d_{e\bi}\,\bigr)\ \times
\end{gather*}
\begin{gather*}
\bigl(\,\de_{hj}+
\sum_{a,b}\,
\textstyle
(\ts u\mp\frac12-m+\Ep)^{-1}_{ba}
\displaystyle
\ts\ot\ts x_{\ts bh}\,\d_{aj}\ts\bigr)
\end{gather*}
with coefficients in the algebra
$\U(\ts\f_{m}\op\gl_{\ts l})\ot\PD\ts(\CC^{\ts m+l}\ot\CC^{\ts n})\ts$.
By mapping the series $S_{ij}(u)$ to this sum 
we describe the action of the extended twisted Yangian 
$\X(\g_n)$ on the subspace \eqref{1vu} modulo $\q\cdot M$.
By comparing this sum with the product of the series
\eqref{baser} and \eqref{zuser} where $z=m\pm\frac12\,$,
we now prove that the map $\ph$  intertwines the actions of $\X(\g_n)\ts$;
here we use \eqref{ehom} and \eqref{fhom}.
This completes the proof of Theorem \ref{parind}.


\section{Zhelobenko operators}
\setcounter{section}{4}
\setcounter{equation}{0}
\setcounter{theorem*}{0}

Let us consider the \textit{hyperoctahedral group\/} $\H_m\ts$. This
is the semidirect product of the symmetric group $\Sym_m$ and
the Abelian group $\ZZ_2^m\ts$, where $\Sym_m$ acts by permutations
of the $m$ copies of $\ZZ_2\ts$. 
In this section, we assume that $m>0\ts$.
The group $\H_m$ is generated by the elements $\si_a$ with $a=1\lcd m\ts$.
The elements $\si_a$ with the indices $a=1\lcd m-1$ are elementary
transpositions generating the symmetric group $\Sym_m\ts$, so that
$\si_a=(a\com a+1)\ts$. Then $\si_m$ is the generator of the
$m\ts$-th factor $\ZZ_2$ of $\ZZ_2^m\ts$.
The elements 
$\si_1\lcd\si_m\in\H_m$ are involutions and satisfy the braid relations
\begin{align*}
\si_a\,\si_{a+1}\,\si_a
&\,=\,
\si_{a+1}\,\si_{a}\,\si_{a+1}
\!\!\!\quad\quad\textrm{for}\ \quad
a=1\lcd m-2\ts;
\\
\si_a\,\si_{\ts b}
&\,=\,
\si_{\ts b}\,\si_a
\hspace{34pt}
\ \quad\textrm{for}\ \quad
a=1\lcd b-2\ts;
\\
\si_{m-1}\,\si_m\,\si_{m-1}\,\si_m
&\,=\,
\si_{m}\,\si_{m-1}\,\si_{m}\,\si_{m-1}.
\end{align*}

Note that $\H_m$ is the Weyl group of the simple Lie algebra $\sp_{2m}\ts$.
Let $\Hh_m$ be the braid group corresponding to $\sp_{2m}\ts$. 
It is generated by the elements $\sih_1\lcd\sih_m$ 
which by definition satisfy the above displayed relations, instead of
the involutions $\si_1\lcd\si_m$ respectively. For any reduced decomposition
$\si=\si_{a_1}\ldots\si_{a_K}$ in $\H_m$ put
\begin{equation}
\label{sih}
\sih=\sih_{a_1}\ldots\,\sih_{a_K}.
\end{equation}
The definition of $\,\sih$
is independent of the choice of a reduced decomposition~of~$\si\ts$.

The hyperoctahedral group $\H_m$ also contains
the Weyl group of the reductive 
Lie algebra $\so_{2m}$ as a subgroup of index two.
Denote this subgroup by $\Hp_m\ts$, 
it is generated by the elementary transpositions
$\si_1\lcd\si_{m-1}$ and by the involution
$\sip_m=\si_m\,\si_{m-1}\,\si_m\ts$. Along with the braid relations
between $\si_1\lcd\si_{m-1}\ts$, we also have the braid relations
involving $\sip_m\ts$, 
\begin{align*}
\si_a\,\sip_{\ts m}
&\,=\,
\sip_{\ts m}\,\si_a
\ \quad\textrm{for}\ \quad
a=1\lcd m-3\com m-1\ts;
\\
\si_{m-2}\,\sip_m\,\si_{m-2}
&\,=\,
\sip_m\,\si_{m-2}\,\sip_m\ts.
\end{align*}
When $m>1\ts$, the braid group corresponding to $\so_{2m}$ 
is generated by $m$ elements
satisfying the same braid relations instead of the
$m$ involutions $\si_1\lcd\si_{m-1}\com\sip_m$ respectively.
When $m=1\ts$, the braid group corresponding to $\f_m=\so_2$ consists
of the identity element only.

Now let the indices
$c\com d$ run through $-m\lcd-1\com1\lcd m\ts$.
For $c>0$ put $\cb=m+1-c\ts$; for $c<0$ put $\cb=-\,m-1-c\ts$. 
Consider a representation $\si\mapsto\sib$ of the
group $\H_m$ by permutations of $-m\lcd-1\com1\lcd m$ such that
\begin{equation}
\label{sib}
\sib\ts(c)=\overline{\si(\ts\cb\ts)\!}\,
\quad\text{for}\quad\si\in\Sym_m
\end{equation}
and 
$\sib_m\ts(c)=-c$ if $|c|=1\ts$, 
while
$\sib_m\ts(c)=c$ if $|c|>1\ts$. 
We can define an action of the braid group $\Hh_m$ by automorphisms of
the Lie algebra $\f_m\ts$, by the assignments
\begin{align}
\label{siact}
\sih:\,F_{cd}
&\,\mapsto\,
F_{\sib(c)\ts\sib(d)}
\quad\text{for}\quad\si\in\Sym_m\ts,
\\[2pt]
\label{simact}
\sih_m:\,F_{cd}
&\,\mapsto\,
(\ts\mp\ts1\ts)^{\ts\de_{c1}\ts+\,\de_{d1}}F_{\sib_m(c)\ts\sib_m(d)}\,.
\end{align}
Here the upper sign in $\mp$ corresponds to the case $\f_m=\sp_{2m}\ts$,
while the lower sign corresponds to $\f_m=\so_{2m}\ts$.
The automorphism property can be checked
by using the relations \eqref{ufmrel}, see also
the proof of Part~(i) of Lemma~\ref{lemma41} below.
This action of the group $\Hh_m$ on $\f_m$
extends to an action of $\Hh_m$ by automorphisms of
the associative algebra $\U(\ts\f_m)\ts$. 
Note that in the case $\f_m=\so_{2m}$
the action of $\Hh_m$ on $\U(\f_m)$
factors to an action of the group $\H_m\ts$.

Further, one can define an action
of the braid group $\Hh_m$ by automorphisms of the algebra 
$\PD\ts(\CC^{\ts m}\ot\CC^{\ts n})$ in the following way. Put
\begin{align}
\nonumber
\sih\ts(x_{ai})&=x_{\ts\sib(a)\ts i}
\hspace{-16pt}
&
\textrm{and}
&
&
\hspace{-16pt}
\sih\ts(\d_{ai})&=\d_{\ts\sib(a)\ts i}
\hspace{-16pt}
&
\textrm{for}
&
&
&
\hspace{-16pt}
\si\in\Sym_m\,,
\\
\nonumber
\sih_m(x_{ai})&=x_{ai}
\hspace{-16pt}
&
\textrm{and}
&
&
\hspace{-16pt}
\sih_m(\d_{ai})&=\d_{ai}
\hspace{-16pt}
&
\textrm{for}
&
&
&
\hspace{-16pt}
a>1\ts,
\\
\label{fourier}
\sih_m\ts(x_{1i})&=-\,\th_i\,\d_{\ts1\bi}
\hspace{-16pt}
&
\textrm{and}
&
&
\hspace{-16pt}
\sih_m\ts(\d_{\ts1i})&=\th_i\,x_{1\bi}
\hspace{-16pt}
&
\phantom{\textrm{for}}
&
&
&
\hspace{-16pt}
\end{align}
where $i=1\lcd n\ts$.
Note that in the case $\f_m=\so_{2m}$ the element 
$\sih_m^{\,2}\in\Hh_m$
acts on $x_{1i}$ and on $\d_{\ts1i}$ as the identity, so that
the action of $\Hh_m$ on $\PD\ts(\CC^{\ts m}\ot\CC^{\ts n})$
factors to an action of the group $\H_m\ts$. But in the case
$\f_m=\sp_{2m}$ the element $\sih_m^{\,2}$
acts on $x_{1i}$ and on $\d_{\ts1i}$ as minus the identity,
because $\th_i=1$ in this case. This is
why we use the braid group, rather than the Weyl group $\H_m$ 
of the simple Lie algebra~$\sp_{2m}\ts$.
Taking the tensor product of the actions of $\Hh_m$
on the algebras $\U(\ts\f_m)$ and $\PD\ts(\CC^{\ts m}\ot\CC^{\ts n})\ts$,
we get an action of $\Hh_m$ by automorphisms
of the algebra $\B_m=\U(\ts\f_m)\ot\PD\ts(\CC^{\ts m}\ot\CC^{\ts n})\ts$.

\begin{lemma*} 
\label{lemma41}
{\rm\,\,(i)} 
The map $\zeta_n:\U(\ts\f_m)\to\PD\ts(\CC^{\ts m}\ot\CC^{\ts n})$
is\/ $\Hh_m\ts$-equivariant.
\\
{\rm(ii)}
The action of\/ $\Hh_m$ on $\B_m$ leaves invariant
any element of the image of\/ $\X(\g_n)$
under the homomorphism $\be_m$.
\end{lemma*}

\begin{proof}
Let us employ the elements $p_{ci}$ and $q_{ci}$ 
of the algebra $\PD\ts(\CC^{\ts m}\ot\CC^{\ts n})\ts$, 
introduced immediately after stating Proposition~\ref{xb}.
In terms of these elements, the action of $\Hh_m$ on the algebra
$\PD\ts(\CC^{\ts m}\ot\CC^{\ts n})$ can be described by setting
\begin{equation*}
\begin{aligned}
\sih\ts(\ts p_{ci})&=p_{\ts\sib(c)\ts i}
&
\textrm{and}&
&
\sih\ts(q_{ci})&=q_{\ts\sib(c)\ts i}
\quad\textrm{for}\quad
\si\in\Sym_m\,,
\\
\sih_m\ts(\ts p_{ci})&=(\ts\mp\ts1\ts)^{\ts\de_{c1}}p_{\ts\sib_m(c)i}
&
\textrm{and}&
&
\sih_m\ts(q_{ci})&=(\ts\mp\ts1\ts)^{\ts\de_{c1}}q_{\ts\sib_m(c)i}
\end{aligned}
\end{equation*}
for every 
$c=-m\lcd-1\com1\lcd m\ts$.
Part (i) follows
by comparing our definition of the action of $\Hh_m$ on $\f_m$
with the description \eqref{ganpq} of the homomorphism~$\zeta_n\ts$.
Part~(ii) follows similarly, using 
the description \eqref{fhompq} of 
$\be_m\ts$.
\end{proof}

Consider the Cartan subalgebra $\h$ from
the triangular decomposition \eqref{tridec}.
In the notation of this section,
our chosen basis of $\h$ is $(\ts F_{-\ab,-\ab}\,|\,a=1\lcd m\ts)\ts$.
Now let $(\ts\ep_a\,|\,a=1\lcd m\ts)\subset\h^\ast$ be the dual basis,
so that $\ep_b\ts(F_{-\ab,-\ab})=\de_{ab}\ts$.
For $c<0$ put $\ep_c=-\ts\ep_{-c}\ts$.
Thus the element $\ep_c\in\h^\ast$ is defined for every index
$c=-\ts m\lcd-1\com1\lcd m\ts$.

Consider the root system of the Lie algebra $\f_m$ in $\h^\ast\ts$. Put
\begin{equation*}
\eta_{\ts a}=\ep_a-\ep_{a+1}
\ \quad\textrm{for}\ \quad
a=1\lcd m-1\ts.
\end{equation*}
Also put $\eta_m=2\ts\ep_m$ in the case $\f_m=\sp_{2m}\ts$, and
$\eta_m=\ep_{m-1}+\ep_m$ in the case $\f_m=\so_{2m}$.
Then $\eta_1\lcd\eta_m$ are the \textit{simple roots\/} of $\f_m\ts$.
Denote by $\De^+$ 
the set of \textit{positive roots\/} of $\f_m\ts$.
These are the weights $\ep_a-\ep_b$ and
$\ep_a+\ep_b$ where $1\le a<b\le m$ in the case $\f_m=\so_{2m}\ts$,
and the same weights together with $2\ts\ep_a$
where $1\le a\le m$ in the case $\f_m=\sp_{2m}\ts$.
We assume that in the case $\f_m=\so_2$ the root system of $\f_m$ is empty.
Let $\rho$ be halfsum of positive roots of $\f_m\ts$,
so that its sequence of labels
$(\ts\rho_1\lcd\rho_m\ts)$ is $(\ts m\lcd1\ts)$
in the case $\f_m=\sp_{2m}\ts$, and is $(\ts m-1\lcd0\ts)$
in the case $\f_m=\so_{2m}\ts$.
For each $a=1,...,m-1$ put
\begin{equation}
\label{Fc}
E_a=F_{-\ab,-\overline{a+1}}\,,\quad
F_a=F_{-\overline{a+1},-\ab}\,,\quad
H_a= F_{-\ab,-\ab}-F_{-\overline{a+1},-\overline{a+1}}\,\ts.
\end{equation}
Put
\begin{equation}
\label{Fm}
E_m=F_{-\overline{m},\overline{m}}\,/\ts2\ts,\quad
F_m=F_{\overline{m},-\overline{m}}\,/\ts2\ts,\quad
H_m=F_{-\overline{m},-\overline{m}}
\end{equation}
in the case $\f_m=\sp_{2m}\ts$. In the case when 
$\f_m=\so_{2m}$ with $m>1\ts$, put
\begin{equation}
\label{Fmm}
E_m=F_{-\overline{m-1},\overline{m}}\,,\quad
F_m=F_{\overline{m},-\overline{m-1}}\,,\quad
H_m=F_{-\overline{m-1},-\overline{m-1}}+F_{-\overline{m},-\overline{m}}\,.
\end{equation}
For every possible index 
$a$ the three elements $E_a\ts,F_a\ts,H_a$ of the Lie algebra $\f_m$
span a subalgebra isomorphic to $\mathfrak{sl}_2\ts$.
They satisfy the commutation relations
\begin{equation}
\label{sltwo}
[\ts E_a\com F_a\ts]=H_a\ts,
\quad
[\ts H_a\com E_a\ts]=2\ts E_a\ts,
\quad
[\ts H_a\com F_a\ts]=-\ts2\ts F_a\ts.
\end{equation}

So far we denoted by $\B_m$ the associative algebra
$\U(\ts\f_m)\ot\PD\ts(\CC^{\ts m}\ot\CC^{\ts n})\ts$.
Let us now use a different presentation
of the same algebra. Namely, from now until the end of the next section,
on we will regard $\B_m$ 
as the associative algebra generated by the algebras $\U(\ts\f_m)$ 
and $\PD\ts(\CC^{\ts m}\ot\CC^{\ts n})$ with the cross relations
\begin{equation}
\label{defar} 
[\ts X\com Y\ts]=[\ts\zeta_n(X)\com Y\ts]
\end{equation}
for any $X\in\f_m$ and $Y\in\PD\ts(\CC^{\ts m}\ot\CC^{\ts n})\ts$.
The brackets at the left hand side of the relation
\eqref{defar} denote the commutator in $\B_m\ts$, 
while the brackets at the right hand side denote the commutator in
the algebra $\PD\ts(\CC^{\ts m}\ot\CC^{\ts n})$ embedded to $\B_m\ts$.
In particular, we will regard
$\U(\ts\f_m)$ as a subalgebra of $\B_m\ts$.
An isomorphism of this $\B_m$ with the 
tensor product $\U(\ts\f_m)\ot\PD\ts(\CC^{\ts m}\ot\CC^{\ts n})$ can be
defined by mapping the elements 
$X\in\f_m$ and $Y\in\PD\ts(\CC^{\ts m}\ot\CC^{\ts n})$ of $\B_m$ 
respectively to the elements
$$
X\ot1+1\ot\zeta_n(X)
\quad\textrm{and}\quad 
1\ot Y
$$
of $\U(\ts\f_m)\ot\PD\ts(\CC^{\ts m}\ot\CC^{\ts n})\ts$.
Here we use \eqref{gan}.
The action of the braid group $\Hh_m$ on $\B_m$ is defined
via its isomorphism of $\B_m$ with 
$\U(\ts\f_m)\ot\PD\ts(\CC^{\ts m}\ot\CC^{\ts n})\ts$. 
Since the map $\zeta_n$ is
$\Hh_m\ts$-equivariant, the same action of 
$\Hh_m$ is obtained by extending
the actions of $\Hh_m$ from the subalgebras
$\U(\ts\f_m)$ and $\PD\ts(\CC^{\ts m}\ot\CC^{\ts n})$ to $\B_m\ts$.

Now consider the following two sets of elements
of the algebra $\U(\h)\subset\U(\ts\f_m)\ts$:
\begin{gather}
\label{denset1}
\{\,F_{aa}-F_{\ts bb}+z\ts,\ F_{aa}+F_{\ts bb}+z
\ |\ 
1\leqslant a<b\leqslant m\ts,\ z\in\ZZ\,\}\,,
\\[2pt]
\label{denset2}
\{\,F_{aa}+z
\ |\ 
1\leqslant a\leqslant m\ts,\ z\in\ZZ\,\}\,.
\end{gather}
In the case $\f_m=\so_{2m}\ts$,
denote by $\Uhb$ the ring of fractions of the commutative
algebra $\U(\h)$ relative to the set of denominators \eqref{denset1} .
In the case $\f_m=\sp_{2m}\ts$,
denote by $\Uhb$ the ring of fractions of $\U(\h)$ 
relative to the union of sets \eqref{denset1} and \eqref{denset2}.
The elements of the ring $\Uhb$ can also
be regarded as rational functions on the vector space
$\h^\ast\ts$. The elements of the subalgebra $\U(\h)\subset\Uhb$
are then regarded as polynomial functions~on~$\h^\ast$.

Let us denote by $\Bb_m$ the ring of fractions of $\B_m$ 
relative to the same set of denominators as was used to
define the ring of fractions $\Uhb$. But now we regard these
denominators as elements of $\B_m$ using the 
embedding of $\h\subset\f_m$ into $\B_m\ts$. The ring $\Bb_m$ 
is defined due to the following relations
in $\B_m\ts$. For $c<0$ put $\ep_c=-\ts\ep_{-c}\ts$.
Thus the element $\ep_c\in\h^\ast$ is defined for every
$c=-m\lcd-1\com1\lcd m\ts$. Then for any element $H\in\h$ we have
\begin{gather*}
[\ts H\ts,F_{cd}\ts]=(\ts\ep_{\ts\db}-\ep_{\ts\cb}\ts)(H)\ts F_{cd}
\quad\text{for}\quad
c\com d=-m\lcd-1\com1\lcd m\ts;
\\[2pt]
[\ts H\ts,x_{ci}\ts]=
-\ts\ep_{\ts\cb}\ts(H)\,x_{ci}
\quad\text{and}\quad
[\ts H\ts,\d_{\ts ci}\ts]=
\ep_{\ts\cb}\ts(H)\,\d_{\ts ci}
\quad\text{for}\quad
c=1\lcd m\ts.
\end{gather*}
So the ring $\B_m$ obeys the Ore condition relative
to our set of denominators.
Using left multiplication by elements of
$\,\overline{\!\U(\h)\!\!\!}\,\,\,$,
the ring 
$\Bb_m$ becomes a module~of~$\Uhb\ts$.

The ring $\Bb_m$ is also an associative algebra over $\CC\ts$.
The action of the braid group $\Hh_m$ on $\B_m$
preserves the set of denominators, so that $\Hh_m$ also
acts by automorphisms of the algebra $\Bb_m\ts$.
Using 
\eqref{Fc} and \eqref{Fm} when $\f_m=\sp_{2m}\ts$, or 
\eqref{Fc} and \eqref{Fmm} when $\f_m=\so_{2m}\ts$, 
for every simple root $\eta_{\ts a}$ of $\f_m$ 
define a linear map $\xi_{\ts a}:\B_m\to\Bb_m$ by setting
\begin{equation}
\label{q1}
\xi_{\ts a}(Y)=
Y+\,\sum_{s=1}^\infty\,\,
(\ts s\ts!\,H_a^{\ts(s)}\ts)^{-1}\ts E_a^{\ts s}\,
\widehat{F}_a^{\ts s}(\ts Y)\ts.
\end{equation}
Here
$$
H_a^{\ts(s)}=(H_a)(H_a-1)\cdots(H_a-s+1)
$$
and $\widehat{F}_a$ is the operator of adjoint action
corresponding to the element $F_a\in\B_m\ts$,
$$
\widehat{F}_a(\ts Y)=[\ts F_a\ts\com Y\ts]\ts.
$$
For a given element $Y\in\B_m$ only finitely many terms of the sum
\eqref{q1} differ from zero. In the case $\f_m=\so_{2}$ there is
no roots of $\f_m\ts$, and no corresponding operators $\B_m\to\Bb_m\ts$.
On the other hand, in the case when $\f_m=\so_{2m}$ with $m>1\ts$, 
by \eqref{simact}
\begin{equation}
\label{xim}
\xi_m\,\sih_m=\sih_m\,\xi_{m-1}\ts,
\end{equation}
because
$$
\sih_m:\,
E_{m-1}\mapsto E_m\ts,\ 
F_{m-1}\mapsto F_m\ts,\ 
H_{m-1}\mapsto H_m\ts.
\vspace{2pt}
$$

Let $\J$ and $\Jb$ be the right ideals of algebras $\B_m$ and $\Bb_m$
respectively, generated by all elements of the subalgebra $\n\subset\f_m\ts$.
The following four properties of the linear operator $\xi_{\ts a}$
go back to \cite[Section~2]{Z1}.
For any elements $X\in\h$ and $Y\in\B_m\ts$, 
\begin{gather}
\label{q1-1}
\xi_{\ts a}(E_{\ts a}\ts Y)
\,\in\,
\ts\Jb\ts,
\\[2pt]
\label{q10}
\xi_{\ts a}(Y\ts F_{\ts a})
\,\in\,
\ts\Jb\ts,
\\[2pt]
\label{q11}
\xi_{\ts a}(X\ts Y)
\,\in\,
(\ts X+\eta_{\ts a}\ts(X))\,\ts\xi_{\ts a}(\ts Y)\ts+\ts\Jb\ts,
\\[2pt]
\label{q12}
\xi_{\ts a}(\ts Y X)
\,\in\,
\,\xi_{\ts a}(\ts Y)\ts(\ts X+\eta_{\ts a}\ts(X))\ts+\ts\Jb\ts.
\end{gather}
See \cite[Section 3]{KN1} for detailed proofs of these four properties.
They use only the commutation relations \eqref{sltwo},
not the actual form of elements $E_a\ts,F_a\ts,H_a\ts$. 

The property \eqref{q11} allows us to define a linear map
$
\bar\xi_{\ts a}:\Bb_m\to\Jb\,\ts\backslash\,\Bb_m
$
by 
\begin{equation}
\label{N14}
\bar\xi_{\ts a}(X\,Y)=Z\,\xi_{\ts a}(\ts Y)\ts+\ts\Jb
\quad\text{for}\quad
X\in\Uhb
\quad\text{and}\quad
Y\in\B_m\ts,
\end{equation}
where the element $Z\in\Uhb$ is defined by the equality
\begin{equation}
\label{N1414}
Z(\mu)=X(\ts\mu+\eta_{\ts a})
\quad\text{for}\quad
\mu\in\h^\ast
\end{equation}
when both $X$ and $Z$ are regarded as rational functions on $\h^\ast\ts$.
The backslash in $\Jb\,\ts\backslash\,\Bb_m$ indicates that the quotient
is taken relative to a \textit{right\/} ideal of $\Bb_m\ts$.

\begin{proposition*}
\label{p2}
For any simple root\/ $\eta_{\ts a}\!$ of\/ $\f_m$ we have
the inclusion\/ 
$
\sih\ts(\ts\Jb\ts)\ts\subset\ts\ker\ts\bar\xi_{\ts a}
$
where $\si=\si_a$ unless\/ $\f_m=\so_{2m}$ and\/ $a=m$,
in which case $\si=\sip_m\ts$.
\end{proposition*}

\begin{proof}
By definition, the nilpotent subalgebra $\n\subset\f_m$ is spanned by 
the elements $F_{cd}$ with $c>d\ts$. Due to the relation
$F_{cd}=-\,\ep_{cd}\,F_{-d,-c}$ the subalgebra $\n$ is also
spanned by the elements $F_{cd}$ with $c>d\ts$ and $d<0\ts$.
Now for any $a=1\lcd m$ denote by
$\n_{\ts a}$ the vector subspace of $\f_m$ spanned by all 
the elements $F_{cd}$ with $c>d\ts$ and $d<0\ts$, except the element
$F_a\ts$; see the definitions \eqref{Fc},\eqref{Fm} and \eqref{Fmm}.
Then $\sih\,(\ts\Jb\ts)$ above
is the right ideal of $\Bb_m$ generated
by the subspace $\n_{\ts a}\subset\f_m$ and by 
the element $E_a\ts$.
Due to \eqref{q1-1}, it is sufficient to prove that 
$\bar\xi_{\ts a}\ts(X\ts Y)\in\Jb\ts$ for 
any $X\in\n_{\ts a}$ and $Y\in\Bb_m\ts$.
Observe that
the subspace $\n_{\ts a}\subset\f_m\ts$ is preserved by the
adjoint action of the elements $E_a,F_a$ and $H_a$ on $\f_m\ts$.
Therefore $\bar\xi_{\ts a}\ts(X\ts Y)\in\Jb\ts$ for 
any $X\in\n_{\ts a}$ and $Y\in\Bb_m$ 
by the definition of $\bar\xi_{\ts a}\ts$.
\end{proof}

Recall that $\np$ denotes the nilpotent subalgebra of $\f_m$ spanned by 
all the elements $F_{cd}$ with $c<d\ts$. Due to the relation
$F_{cd}=-\,\ep_{cd}\,F_{-d,-c}$ the subalgebra $\np$ is also
spanned by the elements $F_{cd}$ with $c<d$ and $c<0\ts$.
Now for any $a=1\lcd m$ denote by
$\np_{\ts a}$ the vector subspace of $\f_m$ spanned by all 
the elements $F_{cd}$ with $c<d$ and $c<0\ts$, except the element $E_a\ts$. 

Denote by $\Jp$ the left ideal of $\B_m\ts$, generated by the
elements $X-\zeta_n(X)$ where $X\in\np\ts$. Under the
isomorphism of $\B_m$ with 
$\U(\ts\f_m)\ot\PD\ts(\CC^{\ts m}\ot\CC^{\ts n})\ts$,
for any $X\in\f_m$
the difference $X-\zeta_n(X)\in\B_m$ is mapped to the element 
\begin{equation}
\label{differ}
X\ot1\in\U(\ts\f_m)\ot1
\subset
\U(\ts\f_m)\ot\PD\ts(\CC^{\ts m}\ot\CC^{\ts n})\ts.
\end{equation}
Let $\Jp_{\ts a}$ be
the left ideal of $\B_m\ts$, generated by the elements
$X-\zeta_n(X)$ with $X\in\np_{\ts a}\ts$, and by the element $E_a\in\B_m\ts$.
Denote $\Jpb=\Uhb\ts\Jp$ and 
$\Jpb_{\ts a}=\Uhb\ts\Jp_{\ts a}\ts$. 
Then both $\Jpb$ and $\Jpb_{\ts a}$ are left ideals of the algebra $\Bb_m\ts$.

\begin{proposition*}
\label{prop3N}
For any simple root\/ $\eta_{\ts a}\!$ of\/ $\f_m$ we have\/
$
\bar\xi_{\ts a}\ts(\,\sih\,(\ts\Jpb_{\ts a}\ts))
\ts\subset\ts\Jpb\ns+\Jb
$
where $\si=\si_a$ unless $\f_m=\so_{2m}$ and $a=m\ts$,
in which case $\si=\sip_m\ts$.
\end{proposition*}

\begin{proof} 
The left ideal $\sih\,(\ts\Jpb_{\ts a}\ts)$ is generated by the
element $F_a\in\B_m\ts$, and by the subspace formed by
all elements $X-\zeta_n(X)$ with $X\in\np_{\ts a}\ts$.
Observe that the latter subspace in $\B_m$ is preserved by the
adjoint action of the element $F_a\in\B_m\ts$.
Hence for any $Y\in\Bb_m$ and any $Z$ from that subspace,
$\bar{\xi}_{\ts a}\ts(\ts YZ)\in\ts\Jpb\ns+\ts\Jb$ by the definition
of the operator $\bar{\xi}_{\ts a}\ts$. Further,
$\bar{\xi}_{\ts a}\ts(\ts YF_a)\in\ts\Jb$
for any $Y\in\bar{\B}_m\ts$ due to \eqref{q10}.
\end{proof}

Proposition \ref{p2} 
allows us to define for any simple root $\eta_{\ts a}$ 
a linear map
\begin{equation*}
\xic_{\ts a}:\,\Jb\,\ts\backslash\,\Bb_m\to\Jb\,\ts\backslash\,\Bb_m
\end{equation*}
as the composition $\bar\xi_{\ts a}\,\sih$
applied to the elements of $\Bb_m$
taken modulo $\Jb\ts$. Here the simple reflection
$\si\in\H_m$ is chosen as in Proposition \ref{p2}.
In their present form, the operators $\xic_{\ts 1}\lcd\xic_{\ts m}$  
on the vector space 
$\Jb\,\ts\backslash\,\Bb_m$ have been defined in \cite{KO}.
We call them the \textit{Zhelobenko operators}.
The next proposition states their key property. 

\begin{proposition*}
\label{p3}
The Zhelobenko operators
satisfy the braid relations corresponding to the Lie algebra\/ $\f_m\ts$.
Namely, in the case\/ $\f_m=\sp_{2m}$ we have
\begin{align}
\label{braid1}
\xic_{\ts a}\,\xic_{\ts a+1}\,\xic_{\ts a}
&\,=\,
\xic_{\ts a+1}\,\xic_{\ts a}\,\xic_{\ts a+1}
\!\!\!\quad\quad\textit{for}\ \quad
a=1\lcd m-2\ts;
\\
\label{braid2}
\xic_{\ts a}\,\xic_{\ts b}
&\,=\,
\xic_{\ts b}\,\xic_{\ts a}
\hspace{38pt}
\quad\textit{for}\ \quad
a=1\lcd b-2\ts;
\\[2pt]
\nonumber
\xic_{\ts m-1}\,\xic_{\ts m}\,\xic_{\ts m-1}\,\xic_{\ts m}
&\,=\,
\xic_{\ts m}\,\xic_{\ts m-1}\,\xic_{\ts m}\,\xic_{\ts m-1}\ts.
\end{align}
In the case when\/ $\f_m=\so_{2m}$ and\/ $m>1$, we have the same
relations \eqref{braid1} and\/ \eqref{braid2}
between\/ $\xic_{\ts1}\lcd\xic_{\ts m-1}$ as in the case\/ $\f_m=\sp_{2m}$
above, and also the relations 
\begin{align}
\label{xiam}
\xic_a\,\xic_{\ts m}
&\,=\,
\xic_{\ts m}\,\xic_{\ts a}
\ \quad\textit{for}\ \quad
a=1\lcd m-3\com m-1\ts;
\\[2pt]
\nonumber
\xic_{\ts m-2}\,\xic_{\ts m}\,\xic_{\ts m-2}
&\,=\,
\xic_{\ts m}\,\xic_{\ts m-2}\,\xic_{\ts m}\ts.
\end{align}
\end{proposition*}

\begin{proof}
Let $\e$ be 
the Lie algebra spanned 
by a central element $C\ts$, and by
the elements $x_{ai}\ts,\d_{ai}$ 
with $a=1\lcd m$ and $i=1\lcd n$ so that
Lie brackets in $\e$ are
$$
[\ts x_{ai}\com\ts x_{bj}\ts ]=0\ts,
\ \quad
[\ts \d_{ai}\com\ts\d_{\ts bj}\ts ]=0\ts,
\ \quad
[\ts x_{ai}\com\ts\d_{\ts bj}\ts]=\de_{ab}\,\de_{ij}\,C\ts.
$$
Thus $\e$ is the 
\textit{Heisenberg algebra\/} of dimension $2\ts m\ts n+1\ts$.
Let $\mathfrak{a}$ be the semidirect sum of the Lie algebras
$\f_m$ and $\e$ with the cross brackets
$\eqref{defar}$ for any $X\in\f_m$ and $Y\in\e\ts$,
where $[\,\zeta_n(X)\com Y\ts]\in\e$ is a commutator 
taken in $\PD\ts(\CC^{\ts m}\ot\CC^{\ts n})\ts$.
Then $\B_m$ is isomorphic to the quotient
of the universal enveloping algebra $\U(\mathfrak{a})$
relative to the twosided ideal generated by the element $C+1\ts$.   
Proposition \ref{p3} now follows from the
\textit{cocycle properties} of the maps
$\bar\xi_{\ts 1}\lcd\bar\xi_{\ts m}$ established in 
\cite[Section~6]{Z}\ts; 
see also \cite[Sections 4 and 6]{KO}.
\end{proof}

For $\f_m=\sp_{2m}\ts$, by using 
any reduced decomposition of an element $\si\in\H_m$ in terms
of the involutions $\si_1\lcd\si_m\ts$, we can
define a linear operator 
\begin{equation}
\label{xis}
\xic_{\ts\si}:\,\Jb\,\ts\backslash\,\Bb_m\to\Jb\,\ts\backslash\,\Bb_m
\end{equation}
in the usual way, like in \eqref{sih}.
This definition of $\xic_{\ts\si}$
is independent of the choice of a reduced decomposition of $\si$
due to Proposition \ref{p3}. 

When $\f_m=\sp_{2m}\ts$, the total number
of the factors $\si_1\lcd\si_m$ in any reduced decomposition
$\si\in\H_m$ will be denoted $\ell\ts(\si)\ts$. This number
is also independent of the choice 
of the decomposition, and is equal to
the number of elements in the set
\begin{equation}
\label{posneg}
\De_{\ts\si}=\{\ts\eta\in\De^+\,|\,\si\ts(\eta)\notin\De^+\,\}
\end{equation}
where $\De^+$ denotes the set of positive roots of
the Lie algebra $\sp_{2m}\ts$.

Now suppose that $\f_m=\so_{2m}\ts$.
Then by using any reduced decomposition in terms of 
$\si_1\lcd\si_{m-1}\com\sip_m\ts$, 
we can define a linear operator \eqref{xis} for every element
$\si\in\Hp_m\ts$. Again, this definition
is independent of the choice of a reduced decomposition of $\si$
due to Proposition \ref{p3}.
It turns out that in this case
we can extend the definition of the operator \eqref{xis}
to any element $\si\in\H_m\ts$, where $m\geqslant1$. Note that 
in this case the action of the element 
$\sih_m$ on $\Bb_m$ preserves the ideal $\Jb\ts$,
and therefore induces a linear operator
on the quotient vector space $\Jb\,\ts\backslash\,\Bb_m\ts$.
This operator will be again denoted by $\sih_m\ts$. 
The extension of the definition of the operators \eqref{xis} to
$\si\in\H_m$ is based on the next lemma.

\begin{lemma*}
\label{lemma45}
When $\f_m=\so_{2m}$ and\/ $m>1\ts$, the operators\/ 
$\xic_{\ts1}\lcd\xic_{\ts m-1}\com\sih_m$ on $\Jb\,\ts\backslash\,\Bb_m$
satisfy the same relations, as the\/ $m$
generators of the braid group $\Hh_m$ respectively.
Then we also have the relation
\begin{equation}
\label{ximh}
\xic_{\ts m}=
\sih_m\,\xic_{\ts m-1}\,\sih_m\ts.
\end{equation}
\end{lemma*}

\begin{proof}
On the top of the relations \eqref{braid1},\eqref{braid2}
between 
$\xic_{\ts1}\lcd\xic_{\ts m-1}$ which we 
have by Proposition \ref{p3}, we need the following
relations involving the operator $\sih_m$ on $\Jb\,\ts\backslash\,\Bb_m\ts$:
\begin{align*}
\xic_{\ts a}\,\sih_m
&\,=\,
\xic_{\ts a}\,\sih_m
\!\!\!\!\quad\quad\textit{for}\ \quad
a=1\lcd m-2\ts;
\\
\xic_{\ts m-1}\,\sih_m\,\xic_{\ts m-1}\,\sih_m
&\,=\,
\sih_m\,\xic_{\ts m-1}\,\sih_m\,\xic_{\ts m-1}\ts.
\end{align*}
The first equality here immediately follows from the definition
of the operator $\xic_{\ts a}\ts$. The equality in the last 
displayed line follows by choosing $a=m-1$ in \eqref{xiam}, and 
then using \eqref{ximh}. Further, from \eqref{xim} we derive
an equality of operators on $\Bb_m\ts$,
$$
\bar\xi_{\ts m}\,\sih_m\,\sih_{m-1}\,\sih_m
=
\sih_m\,\bar\xi_{\ts m-1}\,\sih_{m-1}\,\sih_m\ts,
$$
which implies the equality \eqref{ximh} 
of operators on the quotient space $\Jb\,\ts\backslash\,\Bb_m\ts$.
\end{proof}

Now for $\f_m=\so_{2m}$ with any $m\geqslant1$,
take any decomposition of an element $\si\in\H_m$ in terms
of the involutions $\si_1\lcd\si_m$ such that the
number of occurencies of $\si_1\lcd\si_{m-1}$ in the decomposion
is minimal possible. For $\f_m=\so_{2m}$ the symbol
$\ell\ts(\si)$ will denote this minimal number. 
Note that unlike for $\f_m=\sp_{2m}\ts$,
here we do not count the occurencies of $\si_m$ in the
decomposition. All the decompositions of $\si\in\H_m$ with
the minimal number of occurencies of $\si_1\lcd\si_{m-1}$
can be obtained from each other by using the braid
relations between $\si_1\lcd\si_m\in\H_m$ along with
the relation $\si_m^{\ts2}=1\ts$. 

By substituting the operators
$\xic_{\ts1}\lcd\xic_{\ts m-1}\com\sih_m$ on 
$\Jb\,\ts\backslash\,\Bb_m$ for
involutions $\si_1\lcd\si_m$ in such a decomposition of 
$\si\in\H_m\ts$, we obtain another operator on 
$\Jb\,\ts\backslash\,\Bb_m\ts$. The latter operator
does not depend on the choice of a decomposition
because of the first statement of Lemma \ref{lemma45},
and because the operator $\sih_m^{\,2}$ on 
the vector space $\Jb\,\ts\backslash\,\Bb_m$ is the identity
in the case $\f_m=\so_{2m}$ considered here.
Moreover for $\si\in\Hp_m\subset\H_m\ts$,
the operator on $\Jb\,\ts\backslash\,\Bb_m$
obtained by the latter substitution coincides with the operator \eqref{xis}.
Indeed, for $\f_m=\so_{2m}$
the operator \eqref{xis} has been defined by
substituting the Zhelobenko operators 
$\xic_{\ts1}\lcd\xic_{\ts m-1}\com\xic_{\ts m}$ 
for $\si_1\lcd\si_{m-1}\com\sip_m$
in any reduced decomposition of $\si\in\Hp_m\ts$.
The coincidence of the two operators for $\si\in\Hp_m$
now follows from the relation \eqref{ximh}.
Thus we have extended the definition of 
the operator \eqref{xis} from $\si\in\Hp_m$ to all $\si\in\H_m\ts$.

Note that for $\f_m=\so_{2m}$ and
$\si\in\Hp_m\ts$, the number $\ell\ts(\si)$ is equal to
the length of a reduced decomposition of $\si$ in terms of 
$\si_1\lcd\si_{m-1}\com\sip_m\ts$. Thus
we have also extended the standard length function
from the Weyl group $\Hp_m$ of $\so_{2m}$ to
the hyperoctahedral group $\H_m\ts$.
Moreover for any $\si\in\H_m\ts$,
not only for $\si\in\Hp_m\ts$,
the number $\ell\ts(\si)$ equals
the number of elements in the set \eqref{posneg}, where
$\De^+$ is the set of positive roots of $\so_{2m}\ts$.

From now we shall consider
$\f_m=\sp_{2m}$ and $\f_m=\so_{2m}$ simultaneously,
and will work with the operators \eqref{xis} for all elements
$\si\in\H_m\ts$.
In particular, in the case $\f_m=\so_{2m}$
we will assume that the operator \eqref{xis} with $\si=\si_m$
acts as~$\sih_m\ts$.

The restriction of
the action \eqref{siact}\com\ts\eqref{simact} of the braid group $\Hh_m$
on $\f_m$ to the Cartan subalgebra $\h$ factors 
to an action of the hyperoctahedral group $\H_m\ts$. This is the 
standard action of the Weyl group of $\f_m=\sp_{2m}\ts$.
The resulting action of the subgroup $\Hp_m\subset\H_m$ on $\h$
is the standard action of the Weyl group of $\f_m=\so_{2m}\ts$.
The group $\H_m$ also acts on the dual vector space $\h^\ast\ts$, so that
$\si\ts(\ep_c)=\ep_{\si(c)}$ for any $\si\in\H_m$ and any
$c=-\ts m\lcd-1\com1\lcd m\ts$. Unlike in \eqref{sib}, here we use
the natural action of the group $\H_m$ by permutations
of $-\ts m\lcd-1\com1\lcd m\ts$. Thus $\si_a\in\H_m$ 
with $1\le a<m$ exchanges $a\com a+1$ and also exchanges
$-a\com-a-1\ts$ while $\si_m\in\H_m$ exchanges $m\com-m\ts$.
Note that we always have $\si\ts(-c)=-\si\ts(c)\ts$.
If we identify each weight $\mu\in\h^\ast$ with
the sequence $(\ts\mu_1\lcd\mu_m)$ of its labels, then
\begin{align*}
\si:(\ts\mu_1\lcd\mu_m)
&\,\mapsto\,
(\,\mu_{\ts\si^{-1}(1)}\lcd\mu_{\ts\si^{-1}(m)}\ts)
\quad\text{for}\quad
\si\in\Sym_m\ts,
\\
\si_m:(\ts\mu_1\lcd\mu_m)
&\,\mapsto\,
(\,\mu_{1}\lcd\mu_{m-1},-\mu_{m}\ts)\ts.
\end{align*}

The \textit{shifted\/} action of the
group $\H_m$ on the set $\h^\ast$ is defined by the assignment
\begin{equation*}
\mu\,\mapsto\,\si\circ\mu=\si\ts(\mu+\rho)-\rho
\ \quad\textrm{for}\quad\ 
\si\in\H_m\ts.
\end{equation*}
By regarding the elements of the commutative algebra $\Uhb$
as rational functions on the vector space $\h^\ast$
we can also define an action of the group $\H_m$ on this \textrm{algebra:}
\begin{equation}
\label{saction2}
(\si\circ X)(\mu)=X(\si^{-1}\ns\circ\mu)
\,\quad\textrm{for}\quad
X\in\Uhb\ts.
\end{equation}

\begin{proposition*}
\label{saction}
For any $\si\in\H_m\ts$, $X\in\Uhb$ and 
$Y\in\Jb\,\ts\backslash\,\Bb_m$ we have the relations
\begin{align}
\label{q121}
\xic_{\ts\si}(X\ts Y)
&\,=\,
(\ts\si\ts\circ X)\,\ts\xic_{\ts\si}(\ts Y)\ts,
\\
\label{q122}
\xic_{\ts\si}(\ts YX)
&\,=\,
\,\xic_{\ts\si}(\ts Y)\ts(\ts\si\ts\circ X)\ts.
\end{align}
\end{proposition*}

\begin{proof}
It suffices to check the relations \eqref{q121} and \eqref{q122}
only for $\si=\si_1\lcd\si_m$ and $X\in\h\ts$. When
$\f_m=\so_{2m}$ and $\si=\si_m\ts$, the operator $\xic_{\ts\si}$
acts on $\Jb\,\ts\backslash\,\Bb_m$ as $\sih_m$ by our definition,
while $\si_m\ts\circ X=\si_m\ts(X)$ because $\rho_m=0\ts$. So
the relations \eqref{q121} and \eqref{q122} are obvious then.
In the remaining cases the equalities \eqref{q121} and \eqref{q122}
follow respectively from \eqref{q11} and \eqref{q12}.
Indeed, in the remaining cases for every possible index $a$
we have $\si_a(\ts\mu+\eta_{\ts a})=\si_a\circ\mu\,$;
see also 
\eqref{N14},\eqref{N1414} and \eqref{saction2}.
\end{proof}


\section{Intertwining operators}
\setcounter{section}{5}
\setcounter{equation}{0}
\setcounter{theorem*}{0}

Let $\de=(\ts\de_1\lcd\de_m)$ be any sequence of $m$ 
elements from the set $\{1\com-1\}\ts$. The hyperoctahedral
group $\H_m$ acts on the set of all these sequences naturally,
so that the generator $\si_a\in\H_m$ with $a<m$ acts on $\de$ as 
the transposition
of $\de_a$ and $\de_{a+1}\ts$, while the generator $\si_m\in\H_m$
changes the sign of $\de_m\ts$. Let
$\de_+=(1\lcd1)$ be
the sequence of $m$ elements $1$.
Given any sequence $\de\ts$, take the composition of the automorphisms
of the ring $\PD\ts(\CC^{\ts m}\ot\CC^{\ts n})\ts$,
\begin{equation}
\label{compfour}
x_{\ts\ab i}\mapsto-\,\th_i\,\d_{\,\ab\ts\bi}
\quad\text{and}\quad
\d_{\,\ab i}\mapsto\th_i\,x_{\ts\ab\ts\bi}
\quad\text{whenever}\quad
\de_a=-1\ts.
\end{equation}
Here $a\ge1$ and $i=1\lcd n\ts$.
Let us denote by $\varpi$ this composition.
In particular, the automorphism $\varpi$
corresponding to $\de=(1\lcd1\com-1)$ coincides
with the action of $\sih_m$ on $\PD\ts(\CC^{\ts m}\ot\CC^{\ts n})\ts$,
see \eqref{fourier}. In the case $\f_m=\so_{2m}\ts$,
the automorphism $\varpi$ is involutive for any
$\de\ts$. But in the case $\f_m=\sp_{2m}$ the square $\varpi^2$ maps
$$
x_{\ts\ab i}\mapsto-\ts x_{\,\ab i}
\quad\text{and}\quad
\d_{\,\ab i}\mapsto-\ts\d_{\,\ab i}
\quad\text{whenever}\quad
\de_a=-1\ts.
$$

For any $\f_m\ts$-module $V$, the action of
the extended twisted Yangian $\X(\g_n)$ on its module
$\F_m(V)=V\ot\P\ts(\CC^{\ts m}\ot\CC^{\ts n})$
is defined by means of the homomorphism 
$\be_m:\X(\g_n)\to\U(\ts\f_m)\ot\PD\ts(\CC^{\ts m}\ot\CC^{\ts n})\ts$,
see Proposition \ref{xb}. Further,
the action of Lie algebra $\f_m$ on the second tensor factor
$\P\ts(\CC^{\ts m}\ot\CC^{\ts n})$ of $\F_m(V)$
is defined by means of homomorphism
$\zeta_n:\U(\ts\f_m)\to\PD\ts(\CC^{\ts m}\ot\CC^{\ts n})\ts$, see 
definition \eqref{gan}.
Here any element of the ring 
$\PD\ts(\CC^{\ts m}\ot\CC^{\ts n})$ acts on the vector space
$\P\ts(\CC^{\ts m}\ot\CC^{\ts n})$ naturally. 
We can modify the latter action,
by making any element $Y\in\PD\ts(\CC^{\ts m}\ot\CC^{\ts n})$
act on $\P\ts(\CC^{\ts m}\ot\CC^{\ts n})$ via
the natural action of $\varpi\ts(Y)\ts$.
Then we get another $\PD\ts(\CC^{\ts m}\ot\CC^{\ts n})\ts$-module,
with the same underlying vector space 
$\P\ts(\CC^{\ts m}\ot\CC^{\ts n})$ for every~$\de\ts$.

For any $\f_m\ts$-module $V$, we can now define a bimodule
$\F_\de\ts(V)$ of $\f_m$ and $\X(\g_n)\ts$.
Its underlying vector space is the same
$V\ot\P\ts(\CC^{\ts m}\ot\CC^{\ts n})$ for every~$\de\ts$.
The action of $\X(\g_n)$ on $\F_\de\ts(V)$
is defined by pushing the homomorphism $\be_m$
forward through the automorphism $\varpi$, 
applied to $\PD\ts(\CC^{\ts m}\ot\CC^{\ts n})$
as to the second tensor factor of the target of $\be_m\ts$.
The action of $\f_m$ on $\F_\de\ts(V)$
is also defined by pushing the homomorphism $\zeta_n$
forward through the automorphism $\varpi\ts$. 
Note that the action of $\f_m$ on the tensor factor $V$ 
of $\F_\de\ts(V)$ is not modified.
We have $\F_m(V)=\F_{\ts\de_+}(V)\ts$.

Now let $\mu\in\h^\ast$ be any weight of $\f_m\ts$, such that
\begin{equation}
\label{notinso}
\mu_a-\mu_b\notin\ZZ
\quad\text{and}\quad
\mu_a+\mu_b\notin\ZZ
\quad\text{whenever}\quad
1\le a<b\le m\ts.
\end{equation}
In the case $\f_m=\sp_{2m}$ also suppose that,
in addition to \eqref{notinso},
\begin{equation}
\label{notinsp}
2\ts\mu_a\notin\ZZ
\quad\text{whenever}\quad
1\le a\le m\ts.
\end{equation}
We shall now proceed to show how for every element $\si\in\H_m\ts$,
the Zhelobenko operator \eqref{xis} 
determines an $\X(\g_n)$-intertwining operator
\begin{equation}
\label{distoper}
\F_{m}\ts(\ts M_{\ts\mu})_{\ts\n}
\,\to\,
\F_{\ts\de}\ts(\ts M_{\ts\si\,\circ\ts\mu})_{\ts\n}
\ \quad\text{where}\quad\ 
\de=\si\ts(\de_+)\ts.
\end{equation}

In this section we keep regarding $\B_m$ 
as the associative algebra generated by
$\U(\ts\f_m)$ and $\PD\ts(\CC^{\ts m}\ot\CC^{\ts n})$ with the cross relations
\eqref{defar}.
Let $\I_{\ts\de}$ be the left ideal of algebra $\B_m$
generated by the elements $x_{\ts\ab\ts i}$ with $\de_a=-1\ts$,
and the elements $\d_{\ts\ab\ts i}$ with $\de_a=1\ts$.
Here $a=1\lcd m$ and $i=1\lcd m\ts$. 
Note that in terms of the elements $q_{\ts ci}$ 
introduced immediately after stating Proposition~\ref{xb},
the left ideal $\I_{\ts\de}$ is generated by the elements 
$q_{\ts-\de_a\ab,i}$ where again $a=1\lcd m$ and $i=1\lcd m\ts$.
In particular, the ideal $\I_{\ts\de_+}$ is generated by
all the partial derivations $\d_{\ts ai}\ts$.
Let $\Ib_{\ts\de}$ be the left ideal of 
$\Bb_m$
generated by the same elements as the ideal of $\I_{\ts\de}$ of $\B_m\ts$.

Consider the image of the ideal $\Ib_{\ts\de}$
in the quotient space $\Jb\,\ts\backslash\,\Bb_m\ts$,
that is the subspace $\Jb\,\ts\backslash\,(\ts\Ib_{\ts\de}+\Jb\ts)$
in the quotient space. The image will be occasionally 
denoted by the same symbol
$\Ib_{\ts\de}\ts$. In the context of 
next proposition, this notation should cause no confusion.

\begin{proposition*}
\label{propznak}
For any $\si\in\H_m$ the operator $\xic_{\ts\si}$ maps the subspace\/ 
$\I_{\ts\de_+}$ to\/ 
$\ts\I_{\ts\si(\de_+)}\ts$.
\end{proposition*}

\begin{proof}
For any $a=1\lcd m-1$ consider the operator $\widehat{F}_a$
corresponding to the element $F_a\in\B_m\ts$. By 
\eqref{Fc} and \eqref{gan},\eqref{defar} for any
$Y\in\PD\ts(\CC^{\ts m}\ot\CC^{\ts n})$ we have
$$
\widehat{F}_a(\ts Y)=\,
-\,\sum_{k=1}^n\,\,
[\ts x_{\ts\ab\ts k}\,\d_{\,\ts\overline{a+1}\ts k}\,,Y\ts]\ts.
$$
Similarly, in the case $\f_m=\sp_{2m}$ by 
\eqref{Fm} for any $Y\in\PD\ts(\CC^{\ts m}\ot\CC^{\ts n})$ we have
$$
\widehat{F}_m(\ts Y)=\,
-\,\sum_{k=1}^n\,\,
[\ts x_{\ts\overline{m}\ts\bk}\,x_{\ts\overline{m}\ts k}\,,Y\ts]\ts/\ts2\ts.
$$
In the case $\f_m=\so_{2m}$ we do not need to consider
the operator $\widehat{F}_m\ts$, because in this case
the operator \eqref{xis} corresponding to $\si=\si_m$ 
acts on $\Jb\,\ts\backslash\,\Bb_m$ as $\sih_m$ by our definition.

The above description of the action of $\widehat{F}_a$ with $a<m$ on 
the vector space $\PD\ts(\CC^{\ts m}\ot\CC^{\ts n})$ shows that this action 
preserves each of the two $2n$ dimensional subspaces, spanned by the vectors
\begin{gather}
\label{xxideal}
x_{\ts\ab\ts i}
\quad\text{and}\quad
x_{\ts\overline{a+1}\ts i}
\quad\text{where}\quad
i=1\lcd n\ts;
\\[2pt]
\label{ddideal}
\d_{\,\ab\ts i}
\quad\text{and}\quad
\d_{\,\overline{a+1}\ts i}
\quad\text{where}\quad
i=1\lcd n\ts.
\end{gather}
This action also maps to zero the $2n$ dimensional subspace, spanned by
\begin{equation}
\label{xdideal}
x_{\ts\ab\ts i}
\quad\text{and}\quad
\d_{\,\overline{a+1}\ts i}
\quad\text{where}\quad
i=1\lcd n\ts.
\end{equation}
Therefore for any $\de\ts$, the operator $\bar\xi_a$ with $a<m$ maps
the left ideal $\Ib_{\ts\de}$ of $\Bb_m$ to the image of $\Ib_{\ts\de}$ in
$\Jb\,\ts\backslash\,\Bb_m\ts$, unless $\de_a=1$ and $\de_{a+1}=-1\ts$.
The operator $\xic_{\ts a}$ on 
$\Jb\,\ts\backslash\,\Bb_m$ was defined by
taking the composition of 
$\bar\xi_a$ and
$\sih_a\ts$. Hence $\xic_{\ts a}$ with $a<m$ maps the image of 
$\Ib_{\ts\de}$ to the image of $\Ib_{\ts\si_a(\de)}\ts$, unless
$\de_a=-1$ and $\de_{a+1}=1\ts$.

In the case $\f_m=\sp_{2m}\ts$,
the action of $\widehat{F}_m$ 
on the vector space $\PD\ts(\CC^{\ts m}\ot\CC^{\ts n})$
maps to zero the $n$ dimensional subspace spanned by the elements
\begin{equation}
\label{xideal}
x_{\ts\overline{m}\ts i}\ts=\ts x_{1i}
\quad\text{where}\quad
i=1\lcd n\ts.
\end{equation}
Therefore the operator $\bar\xi_m$ maps
the left ideal $\Ib_{\ts\de}$ of $\Bb_m$ to the image of $\Ib_{\ts\de}$ in
$\Jb\,\ts\backslash\,\Bb_m\ts$, unless $\de_m=1$.
Hence the operator $\xic_{\ts m}$ on 
$\Jb\,\ts\backslash\,\Bb_m$ maps the image of 
$\Ib_{\ts\de}$ 
to the image of $\Ib_{\ts\si_m(\de)}\ts$,
unless $\de_m=-1\ts$. In the case $\f_m=\so_{2m}\ts$,
we just note that 
$\sih_m$ maps the image of
$\Ib_{\ts\de}$ in $\Jb\,\ts\backslash\,\Bb_m$ to the image of
$\Ib_{\ts\si_m(\de)}\ts$.

From now on we will denote
the image of the ideal $\Ib_{\ts\de}$ 
in the quotient space $\Jb\,\ts\backslash\,\Bb_m\ts$
by the same symbol. Put
$$
\widehat{\de}=\sum_{a=1}^m\,\de_a\ts\ep_a\ts.
$$
Then for each $\si\in\H_m$ we have the equality 
$\widehat{\si\ts(\de\ts)}=\si\ts(\,\widehat{\de}\,)$
where at the right hand side we use the action of
the group $\H_m$ on $\h^\ast$. Let $(\ ,\,)$ be the standard bilinear
form on $\h^\ast$, so that the basis of weights $\ep_a$ with $a=1\lcd m$
is orthonormal. The above remarks on the action
of the Zhelobenko operators on $\Ib_{\ts\de}$ can now be restated 
as follows:
\begin{align}
\label{rest1}
\textrm{if}\quad(\ts\widehat{\de}\com\ep_a-\ep_{a+1})&\ge0
&
\textrm{then}&
&
&\ts\xic_{\ts a}\ts(\,\Ib_\de\ts)\subset\ts\Ib_{\ts\si_a(\de)}
\,\quad\textrm{for}\quad
a=1\lcd m-1\,;
\\
\label{rest2}
\textrm{if}\quad(\ts\widehat{\de}\com\ep_m)&>0
&
\textrm{then}&
&
&\ns\xic_{\ts m}(\,\Ib_\de\ts)\subset\Ib_{\ts\si_m(\de)}
\quad\textrm{for}\quad
\f_m=\sp_{2m}\ts.
\end{align}

We shall prove Proposition \ref{propznak}
by the induction on the length of a reduced decomposition
of $\si\in\H_m$ in terms of $\si_1\lcd\si_m\ts$.
This number has been denoted by $\ell(\si)$ in the case
$\f_m=\sp_{2m}\ts$, but may be different from the number
denoted by $\ell(\si)$ in the case $\f_m=\so_{2m}\ts$.
Recall that in both cases $\ell\ts(\si)$ equals
the number of elements in the set \eqref{posneg},
where $\De^+$ is the set of positive roots of $\f_m\ts$.

If $\si$ is the identity element of $\H_m\ts$,
Proposition \ref{propznak} is tautological.
Suppose that for some $\si\in\H_m\ts$,
$$
\xic_{\ts\si}(\,\Ib_{\ts\de_+}\ts)\subset\ts\Ib_{\ts\si(\de_+)}\ts.
$$
Take $\si_a\in\H_m$ with $1\le a\le m\ts$,
such that $\si_a\ts\si$ has a longer reduced decomposition
in terms of $\si_1\lcd\si_m$ than $\si\ts$.
If $\f_m=\so_{2m}$ and $a=m\ts$, then
$\xic_{\ts\si_m\ts\si}=\sih_m\,\xic_{\ts\si}$ and we need the inclusion
\begin{equation}
\label{obvinc}
\sih_m\ts(\,\Ib_{\ts\si(\de_+)}\ts)\subset\ts\Ib_{\ts\si_m\si(\de_+)}\ts,
\end{equation}
which holds by the definition of the action of $\H_m$
on $\Jb\,\ts\backslash\,\Bb_m\ts$.

We may exclude the case when $\f_m=\so_{2m}$ and $a=m\ts$, and assume that
\begin{equation}
\label{lengthup}
\ell\ts(\ts\si_a\ts\sigma)=\ell\ts(\sigma)+1\ts.
\end{equation}
Firstly, suppose that $a<m$ here. Let us then prove the inclusion 
$$
\xic_{\ts a}(\,\Ib_{\ts\si(\de_+)}\ts)\subset\ts\Ib_{\ts\si_a\si(\de_+)}\ts.
$$
By \eqref{rest1}, the latter inclusion will have place if 
\begin{equation*}
(\ts\widehat{\si\ts(\de_+)}\ts\com\ts\ep_a-\ep_{a+1})=
(\ts\si\ts(\ts\widehat{\de}_+)\ts\com\ts\ep_a-\ep_{a+1})\ge0\ts.
\end{equation*}
But the condition \eqref{lengthup} for $a<m$ implies that
$\ep_a-\ep_{a+1}\in\sigma\ts(\De^+)\ts$. Indeed, because the root
$\ep_a-\ep_{a+1}$ of $\f_m$ is simple,
$\si_a(\eta)\in\De^+$ for any
$\eta\in\De^+$ such that $\eta\neq\ep_a-\ep_{a+1}\ts$.
Since $\ell\ts(\sigma)$ and $\ell\ts(\ts\si_a\ts\sigma)$ 
are the numbers of elements in
$\De_{\ts\si}$ and $\De_{\ts\si_a\ts\si}$ respectively,
here $\ep_a-\ep_{a+1}\in\sigma\ts(\De^+)\ts$.
So $\ep_a-\ep_{a+1}=\si\ts(\ts\ep_b-\ep_c)$ where
$1\le b\le m$ and $1\le|c|\le m\ts$. Thus
$$
(\ts\si\ts(\ts\widehat{\de}_+)\ts\com\ts\ep_a-\ep_{a+1})
=
(\ts\si\ts(\ts\widehat{\de}_+)\ts\com\ts\si\ts(\ts\ep_b-\ep_c)\ts)
=
(\ts\widehat{\de}_+\ts\com\ts\ep_b-\ep_c\ts)\ge0\ts.
$$

Now suppose that $a=m\ts$. Here we assume that
$\f_m=\sp_{2m}\ts$. We need the inclusion
$$
\xic_{\ts m}\ts(\,\Ib_{\ts\si(\de_+)}\ts)
\subset\ts
\Ib_{\ts\si_m\si(\de_+)}\ts.
$$
It will have place if 
\begin{equation*}
(\ts\widehat{\si\ts(\de_+)}\com\ep_m\ts)
=
(\ts\si\ts(\ts\widehat{\de}_+)\ts\com\ep_m\ts)>0\ts.
\end{equation*}
But the condition \eqref{lengthup} for $a=m$ implies 
that $2\ts\ep_m\in\sigma\ts(\De^+)\ts$,
where $\De^+$ is the set of positive roots of $\sp_{2m}\ts$.
Indeed, because the root
$2\ts\ep_m$ of $\sp_{2m}$ is simple,
$\si_m(\eta)\in\De^+$ for any
$\eta\in\De^+$ such that $\eta\neq2\ts\ep_m\ts$.
Since $\ell\ts(\sigma)$ and $\ell\ts(\ts\si_m\ts\sigma)$ 
are the numbers of elements in
$\De_\si$ and $\De_{\ts\si_m\ts\si}$ respectively,
here $2\ts\ep_m\in\sigma\ts(\De^+)\ts$.
So $\ep_m=\si\ts(\ts\ep_b)$ where $1\le b\le m\ts$. Thus
$$
(\ts\si\ts(\ts\widehat{\de}_+)\ts\com\ep_m\ts)
=
(\ts\si\ts(\ts\widehat{\de}_+)\ts\com\si\ts(\ts\ep_b)\ts)
=
(\ts\widehat{\de}_+\ts\com\ts\ep_b\ts)>0\ts.
\eqno\qedhere
$$
\end{proof}

\begin{corollary*}
\label{corollary4.6}
For any $\si\in\H_m$ the operator $\xic_{\ts\si}$ on 
$\Jb\,\ts\backslash\,\Bb_m$ maps 
$$
\Jb\,\ts\backslash\,(\Jpb+\Ib_{\ts\de_+}+\Jb\ts)
\,\to\ts
\Jb\,\ts\backslash\,(\Jpb+\Ib_{\ts\si(\de_+)}+\Jb\ts)\ts.
$$
\end{corollary*}

\begin{proof}
We will extend the arguments used in the proof of
Proposition \ref{propznak} above. In particular, we will again use
the length of a reduced decomposition
of $\si$ in terms of $\si_1\lcd\si_m\ts$.
If $\si$ is the identity element of $\H_m\ts$,
then the required statement is tautological.
Now suppose that for some $\si\in\H_m$
the statement of Corollary \ref{corollary4.6} is true.
Take any simple reflection $\si_a\in\H_m$ with $1\le a\le m\ts$,
such that $\si_a\ts\si$ has a longer reduced decomposition
in terms of $\si_1\lcd\si_m$ than $\si\ts$.
In the case $\f_m=\so_{2m}$ we may assume
that $a<m\ts$, because in that case
the required statement for $\si_m\ts\si$ instead of $\si$ 
is provided by the inclusion \eqref{obvinc}.

Thus we have 
\eqref{lengthup}. 
With the above assumption on $a\ts$, we have proved that 
\eqref{lengthup} implies 
\begin{equation}
\label{thischeck}
(\ts\widehat{\si\ts(\de_+)}\ts\com\ts\eta_{\ts a})\ge0\ts.
\end{equation}
Here $\eta_a$ is the simple root corresponding to 
$\si_a\ts$. But \eqref{thischeck} implies the equality
\begin{equation}
\label{N10}
\Jpb+\Ib_{\ts\si(\de_+)}=
\Jpb_{\ts a}+\Ib_{\ts\si(\de_+)}
\end{equation}
of left ideals of $\Bb_m\ts$.
Indeed, the left and right hand sides of \eqref{N10} differ by
the elements $Y\ts\zeta_n(E_a)$ where $Y$ ranges over $\Bb_m\ts$.
The condition \eqref{thischeck} implies that
$\zeta_n(E_a)\in\Ib_{\ts\si(\de_+)}\ts$, see
the definition \eqref{gan} and the arguments in
the beginning of proof of Proposition~\ref{propznak}. 
Using Proposition \ref{prop3N} and the induction step 
from our proof of Proposition~\ref{propznak}, 
$\xic_{\ts a}$ maps
$$
\Jb\,\ts\backslash\,(\Jpb+\Ib_{\ts\si(\de_+)}+\Jb\ts)
=
\Jb\,\ts\backslash\,(\Jpb_{\ts a}+\Ib_{\ts\si(\de_+)}+\Jb\ts)
\,\to\ts
\Jb\,\ts\backslash\,(\Jpb+\Ib_{\ts\si_a\si(\de_+)}+\Jb\ts)\ts.
$$
This makes the induction step of our proof of Corollary \ref{corollary4.6}.
\end{proof}

Let $\I_{\ts\mu,\de}$ be the left ideal of the algebra $\B_m$
generated by $\I_{\ts\de}+\Jp$ and by the elements
$$
F_{-\ab,-\ab}-\zeta_n\ts(F_{-\ab,-\ab})-\mu_a
\quad\text{where}\quad
a=1\lcd m\ts.
$$
Recall that under the isomorphism of the algebra $\B_m$ with 
$\U(\ts\f_m)\ot\PD\ts(\CC^{\ts m}\ot\CC^{\ts n})\ts$,
the difference $X-\zeta_n(X)\in\B_m$ for every $X\in\f_m$
is mapped to the element \eqref{differ}.
Denote by $\Ib_{\ts\mu,\de}$ the subspace
$\Uhb\,\I_{\ts\mu,\de}\ts$ of $\Bb_m\ts$, this
is also a left ideal of 
$\Bb_m\ts$.

\begin{theorem*}
\label{proposition4.5}
For any element $\si\in\H_m$ the operator $\xic_{\ts\si}$ on 
$\Jb\,\ts\backslash\,\Bb_m$ maps 
$$
\Jb\,\ts\backslash\,
(\ts\Ib_{\ts\mu,\de_+}+\Jb\ts)
\,\to\ts
\Jb\,\ts\backslash\,
(\ts\Ib_{\ts\si\ts\circ\ts\mu\ts,\ts\si(\de_+)}+\Jb\ts)\ts.
$$
\end{theorem*}

\begin{proof}
Let $\ka$ be a weight of $\f_m$ with the sequence of
labels $(\ts\ka_1\lcd\ka_m)\ts$. Suppose that the weight $\ka$
satisfies the conditions
\eqref{notinso} instead of $\mu\ts$. In the case $\f_m=\sp_{2m}$
we also suppose that $\ka$ satisfies the conditions
\eqref{notinsp} instead of $\mu\ts$. Denote by
$\tilde\I_{\ts\ka,\de}$ be the left ideal of $\Bb_m$
generated by $\I_{\ts\de}+\Jp$ and by the elements
$$
F_{-\ab,-\ab}-\ka_a
\quad\text{where}\quad
a=1\lcd m\ts.
$$
Proposition \ref{saction} and Corollary \ref{corollary4.6}
imply that the operator $\xic_{\ts\si}$ on 
$\Jb\,\ts\backslash\,\Bb_m$ maps
$$
\Jb\,\ts\backslash\,
(\ts\tilde\I_{\ts\ka,\de_+}+\Jb\ts)
\,\to\ts
\Jb\,\ts\backslash\,
(\ts\tilde\I_{\ts\si\ts\circ\ts\ka\ts,\ts\si(\de_+)}+\Jb\ts)\ts.
$$
Now choose 
\begin{equation}
\label{kappa}
\ka_a=\mu_a-\frac{n}2
\quad\text{for}\quad
a=1\lcd m\ts.
\end{equation}
Then the conditions on $\ka$ stated in the beginning of this proof
are satisfied. For every $\si\in\H_m$ we shall prove the equality 
of left ideals of $\Bb_m\ts$,
\begin{equation}
\label{lastin}
\tilde\I_{\ts\si\ts\circ\ts\ka\ts,\ts\si(\de_+)}=
\Ib_{\ts\si\ts\circ\ts\mu\ts,\ts\si(\de_+)}\ts.
\end{equation}
Theorem \ref{proposition4.5} will then follow. 
Denote $\de=\si\ts(\de_+)\ts$. Then by our choice of $\ka$ we have
\begin{equation}
\label{sika}
\si\ts\circ\ts\ka=\si\ts\circ\ts\mu\ts-\frac{n\ts\de}2
\end{equation}
where the sequence $\de$ is regarded as a weight of $\f_m\ts$,
by identifying the weights with their sequences of labels.
Let the index $a$ run through $1\lcd m\ts$. If 
$\de_a=1$ then by the definition \eqref{gan}, we have
$$
\zeta_n\ts(F_{-\ab,-\ab})+\frac{n}2=
-\,\sum_{k=1}^n\,x_{\ts\ab\ts k}\,\d_{\,\ab\ts k}\in\I_{\ts\de}\ts.
$$
If $\de_a=-1$ then the same definition \eqref{gan} implies that
$$
\zeta_n\ts(F_{-\ab,-\ab})-\frac{n}2=
-\,\sum_{k=1}^n\,\d_{\,\ab\ts k}\,x_{\ts\ab\ts k}\in\I_{\ts\de}\ts.
$$
Hence the relation \eqref{sika} implies the equality \eqref{lastin}.
\end{proof}

Consider the quotient vector space
$\B_m\ts/\,\I_{\ts\mu,\de}$ 
for any sequence $\de\ts$. 
The algebra $\U(\ts\f_m)$ 
acts on this quotient via left multiplication,
being regarded as a subalgebra of $\B_m\ts$. The
algebra $\X(\g_n)$ also acts on this quotient via left multiplication,
using the homomorphism $\be_m:\X(\g_n)\to\B_m\ts$. Recall that
in Section 2, the target algebra $\B_m$ of the homomorphism $\be_m$
was defined as 
$\U(\ts\f_m)\ot\PD\ts(\CC^{\ts m}\ot\CC^{\ts n})\ts$. 
Here we use a different presentation of
the same algebra, by means of the cross relations \eqref{defar}.
In particular, here the image of $\be_m$
commutes with the subalgebra $\U(\ts\f_m)$ of $\B_m\ts$; see 
Proposition \ref{xb}, Part (ii).
Thus here the vector space
$\B_m\ts/\,\I_{\ts\mu,\de}$ becomes a bimodule 
over $\f_m$ and $\X(\g_n)\ts$.

Consider the bimodule $\F_\de\ts(M_\mu)$
over $\f_m$ and $\X(\g_n)\ts$, 
defined in the beginning of this section.
This bimodule is equivalent to
$\B_m\ts/\,\I_{\ts\mu,\de}\ts$. Indeed, let $Z$ run through
$\P\ts(\CC^{\ts m}\ot\CC^{\ts n})\ts$. Then a bijective linear map
$$
\F_\de\ts(M_\mu)\to\B_m\ts/\,\I_{\ts\mu,\de}
$$
intertwining the actions of $\f_m$ and $\X(\g_n)$ can be defined
by mapping the element
$$
1_\mu\ot Z\in M_\mu\ot\P\ts(\CC^{\ts m}\ot\CC^{\ts n})
$$
to the image of 
$$
\varpi^{-1}(Z)\in\PD\ts(\CC^{\ts m}\ot\CC^{\ts n})\subset\B_m
$$
in the quotient $\B_m\ts/\,\I_{\ts\mu,\de}\ts$. 
The intertwining property follows from the definitions
of $\F_\de\ts(M_\mu)$ and $\I_{\ts\mu,\de}\ts$.
The same mapping determines a bijective linear map
\begin{equation}
\label{fbi}
\F_\de\ts(M_\mu)\ts\to\ts\Bb_m\ts/\,\Ib_{\ts\mu,\de}\ts.
\end{equation}

In particular, the space $\F_\de\ts(M_\mu)_{\ts\n}$
of $\n\ts$-coinvariants of $\F_\de\ts(M_\mu)$
is equivalent to the quotient
$\Jb\,\ts\backslash\,\Bb_m\ts/\,\Ib_{\ts\mu,\de}$
as a bimodule over the Cartan subalgebra
$\h\subset\f_m$ and over $\X(\g_n)\ts$. 
But Theorem \ref{proposition4.5} implies
that the operator $\ts\xic_{\ts\si}$ on 
$\Jb\,\ts\backslash\,\Bb_m$ determines a linear operator
\begin{equation}
\label{bbjioper}
\Jb\,\ts\backslash\,\Bb_m\ts/\,\Ib_{\ts\mu,\de_+}
\to\,
\Jb\,\ts\backslash\,\Bb_m\ts/\,\Ib_{\ts\si\ts\circ\ts\mu\ts,\ts\si(\de_+)}
\ts.
\end{equation}
The latter operator intertwines the actions
of $\X(\g_n)$ on the source and the target vector spaces,
because the image of $\X(\g_n)$ in
$\B_m$ relative to $\be_m$ commutes with the subalgebra
$\U(\ts\f_m)\subset\B_m\ts$; see the definition \eqref{q1}.
We also use Lemma \ref{lemma41}, Part (ii).
Recall that $\F_m(V)=\F_{\ts\de_+}(V)\ts$.
Hence by using the equivalences \eqref{fbi} 
for the sequences $\de=\de_+$ and
$\de=\si\ts(\de_+)\ts$, the operator \eqref{bbjioper} becomes
the desired $\X(\g_n)\ts$-intertwining operator \eqref{distoper}.

As usual, for any $\f_m$-module $V$ and any element $\la\in\h^\ast$
let $V^\la\subset V$ be the subspace of vectors \textit{of weight\/}
$\la\ts$, so that any $X\in\h$ acts on $V^\la$
via multiplication by $\la\ts(X)\in\CC\ts$.
It now follows from the property \eqref{q121} of $\ts\xic_{\ts\si}$
that the restriction of
our operator \eqref{distoper} to the subspace of weight $\la$
is an $\X(\g_n)\ts$-intertwining operator  
\begin{equation}
\label{distoperla}	
\F_{m}\ts(\ts M_{\ts\mu})_{\ts\n}^{\ts\la}
\,\to\,
\F_{\ts\de}\ts(\ts M_{\ts\si\,\circ\ts\mu})_{\ts\n}^{\,\si\,\circ\ts\la}
\ \quad\text{where}\quad\ 
\de=\si\ts(\de_+)\ts.
\end{equation}

At the end of Section 2, we defined the modules $P_z$ and
$\Pp_z$ over the Yangian $\Y(\gl_n)\ts$. The underlying vector space
of these modules consists of all polynomial functions on 
$\CC^{\ts n}\ts$. Note that the action of $\Y(\gl_n)$ on each 
of these modules preserves the polynomial degree. Now for any 
$N=1\com2,\ldots$ denote respectively by $P_z^{\ts N}$ and
$P_z^{\ts-N}$ the submodules in $P_z$ and $\Pp_z$ which
consist of the polynomial functions of degree $N\ts$.
Note that $\Y(\gl_n)$ acts on the subspace of
constant functions in $P_z$ trivially, that is via the 
counit homomorphism $\Y(\gl_n)\to\CC\ts$.
That action of $\Y(\gl_n)$ does not depend on $z\ts$.
It will be still convenient to denote by $P_z^{\ts0}$ the
vector space $\CC$ with the trivial action of $\Y(\gl_n)\ts$.

Denote
\begin{equation}
\label{nua}
\nu_a=-\ts\frac{n}2+\mu_a-\la_a 
\quad\textrm{for}\quad
a=1\lcd m\ts.
\end{equation}
Suppose that $\nu_1\lcd\nu_m$ are non-negative integers,
otherwise the source $\X(\g_n)\ts$-module in 
\eqref{distoperla} would be zero thanks to Corollary \ref{verma}.
Under our assumption, Corollary \ref{verma} implies that
the the source $\X(\g_n)\ts$-module in \eqref{distoperla} is equivalent to  
\begin{equation}
\label{munup}
P_{\mu_m+z}^{\,\nu_m}
\ot 
P_{\mu_{m-1}+z+1}^{\,\nu_{m-1}}
\ot\ldots\ot 
P_{\mu_1+z+m-1}^{\,\nu_1}
\end{equation}
pulled back through the automorphism \eqref{fus} of $\X(\g_n)\ts$,
where $f(u)$ is given by \eqref{fuprod} and $z=\pm\tts\frac12\ts$.
A more general results is stated as Proposition \ref{siverma} below.
The tensor product in \eqref{munup} is that of
$\Y(\gl_n)\ts$-modules. Then we employ the embedding 
$\Y(\g_n)\subset\Y(\gl_n)$ and the homomorphism 
$\X(\g_n)\to\Y(\g_n)$ defined by \eqref{xy}.
By using the labels $\rho_1\lcd\rho_m$
of the halfsum $\rho$ of the positive roots of $\f_m\ts$,  
the tensor product \eqref{munup} can be rewritten as
\begin{equation}
\label{pp}
P_{\mu_m-\frac12+\rho_m}^{\,\nu_m}
\ot\ldots\ot 
P_{\mu_1-\frac12+\rho_1}^{\,\nu_1}\ts.
\end{equation}
By using the labels $\rho_1\lcd\rho_m$ 
we can also rewrite the product \eqref{fuprod} as
\begin{equation}
\label{muprod}
\prod_{a=1}^m\,\ts
\frac{u+\mu_a-\frac12+\rho_a}{u+\mu_a+\frac12+\rho_a}\ .
\end{equation}

Let us now consider the target $\X(\g_n)\ts$-module in \eqref{distoperla}.
For each $a=1\lcd m$ denote
$$
\widetilde\mu_a=\mu_{\ts|\si^{-1}(a)|}\ts,
\quad 
\widetilde\nu_a=\nu_{\ts|\si^{-1}(a)|}\ts,
\quad
\widetilde\rho_a=\rho_{\ts|\si^{-1}(a)|}\ts.
$$
The above description of the source $\X(\g_n)\ts$-module in  
\eqref{distoperla} can now be generalized to the target $\X(\g_n)\ts$-module,
which depends on an arbitrary element $\si\in\H_m\ts$.

\begin{proposition*}
\label{siverma}
For $\de=\si\ts(\de_+)$ the $\X(\g_n)\ts$-module 
$\F_{\ts\de}\ts(\ts M_{\ts\si\,\circ\ts\mu})_{\ts\n}^{\,\si\,\circ\ts\la}$
is equivalent to the tensor product
\begin{equation}
\label{ppd}
P_{\widetilde\mu_m-\frac12+\widetilde\rho_m}^{\,\ts\de_m\ts\widetilde\nu_m}
\ot\ldots\ot 
P_{\widetilde\mu_1-\frac12+\widetilde\rho_1}^{\,\ts\de_1\ts\widetilde\nu_1}
\end{equation}
pulled back through the automorphism \eqref{fus} of\/ $\X(\g_n)$ 
where $f(u)$ equals\/ \eqref{muprod}.
\end{proposition*}

\begin{proof}
First consider the bimodule
$\F_m(\ts  M_{\ts\si\,\circ\ts\mu})_{\ts\n}$ of $\h$ and $\X(\g_n)\ts$.
By Corollary~\ref{verma}, this bimodule is equivalent to the tensor product
\begin{equation}
\label{bim1}
P_{\ts\de_m\ts\widetilde\mu_m-\frac12+\ts\de_m\ts\widetilde\rho_m}
\ot\ldots\ot 
P_{\ts\de_1\ts\widetilde\mu_1-\frac12+\ts\de_1\ts\widetilde\rho_1}
\end{equation}
pulled back through the automorphism \eqref{fus} of $\X(\g_n)$ where
$f(u)$ equals
\begin{equation}
\label{deprod}
\prod_{a=1}^m\,\ts\frac
{u+\ts\de_a\ts\widetilde\mu_a-\frac12+\ts\de_a\ts\widetilde\rho_a}
{u+\ts\de_a\ts\widetilde\mu_a+\frac12+\ts\de_a\ts\widetilde\rho_a}\ .
\end{equation}
For any $a=1\lcd m$ the element $F_{-\ab,-\ab}\in\h$ acts on 
the tensor product \eqref{bim1} as
$$
-\ts\frac{n}2\tts-\tts\deg{\nns}_a+(\si\circ\mu)_{\tts a}
$$
where $\deg{\nns}_a$ is the degree operator
on the $a\ts$-th tensor factor, counting the factors from right
to left. It acts on the vector space $\P\ts(\CC^{\ts n})$ of that
tensor factor as the Euler operator
\begin{equation}
\label{euler}
\sum_{k=1}^n\,x_k\ts\d_k\in\PD\ts(\CC^{\ts n})\,.
\end{equation}

A bimodule equivalent to
$\F_\de\ts(\ts  M_{\ts\si\,\circ\ts\mu})_{\ts\n}$ can be obtained
by pushing forward the actions of $\h$ and $\X(\g_n)$ on \eqref{bim1}
through the composition of automorphisms \eqref{onefour},
for every tensor factor with number $a\ts$ such that $\de_a=-1$.
Here we number the $m$ tensor factors of \eqref{bim1} by $1\lcd m$
from right to left. Then we also have to pull the resulting 
$\X(\g_n)\ts$-module back through the automorphism 
\eqref{fus}, where the series $f(u)$ equals \eqref{deprod}.
The automorphism \eqref{onefour} maps the element
\eqref{euler}~to
$$
-\,\sum_{k=1}^n\,\d_{\ts\bk}\ts x_{\ts\bk}=
-\ts n-\,\sum_{k=1}^n\,x_k\ts\d_k\ts.
$$
Hence if $\de_a=-1\ts$, the element 
$F_{-\ab,-\ab}\in\h$ acts on the modified tensor product as
$$
\frac{n}2+(\si\circ\mu)_a+\deg{\nns}_a\ts.
$$
By equating the last displayed expression to $(\si\circ\la)_a$
and by using \eqref{nua} together with the condition $\de_a=-1\ts$,
we get the equation 
$\deg{\nns}_a=\widetilde\nu_a\ts$.
But by Lemma~\ref{ppl}, pushing the $\Y(\gl_n)\ts$-module 
$$
P_{-\widetilde\mu_a-\frac12-\widetilde\rho_a}^{\,\ts\widetilde\nu_a}
$$
forward
through the automorphism \eqref{onefour} of $\PD\ts(\CC^{\ts n})$
yields the same $\Y(\gl_n)\ts$-module as pulling
$$
P_{\widetilde\mu_a-\frac12+\widetilde\rho_a}^{\ts-\ts\widetilde\nu_a}
$$
back through the automorphism \eqref{fut} of $\Y(\gl_n)$ where 
$$
g(u)=
\frac{u+\widetilde\mu_a-\frac12+\widetilde\rho_a}
{u+\widetilde\mu_a+\frac12+\widetilde\rho_a}\ .
$$
Thus the $\X(\g_n)\ts$-module 
$\F_{\ts\de}\ts(\ts M_{\ts\si\,\circ\ts\mu})_{\ts\n}^{\,\si\,\circ\ts\la}$
is equivalent to the tensor product \eqref{ppd}
pulled back through the automorphism \eqref{fus}
where the series $f(u)$ is obtained by multiplying \eqref{deprod}
by $g(-u)\ts g(u)$ for each index $a\ts$ such that $\de_a=-1\ts$;
see the definition \eqref{xy}. But for any 
the element $\si\in\H_m$ the product \eqref{muprod} equals
\begin{equation}
\label{tildaprod}
\prod_{a=1}^m\,\ts
\frac{u+\widetilde\mu_a-\frac12+\widetilde\rho_a}
{u+\widetilde\mu_a+\frac12+\widetilde\rho_a}\ .
\end{equation}
If $\de_a=-1$ then the factors of \eqref{deprod} and \eqref{tildaprod} 
indexed by $a$ are equal to $g(-u)^{-1}$ and $g(u)$ respectively.
If $\de_a=1$ then the factors of \eqref{deprod} and \eqref{tildaprod}
indexed by $a$ coincide. 
This comparison of \eqref{deprod} and \eqref{tildaprod}
completes the proof. 
\end{proof}

The vector spaces of two equivalent
$\X(\g_n)\ts$-modules in Proposition \ref{siverma} are
$$
(M_{\ts\si\,\circ\ts\mu}\ot
\P\ts(\CC^{\ts m}\ot\CC^{\ts n}))_{\ts\n}^{\,\si\,\circ\ts\la}
\ \quad\text{and}\ \quad
\P^{\,\widetilde\nu_m}\tts(\CC^{\ts n})
\ot\ldots\ot 
\P^{\,\widetilde\nu_1}\tts(\CC^{\ts n})
$$
respectively.
Define a linear map from the latter vector space
to the former,
by mapping $f_1\ot\ldots\ot f_m$ to
the class of 
$1_{\ts\si\,\circ\ts\mu}\ot f$ in
the space of $\n\ts$-coinvariants. Here
$$
f_1\in\P^{\,\widetilde\nu_m}\tts(\CC^{\ts n})
\,\lcd\,
f_m\in\P^{\,\widetilde\nu_1}\tts(\CC^{\ts n})
$$
and the polynomial
$f\in\P\ts(\CC^{\ts m}\ot\CC^{\ts n})$ is defined by \eqref{fff}.
This linear map is an equivalence of the $\X(\g_n)\ts$-modules
in Proposition \ref{siverma},
see the remarks made after our proof of Corollary~\ref{verma}. 

Thus for any non-negative integers $\nu_1\lcd\nu_m$
we have demonstrated how the Zhelobenko operator $\ts\xic_{\ts\si}$
on $\Jb\,\ts\backslash\,\Bb_m$ determines an intertwining operator
between the $\X(\g_n)\ts$-modules \eqref{pp} and \eqref{ppd}
pulled back via the automorphism \eqref{fus} of $\X(\g_n)\ts$,
where the series $f(u)$ is the same \eqref{muprod} for both modules.
The same operator intertwines the $\X(\g_n)\ts$-modules
\begin{equation}
\label{ppi}
P_{\mu_m-\frac12+\rho_m}^{\,\nu_m}
\ot\ldots\ot 
P_{\mu_1-\frac12+\rho_1}^{\,\nu_1}
\to\,
P_{\widetilde\mu_m-\frac12+\widetilde\rho_m}^{\,\ts\de_m\ts\widetilde\nu_m}
\ot\ldots\ot 
P_{\widetilde\mu_1-\frac12+\widetilde\rho_1}^{\,\ts\de_1\ts\widetilde\nu_1}
\ts,
\end{equation}
neither of them being pulled back via the authomorphism \eqref{fus}.
It was proved in \cite{MN} that
both $\X(\g_n)\ts$-modules in \eqref{ppi} are irreducible
under our assumptions on $\mu\ts$.
Hence an intertwining operator between them 
is unique up to a multiplier from $\CC\ts$.
For our intertwining operator,
this multiplier is explicitly determined by Proposition~\ref{isis}
below. When $n=1\ts$, the $\Y(\gl_n)\ts$-module
$P^{\ts N}_z$ is one-dimensional for any integer $N$ and $z\in\CC\ts$. 
The operator \eqref{ppi} is then a multiplication by scalar.

Our proof of Proposition \ref{isis} will be based on the following
four lemmas. The proof of the first lemma is similar to the
proof of the second one, and is omitted. Let $s\com
t=0\com1\com2\com\ts\ldots$ and $k=1\lcd n\ts$. 
When $n>1$, we will suppose that $k\neq\bk\ts$.

\begin{lemma*}
\label{ddnorma} 
For any\/ $a=1\lcd m-1$ the
operator\/ $\xic_{\ts a}$ on\/ $\Jb\,\ts\backslash\,\Bb_m$ maps the
image in\/ $\Jb\,\ts\backslash\,\Bb_m$ of\/
$
\d_{\ts\overline{a}\,\bk}^{\,s}\,
\d_{\ts\overline{a+1}\,\bk}^{\,t}
\in\Bb_m
$
to the image in\/ $\Jb\,\ts\backslash\,\Bb_m$ of the product
$$
\widetilde{\si}_a\ts(\ts 
\d_{\ts\overline{a}\,\bk}^{\,s}\,
\d_{\ts\overline{a+1}\,\bk}^{\,t}\ts)
\,\cdot\,
\prod_{r=1}^t\,\frac{H_a+r+1}{\ts H_a+r-s}
$$
plus images in\/ $\Jb\,\ts\backslash\,\Bb_m$ of elements of the
left ideal in\/ $\Bb_m$ generated by\/ $\Jb'$ and \eqref{xxideal}.
\end{lemma*}

\begin{lemma*}
\label{xxnorma} 
For any\/ $a=1\lcd m-1$ the
operator\/ $\xic_{\ts a}$ on\/ $\Jb\,\ts\backslash\,\Bb_m$ maps the
image in\/ $\Jb\,\ts\backslash\,\Bb_m$ of\/ 
$
x_{\ts\overline{a}\ts k}^{\,s}\,
x_{\ts\overline{a+1}\ts k}^{\,t}
\in\Bb_m
$ 
to the image in\/ $\Jb\,\ts\backslash\,\Bb_m$ of the product
$$
\widetilde{\si}_a\ts(\ts 
x_{\ts\overline{a}\ts k}^{\,s}\,
x_{\ts\overline{a+1}\ts k}^{\,t}\ts)
\,\cdot\,
\prod_{r=1}^s\,\frac{H_a+r+1}{\ts H_a+r-t}
$$
plus images in\/ $\Jb\,\ts\backslash\,\Bb_m$ of elements of the
left ideal in\/ $\Bb_m$ generated by\/ $\Jb'$ and  \eqref{ddideal}.
\end{lemma*}

\begin{proof}
Let us use the symbol $\,\equiv\,$ to indicate equalities in 
the vector space 
$\Jb\,\ts\backslash\,\Bb_m$ modulo the subspace, which
is the image of the left ideal in $\Bb_m$ generated by $\Jb'$ and
the elements \eqref{ddideal}. By the definition \eqref{fourier}, we have
$$
\widetilde{\si}_a\ts(\ts 
x_{\ts\overline{a}\ts k}^{\,s}\,
x_{\ts\overline{a+1}\ts k}^{\,t}\ts)
\,=\,
x_{\ts\overline{a}\ts k}^{\,t}\,
x_{\ts\overline{a+1}\ts k}^{\,s}\,.
$$
Hence the operator $\xic_{\ts a}$ maps the
image in $\Jb\,\ts\backslash\,\Bb_m$ of 
$
x_{\ts\overline{a}\ts k}^{\,s}\,
x_{\ts\overline{a+1}\ts k}^{\,t}
\in\Bb_m$ to that of
\begin{align*}
\xi_{\ts a}\ts(\ts 
x_{\ts\overline{a}\ts k}^{\,t}\,
x_{\ts\overline{a+1}\ts k}^{\,s}\ts) 
&\,\,=\,\,
\sum_{r=0}^\infty\,\, 
(\ts r\ts!\,H_a^{\ts(r)}\ts)^{-1}\ts E_a^{\ts r}\, 
\widehat{F}_a^{\ts r}(\ts x_{\ts\overline{a}\ts k}^{\,t}\,
x_{\ts\overline{a+1}\ts k}^{\,s}\ts)
\\
&\,\,\equiv\,\, 
\sum_{r=0}^s\, 
\frac{\ts s\ldots(s-r+1)\,(t+1)\ldots\ts(t+r)}{r\ts!\,H_a\ldots(H_a-r+1)}\,\,
x_{\ts\overline{a}\ts k}^{\,t}\, 
x_{\ts\overline{a+1}\ts k}^{\,s}
\end{align*}
where we use the proof of \cite[Proposition 3.7]{KN1}. In the last
line, the sum of the fractions corresponding to $r=0\lcd s$ is a
particular value $\mathrm{F}\ts(-s\com\ts t+1\com\ts-H_a\,;1)$
of the hypergeometric function. 
Using the well known formula 
\begin{equation}
\label{gaga} 
\mathrm{F}\ts(u\com v\com w\,;1)=
\frac{\,\mathrm\Gamma(w)\,\mathrm\Gamma(w-u-v)}
{\,\mathrm\Gamma(w-u)\,\mathrm\Gamma(w-v)}
\end{equation}
valid for any $u\com v\com w\in\CC$ such that $w\neq0,-1,\ldots$
and $\mathrm{Re}\,(w-u-v)>0\ts$, we get
$$
\mathrm{F}\ts(-s\com v\com w\,;1)= \prod_{r=1}^s\,
\frac{w-v+r-1}{w+r-1}
$$
for any $w\neq0,-1,\ldots$ and any $v\ts$. Therefore
\begin{align*}
\xi_{\ts a}\ts(\ts x_{\ts\overline{a}\ts k}^{\,t}\,
x_{\ts\overline{a+1}\ts k}^{\,s}\ts) 
&\,\,\equiv\,\,
\prod_{r=1}^s\, \frac{H_a+t-s+r+1}{H_a-s+r}
\,\,\cdot\,
x_{\ts\overline{a}\ts k}^{\,t}\, 
x_{\ts\overline{a+1}\ts k}^{\,s}
\\[2pt]
&\,\,=\,\, 
x_{\ts\overline{a}\ts k}^{\,t}\,
x_{\ts\overline{a+1}\ts k}^{\,s}
\,\cdot\,
\prod_{r=1}^s\,\frac{H_a+r+1}{\ts H_a+r-t}
\end{align*}
as required. Here we also used the relation
$$
H_a\,x_{\ts\overline{a}\ts k}^{\,t}\, 
x_{\ts\overline{a+1}\ts k}^{\,s}
\,=\, 
x_{\ts\overline{a}\ts k}^{\,t}\,
x_{\ts\overline{a+1}\ts k}^{\,s}\,(H_a+s-t)
$$
which follows from \eqref{defar}, since
$$
\zeta_n(H_a) \,=\, \zeta_n\ts(\ts
F_{\ts\overline{a+1}\ts,\ts\overline{a+1}}\,-F_{\ts\ab\ts\ab}\ts)
\,=\, \sum_{i=1}^n\, (\ts x_{\ts\overline{a+1}\ts i}\,
\d_{\ts\overline{a+1}\ts i}\, - x_{\ts\overline{a}\ts i}\,
\d_{\ts\overline{a}\ts i}\ts )\ts.
\eqno\qedhere
$$
\end{proof}

\begin{lemma*}
\label{xdnorma}
For any\/ $a=1\lcd m-1$ the
operator\/ $\xic_{\ts a}\!$ on\/ $\Jb\,\ts\backslash\,\Bb_m$ maps
the image in\/ $\Jb\,\ts\backslash\,\Bb_m$ of\/
$
x_{\ts\overline{a}\ts k}^{\,s}\,
\d_{\ts\overline{a+1}\,\bk}^{\,t}
\in\Bb_m
$ 
to the image in\/ $\Jb\,\ts\backslash\,\Bb_m$ of 
$$
\widetilde{\si}_a\ts(\ts 
x_{\ts\overline{a}\ts k}^{\,s}\,
\d_{\ts\overline{a+1}\,\bk}^{\,t}\ts)
\,\cdot\, 
\left\{\begin{array}{cc}
\displaystyle
\prod_{r=1}^s\,\dfrac{H_a+r}{\ts H_a+r+t}
&\quad\text{if}\quad n=1\quad\text{and}\quad k=1\ts,
\\[12pt]
1
&\quad\text{if}\quad n>1\quad\text{and}\quad k\neq\bk\ts,
\end{array}\right.
$$
plus images in\/
$\Jb\,\ts\backslash\,\Bb_m$ of elements of the left ideal in\/
$\Bb_m$ generated by $\Jb'$ and \eqref{xdideal}.
\end{lemma*}

\begin{proof}
By the definition \eqref{fourier}, we have
$$
\widetilde{\si}_a\ts(\ts 
x_{\ts\overline{a}\ts k}^{\,s}\,
\d_{\ts\overline{a+1}\,\bk}^{\,t}\ts)
\,=\,
\d_{\ts\overline{a}\,\bk}^{\,t}\,
x_{\ts\overline{a+1}\ts k}^{\,s}\,.
$$
Hence 
$\xic_{\ts a}$ maps the image in $\Jb\,\ts\backslash\,\Bb_m$ of the element 
$
x_{\ts\overline{a}\ts k}^{\,s}\,
\d_{\ts\overline{a+1}\,\bk}^{\,t}
\in\Bb_m
$ 
to that of
\begin{equation*}
\xi_{\ts a}\ts(\ts \d_{\ts\overline{a}\,\bk}^{\,t}\,
x_{\ts\overline{a+1}\ts k}^{\,s}\ts) 
\,\,=\,\,
\sum_{r=0}^\infty\,\, 
(\ts r\ts!\,H_a^{\ts(r)}\ts)^{-1}\ts E_a^{\ts r}\, 
\widehat{F}_a^{\ts r}(\ts \d_{\ts\overline{a}\,\bk}^{\,t}\,
x_{\ts\overline{a+1}\ts k}^{\,s}\ts)\,.
\end{equation*}
Here
$$
\widehat{F}_a(\ts
\d_{\ts\overline{a}\,\bk}^{\,t}\,
x_{\ts\overline{a+1}\ts k}^{\,s}
\ts) 
\,=\, 
-\,\sum_{i=1}^n\,
[\,\ts
x_{\overline{a}\ts i}\,
\d_{\overline{a+1}\ts i}\,,\ts
\d_{\ts\overline{a}\,\bk}^{\,t}\, 
x_{\ts\overline{a+1}\ts k}^{\,s}\,]\ts.
$$ 
If $\bk\neq k\ts$, then each commutator in the above sum
equals zero modulo the left ideal generated by the elements \eqref{xdideal}.
This proves the lemma for $k\neq\bk\ts$.

Let $n=1$, then $k=\bk=1\ts$.
Let us now indicate by $\,\equiv\,$ the equalities in 
$\Jb\,\ts\backslash\,\Bb_m$ modulo 
the image of the left ideal in $\Bb_m$ generated by $\Jb'$ and by
\eqref{xdideal}. Then
\begin{eqnarray*}
\xi_{\ts a}\ts(\ts
\d_{\ts\overline{a}\,1}^{\,t}\,
x_{\ts\overline{a+1}\,1}^{\,s}
\ts)
&\equiv&
\!\!\!\!\sum_{r=0}^{\min(s,t)}\,
\frac
{(-1)^r\,s\ldots(s-r+1)\,t\ldots(t-r+1)}
{r\ts!\,H_a\ldots(H_a-r+1)}\,\,
\d_{\ts\overline{a}\,1}^{\,t}\,
x_{\ts\overline{a+1}\,1}^{\,s}
\\[6pt]
&=&
\mathrm{F}\ts(-s\com\ts-\ts t\com\ts-H_a\,;1)\,\,
\d_{\ts\overline{a}\,1}^{\,t}\,
x_{\ts\overline{a+1}\,1}^{\,s}\,
\\[6pt]
&=&
\prod_{r=1}^s\,\frac{H_a-s-t+r}{\ts H_a-r+1}\,
\cdot\,
\d_{\ts\overline{a}\,1}^{\,t}\,
x_{\ts\overline{a+1}\,1}^{\,s}
\\[2pt]
&=& 
\d_{\ts\overline{a}\,1}^{\,t}\,
x_{\ts\overline{a+1}\,1}^{\,s}\,
\cdot\,
\prod_{r=1}^s\,\frac{H_a+r}{\ts H_a+r+t}
\end{eqnarray*}
as required. Here we also used a relation which follows from \eqref{defar},
$$
H_a\,
\d_{\ts\overline{a}\,1}^{\,t}\,
x_{\ts\overline{a+1}\,1}^{\,s}\,
\,=\,
\d_{\ts\overline{a}\,1}^{\,t}\,
x_{\ts\overline{a+1}\,1}^{\,s}\,
(H_a+s+t)\,.
\eqno\qedhere
$$
\end{proof}

\begin{lemma*}
\label{xnorma} 
If $\,\f_m=\sp_{2m}\ts$, then the
operator\/ $\xic_{\ts m}$ on\/ $\Jb\,\ts\backslash\,\Bb_m$ maps the
image in\/ $\Jb\,\ts\backslash\,\Bb_m$ of\/
$x_{\ts1k}^{\,s}\in\Bb_m$ to the image in\/
$\Jb\,\ts\backslash\,\Bb_m$ of
$$
\widetilde{\si}_m\ts(\ts x_{\ts1k}^{\,s}\ts)
\,\cdot\, 
\left\{\begin{array}{cc}
\displaystyle
\prod_{r=1}^{[\ts s/2\ts]}\,\frac{H_m+r+1/2}{\ts H_m+s-r+1}
&\quad\text{if}\quad n=1\quad\text{and}\quad k=1\ts,
\\[12pt]
1
&\quad\text{if}\quad n>1\quad\text{and}\quad k\neq\bk\ts,
\end{array}\right.
$$
plus images in\/
$\Jb\,\ts\backslash\,\Bb_m$ of elements of the left ideal in\/
$\Bb_m$ generated by $\Jb'$ and \eqref{xideal}.
Here\/ $[\ts s/2\ts]=t$ whenever $s=2\ts t$ or\/ $s=2\ts t+1\ts$.
\end{lemma*}

\begin{proof}
By \eqref{fourier},
$$
\widetilde{\si}_m\ts(\ts x_{\ts1k}^{\,s}\ts)
\,=\,(-1)^{\ts s}\,\d_{\ts1\bk}^{\,s}\,.
$$
Hence 
$\xic_{\ts m}$ maps the image in $\Jb\,\ts\backslash\,\Bb_m$ of the element 
$x_{\ts1k}^{\,s}\,\in\Bb_m$ to that of 
\begin{equation*}
(-1)^{\ts s}\,
\xi_{\ts m}\ts(\ts\d_{\ts1\bk}^{\,s}\ts) 
\,\,=\,\,
(-1)^{\ts s}\,
\sum_{r=0}^\infty\,\, 
(\ts r\ts!\,H_m^{\ts(r)}\ts)^{-1}\ts E_m^{\ts r}\, 
\widehat{F}_m^{\ts r}\ts(\ts\d_{\ts1\bk}^{\,s}\ts)\,.
\end{equation*}
Here
$$
\widehat{F}_m(\ts\d_{\ts1\bk}^{\,s}\ts) 
\,=\, 
-\,\sum_{i=1}^n\,
[\,\ts
x_{\ts1\bi}\,x_{\ts1i}\,,\ts
\d_{\ts1\bk}^{\,s}\,]\,/\,2\,.
$$
If $\bk\neq k\ts$, then each commutator in the above sum
equals zero modulo the left ideal generated by the elements \eqref{xideal}.
This proves the lemma for $k\neq\bk\ts$.

Let $n=1$, then $k=\bk=1\ts$.
Let us now indicate by $\,\equiv\,$ the equalities in 
$\Jb\,\ts\backslash\,\Bb_m$ modulo 
the image of the left ideal in $\Bb_m$ generated by $\Jb'$ and by
\eqref{xideal}. Then
\begin{eqnarray*}
(-1)^s\,\xi_{\ts m}\ts(\ts\d_{\ts11}^{\,s}\ts)
&\!\!\equiv\!\!&
(-1)^s\,\sum_{r=0}^{[s/2\ts]}\,
\frac
{(-1)^r\,s\ts(s-1)\ldots(s-2r+2)\ts(s-2r+1)}
{4^r\,r\ts!\,H_m\ldots(H_m-r+1)}\,\,
\d_{\ts11}^{\,s}
\\[6pt]
&\!\!=\!\!&
(-1)^s\,\,
\mathrm{F}\ts(-s/2\com\ts(1-s)/2\com\ts-H_m\,;1)\,\,
\d_{\ts11}^{\,s}\,.\phantom{\sum}
\end{eqnarray*}
But the identity \eqref{gaga} also implies that for any $w\not=0,-1,\ldots$
$$
\mathrm{F}\ts(-s/2\com\ts(1-s)/2\com\ts w\,;1)\,\,=
\prod_{r=1}^{[\ts s/2\ts]}\,\frac{w+s-r-1/2}{\ts w+r-1}\ .
$$
Hence
\vspace{-4pt}
\begin{eqnarray*}
(-1)^s\,\xi_{\ts m}\ts(\ts\d_{\ts11}^{\,s}\ts)
&\!\!\equiv\!\!&
(-1)^s\,
\prod_{r=1}^{[\ts s/2\ts]}\,\frac{H_m-s+r+1/2}{\ts H_m-r+1}
\,\cdot\,
\d_{\ts11}^{\,s}
\\[2pt]
&\!\!=\!\!& 
(-1)^s\,\d_{\ts11}^{\,s}
\,\cdot\,
\prod_{r=1}^{[\ts s/2\ts]}\,\frac{H_m+r+1/2}{\ts H_m+s-r+1}
\end{eqnarray*}
as required. Here we also used the relation
$
H_m\,\d_{\ts11}^{\,s}\,=\,\d_{\ts11}^{\,s}\,(H_m+s)
$
which follows from \eqref{defar}, since for $\f_m=\sp_{2m}$ and $n=1$ we have
$$
\zeta_n(H_m)
\,=\,
-\,\zeta_n\ts(F_{11})
\,=\,
-\,\frac{1}2\,-\,
x_{11}\,\d_{11}\,.
\eqno\qedhere
$$
\end{proof}

Let us now state Proposition \ref{isis}.
We keep assuming that the weight $\mu$ satisfies the conditions
\eqref{notinso}, and also satisfies the conditions \eqref{notinsp}
in the case $\f_m=\sp_{2m}\ts$. Let
$(\ts\mu_1^\ast\lcd\mu_m^\ast\ts)$ be the sequence of labels of the
weight $\mu+\rho\ts$. Thus for each $a=1\lcd m$ we have
$\mu_a^\ast=\mu_a+m-a$ if $\f_m=\so_{2m}\ts$, and
$\mu_a^\ast=\mu_a+m-a+1$ if $\f_m=\sp_{2m}\ts$. Let
$(\ts\la_1^\ast\lcd\la_m^\ast\ts)$ be the sequence of labels of
$\la+\rho\ts$. Suppose that for each 
$a=1\lcd m$ the number $\nu_a$ defined by \eqref{nua} is a
non-negative integer.
For each positive root $\eta\in\De^+$ define a number $z_\eta\in\CC$ by
setting
$$
z_\eta\,=\,
\left\{
\begin{array}{ll}
\hspace{3.5pt}
\displaystyle
\prod_{r=1}^{\nu_{\ts b}}\,\hspace{5pt}
\frac{\mu_{\ts b}^\ast-\mu_c^\ast-r}
{\lambda_{\ts b}^\ast-\lambda_c^\ast+r}
&\quad\text{if}\quad
\eta=\ep_{b}-\ep_c\,,
\\[16pt]
\hspace{3.5pt}
\displaystyle
\prod_{r=1}^{\nu_{\ts b}}\hspace{5pt}
\frac{\mu_{\ts b}^\ast+\mu_c^\ast-r}
{\lambda_{\ts b}^\ast+\lambda_c^\ast+r}
&\quad\text{if}\quad
\eta=\ep_{b}+\ep_c\,,
\\[12pt]
\displaystyle
\prod_{r=1}^{[\nu_{\ts b}/2]}\,
\frac{\mu_{\ts b}^\ast-r}{\lambda_{\ts b}^\ast+r}
&\quad\text{if}\quad
\eta=2\ts\ep_{b}\,.
\end{array}
\right.
$$
In the first two cases above we have $1\le b<c\le m\ts$,
while in the third case we have $1\le b\le m$ and $\f_m=\sp_{2m}\ts$.

Denote by $\De^0$ 
the set of \textit{compact roots\/} of $\f_m\ts$,
these are the weights $\ep_b-\ep_c$
where $1\le b\com c\le m$ and $b\neq c\ts$. 
For any element $\si\in\H_m$ consider the set
$$
\De_{\ts\si}^0=\De_{\ts\si}\cap\De^0=
\{\ts\eta\in\De^+\cap\De^0\,|\,\si\ts(\eta)\notin\De^+\ts\}\,.
$$

If $n>1\ts$, choose any index $k\in\{1\lcd n\}$ such that $k\neq\bk\ts$,
that is $k<n$ when $n$ is odd.
If $n=1\ts$, then we must have $k=\bk=1\ts$.
Denote by $v_{\mu}^{\ts\la}$ the image in
$\Jb\,\ts\backslash\,\Bb_m\ts/\,\Ib_{\ts\mu,\de_+}$ of the product
$
x_{\ts\overline{1}\ts k}^{\,\nu_1}
\,\ldots\, 
x_{\ts\overline{m}\ts k}^{\,\nu_m}
\in\Bb_m\ts.
$

\begin{proposition*}
\label{isis} 
{\rm\,\,(i)} 
The vector $v_{\mu}^{\ts\la}$ is not in the zero coset of
$\Jb\,\ts\backslash\,\Bb_m\ts/\,\Ib_{\ts\mu,\de_+}\ts$.
\\
{\rm\,\,(ii)}
Under the action of\/ $\h$ on\/ 
$\Jb\,\ts\backslash\,\Bb_m\ts/\,\Ib_{\ts\mu,\de_+}$ the vector
$v_{\mu}^{\ts\la}$ is of weight $\lambda\ts$.
\\
{\rm\,\,(iii)} 
For any\/ $\si\in\H_m$ the operator \eqref{bbjioper}
determined by\/ $\xic_{\ts\si}$ maps the vector\/ $v_{\mu}^{\ts\la}$
to the image in\/
$\Jb\,\ts\backslash\,\Bb_m\ts/\,\Ib_{\ts\si\ts\circ\ts\mu\ts,\ts\si(\de_+)}$
of\/ $ \sih\ts(\ts x_{\ts\overline{1}\ts k}^{\,\nu_1} \,\ldots\,
x_{\ts\overline{m}\ts k}^{\,\nu_m} \ts) \in\Bb_m $ multiplied by
\begin{equation}
\label{isim} 
\prod_{\eta\ts\in\De_\si}\!z_\eta
\quad\ \,\text{or}\ \quad
\prod_{\eta\ts\in\De_\si^0}\!z_\eta
\end{equation}
when $n=1$ or $n>1$ respectively.
\end{proposition*}

\begin{proof}
Part (i) follows directly from the definition of
the ideal $\Ib_{\ts\mu,\de_+}\ts$. Let us prove Part (ii). The
elements of $\h$ act on
$\Jb\,\ts\backslash\,\Bb_m\ts/\,\Ib_{\ts\mu,\de_+}$ via their left
multiplication on $\Bb_m\ts$. Let us indicate by $\,\equiv\,$ the
equalities in $\Bb_m$  modulo the left ideal
$\Ib_{\ts\mu,\de_+}\ts$. Then by the definition \eqref{gan} for
$a=1\lcd m$ we have the relations in the algebra $\Bb_m\ts$,
\begin{gather*}
F_{-\ab,-\ab}\,\ts 
x_{\ts\overline{1}\ts k}^{\,\nu_1}
\,\ldots\,
x_{\ts\overline{m}\ts k}^{\,\nu_m}
\,=\,
x_{\ts\overline{1}\ts k}^{\,\nu_1}
\,\ldots\,
x_{\ts\overline{m}\ts k}^{\,\nu_m} 
\,F_{-\ab,-\ab} \,- \sum_{i=1}^n\, 
\bigl[\,
x_{\ts\overline{a}\ts i}\,
\d_{\ts\overline{a}\ts i}\ts 
+\frac{n}2\,\com\ts 
x_{\ts\overline{1}\ts k}^{\,\nu_1}
\,\ldots\,
x_{\ts\overline{m}\ts k}^{\,\nu_m} 
\,\bigr]
\\[4pt]
=\, 
x_{\ts\overline{1}\ts k}^{\,\nu_1}
\,\ldots\,
x_{\ts\overline{m}\ts k}^{\,\nu_m}
\,(\ts F_{-\ab,-\ab}-\nu_a\ts) 
\,\equiv\,
x_{\ts\overline{1}\ts k}^{\,\nu_1}
\,\ldots\,
x_{\ts\overline{m}\ts k}^{\,\nu_m}
\,(\,\zeta_n(F_{-\ab,-\ab})+\mu_a-\nu_a\ts)
\\[9pt]
\equiv\, 
x_{\ts\overline{1}\ts k}^{\,\nu_1}
\,\ldots\,
x_{\ts\overline{m}\ts k}^{\,\nu_m}
\,(\ts-\ts\frac{n}2+\mu_a-\nu_a\ts) 
\,=\,
\la_a\,
x_{\ts\overline{1}\ts k}^{\,\nu_1}
\,\ldots\,
x_{\ts\overline{m}\ts k}^{\,\nu_m}
\,.
\end{gather*}
As required, then
$$
F_{-\ab,-\ab}\,\ts v_{\mu}^{\ts\la}
\,=\,
\la_a\,v_{\mu}^{\ts\la}\, \
\quad\text{for}\ \quad a=1\lcd m\ts.
$$

We will prove Part (iii) by induction on the length of a reduced
decomposition of $\si$ in terms of $\si_1\lcd\si_m\ts$. If $\si$ is
the identity element of $\H_m\ts$, then the required statement is
tautological. Now suppose that for some $\si\in\H_m$ the statement
of (iii) is true. Take any simple reflection $\si_a\in\H_m$ with
$1\le a\le m\ts$, such that $\si_a\ts\si$ has a longer reduced
decomposition in terms of $\si_1\lcd\si_m$ than $\si\ts$. If
$\f_m=\so_{2m}$ and $a=m\ts$, then we have
$\xic_{\ts\si_m\ts\si}=\sih_m\,\xic_{\ts\si}$ and
$\De_{\ts\si_m\ts\si}=\De_{\ts\si}\,$, so that the induction step is
immediate. We may now assume that $a<m$ in the case
$\f_m=\so_{2m}\ts$.

Take the simple root $\eta_{\ts a}$ corresponding to the reflection
$\si_a\ts$. Let $\eta=\si^{-1}(\eta_{\ts a})\ts$. Then
$\eta\in\Delta^+$ and
$$
\si_a\ts\si\ts(\eta)=\si_a(\eta_{\ts a})=-\ts\eta_{\ts
a}\notin\De^+.
$$
Hence
$$
\De_{\ts\si_a\ts\si}=\De_{\ts\si}\sqcup\ts\{\eta\tts\}\ts.
$$
Let $\ka\in\h^\ast$ be the weight with labels
\eqref{kappa}. Using the proof of Theorem~\ref{proposition4.5}, we
get the equality of two left ideals of the algebra $\Bb_m\ts$,
$$
\Ib_{\ts(\ts\si_a\ts\si\ts)\ts\circ\ts\mu \ts,\ts
(\ts\si_a\ts\si\ts)\ts(\de_+)}\,=\,
\tilde\I_{\ts(\ts\si_a\ts\si\ts)\ts\circ\ts\ka \ts,\ts
(\ts\si_a\ts\si\ts)\ts(\de_+)}\ts.
$$
But modulo the second of these two ideals, the element $H_a$ equals
\begin{gather}
\nonumber ((\ts\si_a\ts\si\ts)\circ\ka\ts)\ts(H_a)\,=\,
(\ts\si_a\ts\si\ts(\ka+\rho)-\rho\ts)\ts(H_a)\,=\,
(\ka+\rho)\ts(\ts\si^{-1}\si_a(H_a))-\rho\ts(H_a)\,=\,
\\[4pt]
\label{rhs} -\,(\ka+\rho)\ts(\si^{-1}(H_a))-1\,=\,
-\,(\ka+\rho)\ts(H_\eta)-1\,=\,
-\,\frac{\,2\ts(\ts\ka+\rho\com\eta\ts)\,}{(\ts\eta\com\eta\ts)}\ts-1\ts.
\end{gather}
Here $H_\eta=\si^{-1}(H_a)$ is the coroot corresponding to the root
$\eta\,$, and we use the standard bilinear form on $\h^\ast\ts$.
Using only 
\eqref{kappa}, the right hand side of \eqref{rhs}~equals
\begin{align*}
-\mu_{\ts b}^\ast+\mu_c^\ast-1
\ \quad\textrm{if}\ \quad 
\eta&=\ep_{b}-\ep_c\ts;
\\[5pt]
-\mu_{\ts b}^\ast-\mu_c^\ast+n-1
\ \quad\textrm{if}\ \quad 
\eta&=\ep_{b}+\ep_c\ts;
\\[2pt]
-\mu_{\ts b}^\ast+\frac{n}{2}-1
\ \quad\textrm{if}\ \quad 
\eta&=2\ts\ep_{b}\ts.
\end{align*}
Let us use (iii) as the induction assumption.
Denote $\de=\si(\de_+)\ts$. Consider five~cases.

I. Suppose that $\eta=\ep_b-\ep_c$ where $1\le b<c\le m\ts$, while
$\si(\ep_b)=\ep_a$ and $\si(\ep_c)=\ep_{a+1}\ts$. Then
$\si_a=\ep_a-\ep_{a+1}$ and $\de_a=\de_{a+1}=1\ts$. Hence
$$
\sih\ts(\ts x_{\ts\overline{1}\ts k}^{\,\nu_1}
\,\ldots\,
x_{\ts\overline{m}\ts k}^{\,\nu_m}\ts)
\,=\, 
x_{\ts\overline{a}\ts k}^{\,\nu_b}\, 
x_{\ts\overline{a+1}\ts k}^{\,\nu_c}\,Y
$$ 
where $Y$ is a certain element of the subalgebra
of $\PD\ts(\CC^{\ts m}\ot\CC^{\ts n})$
generated by all $x_{dk}$ and $\d_{dk}$
with $d\neq\overline{a}\com\overline{a+1}\ts$.
Now Lemma \ref{xxnorma} with $s=\nu_{\ts b}$ and $t=\nu_c$
applies. By substituting $-\mu_{\ts b}^\ast+\mu_c^\ast-1$ for $H_a$
in the product over the indices $r=1\lcd s$ in that lemma, 
the product becomes
\begin{equation}
\label{sampro}
\prod_{r=1}^{\nu_b}\,
\frac{-\mu_b^\ast+\mu_c^\ast+r}
{-\mu_b^\ast+\mu_c^\ast+\nu_b-\nu_c-r}
\,=\, 
\prod_{r=1}^{\nu_{\ts b}}\,\hspace{5pt} 
\frac{\mu_{\ts b}^\ast-\mu_c^\ast-r} 
{\lambda_{\ts b}^\ast-\lambda_c^\ast+r}\,\,.
\end{equation}

II. Suppose that $\eta=\ep_b-\ep_c$ where $1\le b<c\le m\ts$, but
$\si(\ep_b)=-\ts\ep_{a+1}$ and $\si(\ep_c)=-\ts\ep_{a}\ts$. Then
$\si_a=\ep_a-\ep_{a+1}$ again, but $\de_a=\de_{a+1}=-\ts1\ts$. Hence
$$
\sih\ts(\ts x_{\ts\overline{1}\ts k}^{\,\nu_1} 
\,\ldots\,
x_{\ts\overline{m}\ts k}^{\,\nu_m} \ts)
\,=\, 
\d_{\ts\overline{a}\ts\bk}^{\,\nu_c}\,
\d_{\ts\overline{a+1}\,\bk}^{\,\nu_b}\,Y
$$
where $Y$ is another element of the subalgebra
of $\PD\ts(\CC^{\ts m}\ot\CC^{\ts n})$
generated by all $x_{dk}$ and $\d_{dk}$
with $d\neq\overline{a}\com\overline{a+1}\ts$.
Now Lemma \ref{ddnorma} with $s=\nu_{\ts c}$ and $t=\nu_{\ts b}$
applies. By substituting $-\mu_{\ts b}^\ast+\mu_c^\ast-1$
for $H_a$ in the product over $r=1\lcd t$ in that lemma, 
the product becomes the same number \eqref{sampro}.

III. Suppose that $\eta=\ep_b+\ep_c$ and $1\le b<c\le m\ts$, while
$\si(\ep_b)=\ep_a$ and $\si(\ep_c)=-\ts\ep_{a+1}\ts$. Then
$\si_a=\ep_a-\ep_{a+1}$ again, but $\de_a=1$ and $\de_{a+1}=-1\ts$.
Hence
$$
\sih\ts(\ts 
x_{\ts\overline{1}\ts k}^{\,\nu_1}
\,\ldots\,
x_{\ts\overline{m}\ts k}^{\,\nu_m} 
\ts)
\,=\,
x_{\ts\overline{a}\ts k}^{\,\nu_b}\,
\d_{\ts\overline{a+1}\,\bk}^{\,\nu_c}\,Y
$$
where $Y$ is yet another element of the subalgebra
of $\PD\ts(\CC^{\ts m}\ot\CC^{\ts n})$
generated by all $x_{dk}$ and $\d_{dk}$
with $d\neq\overline{a}\com\overline{a+1}\ts$.
Now Lemma \ref{xdnorma} with $s=\nu_{\ts b}$ and $t=\nu_c$ applies.  
When $n=1$, 
by substituting 
$-\mu_{\ts b}^\ast-\mu_c^\ast$ for $H_a$
in the product over 
$r=1\lcd s$ in that lemma, the product becomes
\begin{equation}
\label{samsam}
\prod_{r=1}^{\nu_b}\,
\frac{-\mu_b^\ast-\mu_c^\ast+r}
{-\mu_b^\ast-\mu_c^\ast+\nu_b+\nu_c+1-r}
\,=\,
\prod_{r=1}^{\nu_{\ts b}}\hspace{5pt}
\frac{\,\mu_{\ts b}^\ast+\mu_c^\ast-r}
{\lambda_{\ts b}^\ast+\lambda_c^\ast+r}\,\,.
\end{equation}

IV. Suppose that $\eta=\ep_b+\ep_c$ where $1\le b<c\le m\ts$, but
$\si(\ep_b)=-\ts\ep_{a+1}$ and $\si(\ep_c)=\ep_{a}\ts$. Then
$\si_a=\ep_a-\ep_{a+1}$ again, but $\de_a=1$ and
$\de_{a+1}=-\ts1\ts$. Hence
$$
\sih\ts(\ts x_{\ts\overline{1}\ts k}^{\,\nu_1}
\,\ldots\,
x_{\ts\overline{m}\ts k}^{\,\nu_m}\ts)
\,=\,
x_{\ts\overline{a}\ts k}^{\,\nu_c}\,
\d_{\ts\overline{a+1}\,\bk}^{\,\nu_b}\,Y
$$
where $Y$ is another element of the subalgebra
of $\PD\ts(\CC^{\ts m}\ot\CC^{\ts n})$
generated by all $x_{dk}$ and $\d_{dk}$
with $d\neq\overline{a}\com\overline{a+1}\ts$.
Now Lemma \ref{xdnorma} with $s=\nu_{\ts c}$ and $t=\nu_b$ applies.  
When $n=1$, 
by substituting 
$-\mu_{\ts b}^\ast-\mu_c^\ast$ for $H_a$
in the product over 
$r=1\lcd s$ in that lemma, the product becomes
\begin{equation*}
\prod_{r=1}^{\nu_c}\,
\frac{-\mu_b^\ast-\mu_c^\ast+r}
{-\mu_b^\ast-\mu_c^\ast+\nu_b+\nu_c+1-r}\,\,.
\end{equation*}
Then the latter product is equal to \eqref{samsam}. 

V. Suppose that $\f_m=\sp_{2m}$ and $\eta=2\ts\ep_{\ts b}$ with
$1\le b\le m\ts$. Then $\si(\ep_{\ts b})=\ep_{\ts m}$ and
$\si_a=\si_m\ts$, while $\de_m=1\ts$. Hence
$$
\sih\ts(\ts 
x_{\ts\overline{1}\ts k}^{\,\nu_1}
\,\ldots\,
x_{\ts\overline{m}\ts k}^{\,\nu_m}\ts)
\,=\,
x_{\ts\overline{m}\ts k}^{\,\nu_b}\,Y
$$
where $Y$ is an element of the subalgebra
of $\PD\ts(\CC^{\ts m}\ot\CC^{\ts n})$
generated by $x_{dk}$ and $\d_{dk}$
with $d\neq\overline{m}=1\ts$.
So Lemma \ref{xnorma} with $s=\nu_{\ts b}$ applies.
When $n=1$, 
by substituting 
$-\mu_{\ts b}^\ast+1/2$ for $H_m$
in the product over $r$ 
in that lemma, the product becomes
$$
\prod_{r=1}^{[\ts \nu_{\ts b}/2\ts]}\,
\frac{-\mu_{\ts b}^\ast+r}
{-\mu_{\ts b}^\ast+1/2+\nu_{\ts b}-r}
\,\,\,=\,
\prod_{r=1}^{[\nu_{\ts b}/2]}\,
\frac{\mu_{\ts b}^\ast-r}{\lambda_{\ts b}^\ast+r}\,\,.
$$

Thus in all the five cases considered above, the operator
$$
\Jb\,\ts\backslash\,\Bb_m\ts/\,\Ib_{\ts\mu,\de_+} \to\,
\Ib_{\ts(\ts\si_a\ts\si\ts)\ts\circ\ts\mu \ts,\ts
(\ts\si_a\ts\si\ts)\ts(\de_+)}
$$
determined by\/ $\xic_{\ts\si_a\ts\si}$ maps $v_{\mu}^{\ts\la}$ to
the image in $ \,\Jb\,\ts\backslash\,\Bb_m\ts/\,
\Ib_{\ts(\ts\si_a\ts\si\ts)\ts\circ\ts\mu \ts,\ts
(\ts\si_a\ts\si\ts)\ts(\de_+)} $
of 
$$
\sih_a\ts\sih\ts(\ts x_{\ts\overline{1}\ts k}^{\,\nu_1} \,\ldots\,
x_{\ts\overline{m}\ts k}^{\,\nu_m} \ts) \in\Bb_m
$$
multiplied by the product \eqref{isim},
and by the extra factor $z_\eta$
when $\eta$ is a compact root or $n=1\ts$. 
This observation makes the induction step.
\end{proof}

The product \eqref{isim} in Proposition \ref{isis}
does not depend on the choice of a reduced decomposition of $\si\in\H_m$
in terms of $\si_1\lcd\si_m\ts$.
The uniqueness of the
intertwining operator \eqref{ppi} thus provides another proof
of the independence of our operator \eqref{distoperla}
on the decomposition of $\si\ts$, 
not involving Proposition~\ref{p3}.
Proposition \ref{isis} also shows that our
intertwining operator \eqref{distoperla} is not zero.


\section{Olshanski homomorphism}
\setcounter{section}{6}
\setcounter{equation}{0}
\setcounter{theorem*}{0}

For a positive integer $l\ts$, take the vector space
$\CC^{\ts n+l}$. 
In the case of an alternating form on
$\CC^{\ts n}$ choose $l$ to be even.
Let $e_1\lcd e_{\ts n+l}$ 
be the vectors of the standard basis in $\CC^{\ts n+l}$.
Consider the decomposition 
$\CC^{\ts n+l}=\CC^{\ts n}\op\CC^{\ts l}$
where the direct summands $\CC^{\ts n}$ and $\CC^{\ts l}$
are spanned by the vectors $e_1\lcd e_n$ and $e_{\ts n+1}\lcd e_{\ts n+l}$
respectively. This defines an embedding of the direct sum 
$\gl_n\op\ts\gl_{\ts l}$ of Lie algebras to $\gl_{\ts n+l}\ts$.
As a subalgebra of $\gl_{\ts n+l}\ts$,
the summand $\gl_n$ is spanned by the matrix units 
$E_{ij}\in\gl_{\ts n+l}$ where
$i\com j=1\lcd n\ts$. The summand $\gl_{\ts l}$ is spanned by
the matrix units $E_{ij}$ where $i\com j=n+1\lcd n+l\ts$.

The subspace $\CC^{\ts n}\subset\CC^{\ts n+l}$ comes
with a bilinear form chosen in Section 1. Now
choose a bilinear form on the subspace $\CC^{\ts l}\subset\CC^{\ts n+l}$
in a similar way. Namely, let 
$i$ be any of the indices $n+1\lcd n+l\ts$. If $i-n$ is even, then put
$\bi=i-1\ts$. If $i-n$ is odd and $i<n+l$,  then put $\bi=i+1\ts$.
If $i=n+l$ and $l$ is odd, then put $\bi=i\ts$.
Further, put $\th_i=1$ or $\th_i=(-1)^{\ts i-n-1}$
in the case of the symmetric or alternating form on $\CC^{\ts n}$. 
For any basis vectors $e_i$ and $e_j$ of the subspace $\CC^{\ts l}$ put 
$\langle\ts e_i\com e_j\ts\rangle=\th_i\,\de_{\ts\bi j}\ts$.
Equip the vector space $\CC^{\ts n+l}$ with the
bilinear form which is the sum of the forms on the direct summands.
The forms on $\CC^{\ts l}$ and $\CC^{\ts n+l}$ are of the same type
(symmetric or alternating) as the form on $\CC^n\ts$.

Now consider the subalgebras $\g_n\com\ts\g_{\ts l}$ and $\g_{n+l}$ of
the Lie algebras $\gl_n\com\ts\gl_{\ts l}$ and $\gl_{\ts n+l}$ respectively.
We have an embedding of the direct sum $\g_n\op\g_l$ 
to the Lie algebra $\g_{n+l}\ts$. We also have an embedding
of the direct product of Lie groups $G_n\times G_l$ to $G_{n+l}\ts$. 
Let $\C_l$ denote the subalgebra of $G_l\ts $-invariants
in the universal enveloping algebra $\U(\g_{n+l})\ts$.
Then $\C_l$ contains the subalgebra $\U(\g_n)\subset\U(\g_{n+l})\ts$.
If $\g_n=\sp_n$ then $\C_l$ coincides with the
centralizer of the subalgebra $\U(\sp_l)\subset\U(\sp_{n+l})\ts$. 
If $\g_n=\so_n$ then $\C_l$ is contained in the
centralizer of $\U(\so_l)\subset\U(\so_{n+l})\ts$,
but may not coincide with the centralizer.

Take the extended twisted Yangian $\X(\g_{n+l})\ts$. 
The subalgebra of $\X(\gl_{n+l})$ generated~by 
$$
S_{ij}^{\ts(1)},S_{ij}^{\ts(2)},\ts\ldots
\quad\text{where}\quad
i\com j=1\lcd n 
$$
is isomorphic to $\X(\g_n)\ts$ as an associative algebra, see
\cite[Section 3.14]{MNO}.  Thus we have a natural embedding
$\X(\g_n)\to\X(\g_{n+l})\ts$, let us denote it by $\io_{\ts l}\ts$.
We also have a surjective homomorphism 
$$
\pi_{n+l}:\ts\X(\g_{n+l})\to\U(\g_{n+l})\ts,
$$
see \eqref{pin}. Note that the composition 
$\pi_{n+l}\,\io_{\ts l}$ coincides with the homomorphism~$\pi_n\ts$.

Further, consider the involutive automorphism $\om_{n+l}$
of the algebra $\X(\g_{n+l})\ts$, see the definition \eqref{sin}.
The image of the composition of homomorphisms
\begin{equation}
\label{ohom}
\pi_{n+l}\,\ts\om_{n+l}\,\ts\io_{\ts l}:\ts
\X(\g_n)\to
\U(\g_{n+l})
\end{equation}
belongs to subalgebra $\C_l\subset\U(\g_{n+l})\ts$.
Moreover, together with the subalgebra of
$G_{n+l}\ts$-invariants in $\U(\g_{n+l})\ts$,
this image generates 
$\C_l\ts$. These two results are due to G.\,Olshanski \cite{O2},
for their detailed proofs see \cite[Section 4]{MO}.
We will use the composition of the homomorphisms
$$
\ga_{\tts l}=
\pi_{n+l}\,\ts
\om_{n+l}\,\ts
\io_{\ts l}\,\ts
\om_n\ts.
$$
We will call it the \textit{Olshanski homomorphism\/}.
The images of the homomorphisms $\ga_{\tts l}$ and \eqref{ohom} 
in $\U(\g_{n+l})$ coincide.
The reason for using the homomorphism $\ga_{\tts l}$ rather than
the homomorphism \eqref{ohom} will become apparent when we state
Theorem~\ref{5.1}.

An irreducible representation of the group $G_n$ is
\textit{polynomial\/} if it appears as a subrepresentation
of some tensor power of the defining representation $\CC^{\ts n}$. 
According to \cite[Sections V.7 and VI.3]{W}
the irreducible polynomial representations of the group $G_n$ 
are parameterized
by all the partitions $\nu$ of $N=0\com1\com2\com\,\ldots$
such that $\nus_1+\nus_2\le n$ in the case $G_n=O_n\ts$, and
$2\ts\nus_1\le n$ in the case $G_n=Sp_n\ts$.
Here $\nus$ is the partition conjugate to $\nu$
while $\nus_1\ts\com\nus_2\ts\com\,\ldots$ are the parts of $\nus$.
We denote by $W_\nu$ the irreducible polynomial representation
of $G_n$ corresponding to $\nu\ts$.

Let $\nu_1\com\ts\nu_2\com\,\ldots$ be the parts of $\nu\ts$.
Let $\nut$ be the weight of $\f_m$ with the sequence of labels
$$
(\ts-\,\frac{n}2-\nu_m\ts\lcd-\,\frac{n}2-\nu_1\ts)\ts.
$$
Consider $\P\ts(\CC^{\ts m}\ot\CC^{\ts n})$ as a bimodule
over $\f_m$ and $G_n\ts$. Then by \cite[Section 6]{KV} 
when $G_n=Sp_n\ts$, or by \cite[Section 9]{EHW} 
when $G_n=O_n\ts$, we have a 
decomposition
\begin{equation}
\label{dirsum}
\P\ts(\CC^{\ts m}\ot\CC^{\ts n})\,=\,\,
\mathop{\op}\limits_\nu\,L_{\ts\nut}\ot W_\nu
\end{equation}
where $\nu$ ranges over all parameters of
irreducible polynomial representations of $G_n$ such that $\nus_1\le m\ts$.
Here $L_{\ts\nut}$ is the irreducible $\f_m\ts$-module
of the highest weight $\nut\ts$.

Let $\la$ and $\mu$ be parameters of any
irreducible polynomial representations of the groups $G_{n+l}$ and $G_l$
respectively. Suppose that $\las_1\com\ts\mus_1\le m\ts$.
Using the action of the group $G_l$ on $W_\la$ 
via its embedding to $G_{n+l}$
as the second direct factor of the subgroup 
$G_n\times G_l$ consider the vector space
\begin{equation}
\label{hll}
\Hom_{\,G_l}(\ts W_\mu\ts,W_\la)\ts.
\end{equation}
The subalgebra $\C_l\subset\U(\g_{n+l}\ts)$ acts on this vector
space through the action of $\U(\g_{n+l})$ on $W_\la\ts$. 
In the case $G_n=Sp_n\ts$,
the vector space \eqref{hll} is irreducible under the action of the
algebra $\C_l\ts$; see \cite[Theorem 9.1.12]{D}.
In the case $G_n=O_n\ts$, the $\C_l\ts$-module \eqref{hll}
is either irreducible or splits to a direct
sum of two irreducible $\C_l\ts$-modules.
It is irreducible if $W_\la$ is irreducible
as a $\so_{n+l}\ts$-module, that is if $2\ts\las_1\neq n+l$
by \cite[Section V.9]{W}.
Note that in the case $G_n=O_n$ the condition $2\ts\las_1\neq n+l$
is sufficient but not necessary for the irreducibility 
of the $\C_l\ts$-module \eqref{hll}\ts; see \cite[Section 1.7]{N2}.

In any case, \eqref{hll} is irreducible under the joint action of
the subalgebra $\C_l\subset\U(\g_{n+l})$ and of
the subgroup $G_n\subset G_{n+l}\ts$.
Hence the following identifications of modules over 
$C_l$ and $G_n$ are unique up to rescaling of their
vector spaces\ts:
$$
\Hom_{\,G_l}(\ts W_\mu\ts,W_\la)\ts=
$$
$$
\Hom_{\,G_l}(\ts W_\mu\ts,
\Hom_{\,\f_m}(\ts L_{\ts\lat}\,,\P\ts(\ts\CC^{\ts m}\ot\CC^{\ts n+l}\ts)))\ts=
\vspace{2pt}
$$
$$
\Hom_{\,G_l}(\ts W_\mu\ts,
\Hom_{\,\f_m}(\ts L_{\ts\lat}\,,
\P\ts(\ts\CC^{\ts m}\ot\CC^{\ts l}\ts)\ot
\P\ts(\ts\CC^{\ts m}\ot\CC^{\ts n}\ts)))\ts=
\vspace{2pt}
$$
\begin{equation}
\label{hlls}
\,\Hom_{\,\f_m}(\ts L_{\ts\lat}\,,
L_{\ts\mut}\ot\P\ts(\ts\CC^{\ts m}\ot\CC^{\ts n}\ts))\ts.
\vspace{4pt}
\end{equation}
We uses the decompositions \eqref{dirsum} for $n+l$
and $l$ instead of $n\ts$, and the identification
\begin{equation}
\label{pmnl}
\P\ts(\ts\CC^{\ts m}\ot\CC^{\ts n+l}\ts)=
\P\ts(\ts\CC^{\ts m}\ot\CC^{\ts l}\ts)\ot
\P\ts(\ts\CC^{\ts m}\ot\CC^{\ts n}\ts)\ts,
\end{equation}
so that the labels of the weights $\lat$ and $\mut$ of $\f_m$
are respectively
$$
(\ts-\,\frac{n+l}2-\la_m\ts\lcd-\,\frac{n+l}2-\la_1\ts)
\ \quad\text{and}\ \quad
(\ts-\,\frac{l}2-\mu_m\ts\lcd-\,\frac{l}2-\mu_1\ts)\ts.
$$

By pulling back via the Olshanski
homomorphism $\ga_{\tts l}:\X(\g_n)\to\C_l\ts$,
the vector space \eqref{hll} becomes a module over the 
extended twisted Yangian $\X(\g_n)\ts$. Using the above identifications,
the vector space \eqref{hlls} than also becomes a module over $\X(\g_n)\ts$.
But the target $\f_m\ts$-module
$L_{\ts\mut}\ot\P\ts(\ts\CC^{\ts m}\ot\CC^{\ts n}\ts)$
in \eqref{hlls} coincides with the $\f_m\ts$-module
$\F_m(L_{\ts\mut})\ts$.

\begin{theorem*}
\label{5.1}
The action of\/ $\X(\g_n)$ on the vector space \eqref{hlls}
via the homomorphism\/ $\ga_{\tts l}$ coincides with the action, obtained by
pulling the action of\/ $\X(\g_n)$ on the bimodule\/ 
$\F_m(L_{\ts\mut})$ back through the homomorphism \eqref{fus} where
\begin{equation}
\label{fuss}
f(u)=1\,+\,m\ts\Bigl(\ts u-\frac{l\mp1}2\ts\Bigr)^{-1}\,.
\end{equation}
\end{theorem*}

\begin{proof}
Consider the action of the subalgebra $\C_{\ts l}\subset\U(\gl_{n+l})$ on 
$\P\ts(\ts\CC^{\ts m}\ot\CC^{\ts n+l}\ts)\ts$.
Then $\X(\g_n)$ acts on this vector space via 
the homomorphism
$\ga_{\tts l}:\X(\g_n)\to\C_{\ts l}\ts$. Using the decomposition \eqref{pmnl}
we will show that then for any $i\com j=1\lcd n$ the generators
$
S_{ij}^{\ts(1)},S_{ij}^{\ts(2)},\ts\ldots
$
of $\X(\g_n)$ act on this
vector space respectively as the coefficients at 
$u^{-1}\com u^{-2}\com\ts \ldots$ 
of the series \eqref{fhom} multiplied by the series \eqref{fuss}. 

For $i\com j=1\lcd n+l$
the element $F_{ij}\in\U(\g_{n+l})$ acts on
$\P(\ts\CC^{\ts m}\ot\CC^{\ts n+l}\ts)$ as 
$$
\sum_{c=1}^m\,\,
(\ts x_{ci}\,\d_{\ts cj}-
\ts\th_i\,\th_j\,x_{c\ts\bj}\,\d_{\ts c\ts\bi}\ts)\,.
$$
Here we use the standard coordinate functions
$x_{ci}$ on $\CC^{\ts m}\ot\CC^{\ts n+l}$
with $c=1\lcd m$ and $i=1\lcd n+l\ts$.
Then $\d_{ci}$ is the partial derivation on 
$\P\ts(\CC^{\ts m}\ot\CC^{\ts n+l})$ relative to $x_{ci}\ts$.
The functions $x_{ci}$ with $c\leqslant n$ and $c>n$
correspond to the direct summands $\CC^{\ts n}$ and 
$\CC^{\ts l}$ of $\CC^{\ts n+l}$.
Consider the $(n+l)\times(n+l)$ matrix whose $i\com j$ entry is
$$
\de_{ij}\,+\,\Bigl(\ts u-\frac{l\mp1}2\ts\Bigr)^{-1}\,
\sum_{c=1}^m\,\,
(\ts x_{ci}\,\d_{\ts cj}-
\ts\th_i\,\th_j\,x_{c\ts\bj}\,\d_{\ts c\ts\bi}\ts)\,.
$$
Write this matrix and its inverse as the block matrices
$$
\begin{bmatrix}
\,A\,&B\,\\\,C\,&D\,
\end{bmatrix}
\quad\textrm{and}\quad
\begin{bmatrix}\ \At\,&\Bt\ \\ \ \Ct\,&\Dt\ 
\end{bmatrix}
$$
where the blocks $A,B,C,D$ and $\At,\Bt,\Ct,\Dt$ are matrices of sizes
$n\times n$, $n\times l$, $l\times n$, $l\times l$ respectively. 
The action of the algebra $\X(\g_n)$ on the vector space
$\P(\ts\CC^{\ts m}\ot\CC^{\ts n+l}\ts)$ via the homomorphism
$\ga_{\tts l}:\X(\gl_n)\to\C_{\ts l}$ can now be described by assigning
to the series $S_{ij}(u)$ with $i\com j=1\lcd n$ the 
$i\com j$ entry of the matrix $\At^{\,-1}\ts$.

Introduce the $(n+l)\times 2\ts m$ matrix whose $i\com c$ entry 
for $c=-\ts m\lcd-1$ is the operator of multiplication by $x_{ci}$ on 
$\P(\ts\CC^{\ts m}\ot\CC^{\ts n+l}\ts)\ts$. 
For $c=1\lcd m$ let the $i\com c$ entry of this matrix be
the differential operator $-\,\th_i\,\d_{c\ts\bi}\,$. Write this matrix as
$$
\begin{bmatrix}
\,P\,
\\
\,\Pb\,
\end{bmatrix}
$$
where the blocks $P$ and $\Pb$ are matrices of sizes
$n\times2\ts m$ and $l\times2\ts m$ respectively. Further,
introduce the $2\ts m\times(n+l)$ matrix whose $c\com j$ entry 
for $c=-\ts m\lcd-1$ is the operator $\d_{cj}\ts$.
For $c=1\lcd m$ let the $c\com j$ entry of this matrix be
the operator of multiplication by $\th_j\,x_{c\ts\bj}\,$.
Write this matrix as
$$
\begin{bmatrix}
\,Q\,\,\Qb\,\ts
\end{bmatrix}
$$
where $Q$ and $\Qb$ are matrices of sizes
$2\ts m\times n$ and $2\ts m\times l$ respectively. Then
$$
\begin{bmatrix}
\,A\,&B\,\\\,C\,&D\,
\end{bmatrix}
\,=\,
1\,+\,\Bigl(\ts u-\frac{l\mp1}2\ts\Bigr)^{-1}
\begin{bmatrix}\,
\,P\ts Q+m\,&\,P\ts\Qb
\\ 
\Pb\ts Q\,&\,\Pb\ts \Qb+m\,
\end{bmatrix}
$$ 
which can be also written as the matrix
$$
1\,+\,\bigl(\ts u-\frac{l\mp1}2+m\ts\bigr)^{-1}
\begin{bmatrix}\,
\,P\ts Q&\ P\ts\Qb\
\\ 
\,\Pb\ts Q&\ \Pb\ts \Qb\ 
\end{bmatrix}
$$
multiplied by the series $f(u)$ determined by \eqref{fuss}.
By using Lemma \ref{lemma0}, we get
$$
\At^{\,-1}\,=\,A-B\ts D^{-1}\ts C
\,=\,f(u)\,\bigl(\ts
1+\bigl(\ts u-\frac{l\mp1}2+m\ts\bigr)^{-1}\ts P\ts Q
$$
$$
-\,\,
\bigl(\ts u-\frac{l\mp1}2+m\ts\bigr)^{-2}\,
P\ts\Qb\,\ts\bigl(\ts1+
\bigl(\ts u-\frac{l\mp1}2+m\ts\bigr)^{-1}\ts\Pb\ts\Qb\,\bigr)^{-1}\ts\Pb\ts Q
\,\bigr)
\vspace{2pt}
$$
\begin{equation}
\label{pulq}
=\,f(u)\,\bigl(\ts1+P\ts
\bigl(\ts u-\frac{l\mp1}2+m+\Qb\ts\Pb\,\bigr)^{-1}
\ts Q\,\bigr)\,.
\vspace{4pt}
\end{equation}

Consider the $2m\times2m$ matrix $\Qb\ts\Pb$ appearing in the last line.
For any indices
$a\com b=-\ts m\lcd-1\com1\lcd m$ the $a\com b$ entry of this matrix 
is the differential operator 
$$
\de_{ab}\,\frac{l}2\ts-\,\bar\zeta_{\ts l}\ts(F_{ab})
$$ 
where $\bar\zeta_{\ts l}:\U(\ts\f_m)\to\PD\ts(\CC^{\ts m}\ot\CC^{\ts n+l})$
is the homomorphism corresponding to the action of the 
Lie algebra $\f_m$ on $\P\ts(\CC^{\ts m}\ot\CC^{\ts n+l})$
via the tensor factor $\P\ts(\CC^{\ts m}\ot\CC^{\ts l})$ in \eqref{pmnl},
similar to the homomorphism \eqref{gan}.
Namely for $a\com b=1\lcd m$ we have
$$
\bar\zeta_{\ts l}\ts(F_{ab})\,=\,
\de_{ab}\,\frac{l}2\ +
\sum_{k=n+1}^{n+l}\,x_{ak}\,\d_{\ts bk}\,,
$$
$$
\bar\zeta_{\ts l}\ts(F_{a,-b})\,=\,-\,\sum_{k=n+1}^{n+l}\,
\th_k\,x_{a\bk}\,x_{\ts bk}\,,
\qquad
\bar\zeta_{\ts l}\ts(F_{-a,b})\,=\sum_{k=n+1}^{n+l}\,
\th_k\,\d_{ak}\,\d_{\ts b\bk}\,.
$$
Hence any entry of the $2m\times2m$ matrix 
$$
\bigl(\ts u-\frac{l\mp1}2+m+\Qb\ts\Pb\,\bigr)^{-1}
$$
can is obtained by applying the homomorphism $\bar\zeta_{\ts l}$
to the respective entry
of the matrix $F\ts(u\pm{\textstyle\frac12}+m)\ts$;
the latter entries are series in $u^{-1}$
with coefficients in $\U(\ts\f_m)\ts$.
Comparing \eqref{fhom} with the $i\com j$ entry of 
$n\times n$ matrix \eqref{pulq} completes the proof.
\end{proof}

Set $\C_0=\U(\g_n)$ and $\ga_{\ts0}=\pi_n\ts$.
Then Theorem \ref{5.1} remains
valid for $l=0\ts$, assuming that $\g_{\ts 0}=\{0\}\ts$.
Further, the homomorphism $\zeta_{\ts l}$ is injective for
any $l\ge m\ts$. Thus when $l\ge m\ts$, 
our proof of Theorem \ref{5.1} also implies Proposition \ref{xb}.

Let $\la$ and $\mu$ be the parameters of any
irreducible polynomial representations of $G_{n+l}$ and $G_l$
respectively. The vector space \eqref{hll} is not zero if and only if 
\begin{equation}
\label{llcon}
\la_k\ge\mu_k
\quad\textrm{and}\quad
\las_k-\mus_k\le n
\quad\textrm{for every}\quad
k=1\com2\com\ts\ldots\ts;
\end{equation}
see \cite[Section 1.3]{N2}. 
Suppose that $\las_1\com\mus_1\le m\ts$.
Then we can identify
the vector spaces \eqref{hll} and \eqref{hlls}.
Then the algebra $\C_l$ acts on \eqref{hlls} irreducibly, if $G_n=Sp_n\ts$.
If $G_n=O_n\ts$, then \eqref{hlls} is irreducible
under the joint action of the algebra $\C_l$ and the group $O_n\ts$.
In both cases, the $G_{n+l}\ts$-invariant elements of $\U(\gl_{n+l})$ 
act on \eqref{hlls} via multiplication by scalars.
Then Theorem \ref{5.1} has a corollary,
which refers to the 
action of $\X(\g_n)$ on the vector space \eqref{hlls}
inherited from the bimodule\/ $\F_m(L_{\ts\mut})\ts$.

\begin{corollary*}
\label{xirred}
The algebra\/ $\X(\g_n)$ acts on space \eqref{hlls} irreducibly,
if\/ $G_n=Sp_n\ts$. 
If $G_n=O_n\ts$, the space \eqref{hlls} is irreducible
under the joint action of\/ $\X(\g_n)$ and\/ $O_n\ts$.
\end{corollary*}

Let $V$ and $V'$ be irreducible highest weight modules
of the Lie algebra $\f_m\ts$,
corresponding to unitary
representations of the real metaplectic group
$Mp_{\ts 2m}$ if $\f_m=\sp_{2m}\ts$, 
or of the real group $SO_{2m}^\ast$ if $\f_m=\so_{2m}\ts$.
Consider the vector space \eqref{vvp}.
According to the results of \cite{EP} or \cite{EHW}
when $\f_m=\sp_{2m}$ or $\f_m=\so_{2m}$ respectively, the
$\f_m\ts$-module $V$ is equivalent to $L_{\ts\mut}$ for some
non-negative integer $l$ and the label $\mu$
of some irreducible polynomial representation of the group $G_l$
with $\mus_1\le m\ts$. If 
the vector space \eqref{vvp} is non zero, then
$V^{\ts\prime}$ has to be equivalent to $L_{\ts\lat}$ for the label $\la$
of some irreducible polynomial representation of $G_{n+l}$ with
$\las_1\le m\ts$. Thus any non-zero vector space \eqref{vvp}
has to be of the form \eqref{hlls}.


\subsection*{Acknowledgments}
We are grateful to V.\,Tarasov and V.\,Tolstoy for friendly discussions
which led us to the results presented in this article.
The first author was supported by the RFBR grant 05-01-01086, 
the grant for Support of Scientific Schools 8065-2006-2,  
the ANR grant 05-BLAN-0029-01,
and by the Atomic Energy Agency of the Russian Federation.
The second author was supported by the EPSRC grant C511166,
and the EC grant MRTN-CT2003-505078. 
This work began when both authors visited MPIM. 
Part of this work was also done
when the first author visited IHES.
We are grateful to the staff of both institutes for
their kind help and generous hospitality.



\end{document}